\documentclass[12pt]{amsart}
\usepackage{amsmath,amssymb,amsfonts,latexsym}
\usepackage{mathrsfs}%pour les caracteres en calligraphie
\usepackage[dvips]{geometry}
\usepackage{array}
\usepackage{epsf}
\usepackage{epsfig}
\newtheorem*{lemma}{Lemma}
\newtheorem*{prop}{Proposition}
\newtheorem*{thm}{Theorem}
\newtheorem*{cor}{Corollary}
\newcommand{\iso}{\overset{\sim}{\rightarrow}}

\newcommand{\twoheaddownarrow}{\overset{\sim}{\twoheaddownarrow}}

\newcommand{\ad}{\operatorname{ad}}

\newcommand{\nc}{\newcommand}
\nc{\Ker}{\operatorname{Ker}} \nc{\rker}{\operatorname{rKer}}
\nc{\im}{\operatorname{Im}}
\nc{\stab}{\operatorname {Stab}}
\nc{\ann}{\operatorname {Ann}}
\nc{\Id}{\operatorname {Id}}
\nc{\Prim}{\operatorname {Prim}}
\nc{\Real}{\operatorname {Re}}
\nc{\Ext}{\operatorname {Ext}}
\nc{\rad}{\operatorname {rad}}

\begin {document}

\title[Trails and $S$-graphs]{Trails, $S$-graphs and identities in Demazure modules}
\author[Anthony Joseph]{Anthony Joseph}

\date{\today}
\maketitle

\vspace{-.9cm}\begin{center}
Donald Frey Professional Chair\\
Department of Mathematics\\
The Weizmann Institute of Science\\
Rehovot, 76100, Israel\\
anthony.joseph@weizmann.ac.il
\end{center}\

%\centerline {SHMUEL ZELIKSON}
%\begin{center}
%D\'epartement de Math\'ematiques\\
%Universit\'e de Caen, B.P. 5186\\
%14032 Caen Cedex, France\\
%Shmuel.Zelikson@unicaen.fr \end{center}
%
%\

\medskip


 \

 \textbf{Abstract}.  Let $\mathfrak g$ be a Kac-Moody algebra.  The Kashiwara $B(\infty)$ crystal parameterizes a basis for the Verma module of highest weight zero.  For every sequence $J$ of reduced decompositions of Weyl group elements it admits a realization $B_J(\infty)$ as a subset of a crystal $B_J$ which as a set identifies with the free additive semi-group $\mathbb N^{|J|}$.   A natural question is whether $B_J(\infty)$ is a polyhedral subset and if so to determine a (minimal) set of linear inequalities which define $B_J(\infty)$ as a subset of $B_J$.

 In earlier work this question led to the notion of an $S$-set associated with a given simple root.  In the present work the notion of a giant $S$-set is studied. Essentially this is for every simple root, a disjoint union of $S$-sets associated to that simple root.  Its elements are to give the linear inequalities which define $B_J(\infty)$.

 A trail is a certain sequence of vectors specified by the choice of $J$ in a fundamental module for Langlands dual of $\mathfrak g$.  Every such sequence gives a linear function on $B_J$ whose elements are to give the linear inequalities which define $B_J(\infty)$.

 The present work inter-relates these two approaches for describing $B_J(\infty)$.  Under the hypothesis that there are no ``false trails'', the set of all trails is determined.  Then the resulting functions are shown to describe the ``envelope'' of a giant $S$-set.  This in turn leads to the description of $B_J(\infty)$ as a polyhedral subset of $B_J$.

\section{Introduction}\label{1}

\subsection{}\label{1.1}

Let $\mathfrak g$ be a Kac-Moody algebra (\textit{not} necessarily symmetrizable). Recall that the construction of $\mathfrak g$ depends on a choice of Cartan subalgebra $\mathfrak h$ and a set $\pi:=\{\alpha \in \mathfrak h^*\}$ (resp. $\pi^\vee:=\{\alpha^\vee \in \mathfrak h\}$) of simple roots (resp. coroots). Let $W$ denote the corresponding Weyl group and $P$ (resp. $P^+$) the corresponding set of integral (resp. integral and dominant) weights. Let $\mathfrak b$ be the Borel subalgebra of $\mathfrak g$ defined by $\pi$ and $\mathfrak n$ its nilradical.

Assume for the moment that $\mathfrak g$ is symmetrizable.  Then after Kashiwara \cite {Ka1} the dual Verma module $\delta M(0)$ for the pair $(\mathfrak g, \mathfrak b)$ of highest weight zero admits a crystal basis.  We call the corresponding crystal, the Kashiwara $B(\infty)$ crystal. It has a combinatorial structure of great intricacy which one would like to understand. This is partly because it can provide useful information on the nature of tensor product decomposition of integrable highest weight modules and on their Demazure submodules and partly because we believe it leads to new combinatorial constructions, specifically that of giant $S$-sets.

\subsection{}\label{1.2}

Again for $\mathfrak g$ symmetrizable, Lusztig \cite {Lu1} showed that the Verma module $M(0)$ admits a canonical basis which coincides with Kashiwara's globalisation of the $B(\infty)$ crystal basis.  In this Lusztig further discussed how this basis behaves under the action of the simple root vectors $e_\alpha: \alpha \in \pi$. In particular for $\mathfrak g$ simply-laced, the resulting coefficients are all non-negative integers.

%One may shift the highest weight of $\delta M(0)$ to an arbitrary weight $\lambda \in P$ without changing its structure as a $[\mathfrak n$ module.

Assume $\lambda \in P^+$ and let $V(\lambda)$ be the integrable module with highest weight $\lambda$. let $k_{-\lambda}$ denote the one-dimensional $\mathfrak b$ module of weight $-\lambda$. As noted in \cite [Lemma 2.4]{J1} there is a unique up to scalars $\mathfrak b$ module embedding of $V(\lambda)\otimes k_{-\lambda}$ into $\delta M(0)$.  It is then a consequence of Kashiwara-Lusztig theory (cf \cite [Thm. 6.2.19]{J2}) that the dual canonical basis of $\delta M(0)$ restricts to a basis of $V(\lambda)$, called its dual canonical basis.  Here one may remark that the extremal vectors $v_{w\lambda}$ of weight $w\lambda: w \in W$, up to an appropriate scalar, belong to this basis.

\subsection{}\label{1.3}

Fix an index set $I$ for $\pi$. Given $i \in I$ let $s_i$ denote the corresponding reflection and $\varpi_i$ the corresponding fundamental weight.  Fix a possibly infinite sequence $J= \ldots,i_j,i_{j-1},\ldots,i_1$ of elements of $I$ so that $w_j:=s_{i_j}s_{i_{j-1}}\cdots s_{i_1}$ is a reduced decomposition for all $j \in J$.  A construction of Kashiwara (cite {Ka2}, \cite [2.4]{J3}) gives a crystal $B_J$ which identifies with $\mathbb N^{|J|}$ as a set.  The element $b_\infty$ in which all entries are zero, generates a subcrystal $B_J(\infty)$.  A remarkable result of Kashiwara is that, as a crystal $B_J(\infty)$, is independent of $J$. This follows from Kashiwara's embedding theorem which we remark is extended in \cite [2.5.7]{J3} to the case when $\mathfrak g$ is not necessarily symmetizable.  Here the way in which the embedding theorem is applied is discussed in \cite [2.5.8]{J3}.  In this there is a minor subtlety concerning the number of appearances of a given element of $I$ in $J$, but it need not detain us here.  Suffice to say that for $B(\infty)$ to be independent of $J$ as a crystal, it is enough that $s \in I$ appears in the set $\{i_j:j \in J\}$ sufficiently many times.  In any case we shall not presently be concerned with the independence of $B_J(\infty)$ on $J$ as a crystal.

%Indeed we may regard our computations as pertaining to the Demazure ``crystals'' $B_{w_j}(\infty)$ as described for example in \cite [3.4.3]{J3} and eventually take $B(\infty)$ to be their limit.  %This yields the crystal $B(\infty)$ referred to in \ref {1.1}.

\subsection{}\label{1.4}

As a subset of $B_J$ (and hence of $\mathbb N^{|J|}$), the crystal $B_J(\infty)$ is \textit{dependent} on $J$ and can be \textit{extremely complicated} to describe.  Yet for $\mathfrak g$ of type $A$, that is to say for $\mathfrak g=\mathfrak {sl}(n)$, Gleizer and Postnikov \cite {GP} gave a beautiful description of $B_J(\infty)$ in a purely combinatorial manner using wiring diagrams.  It is notable that they viewed their result as an extension of the classical Littlewood-Richardson rule for tensor product decomposition.  Moreover they showed that $B_J(\infty)$ is always a polyhedral subset of $B_J$ (in type $A$).  This latter result was extended to the case when $W$ is finite, by Berenstein and Zelevinsky \cite {BZ}.  However this latter description was less precise in that it was given in terms of $\textbf{i}$-trails (which we call simply, trails) which are not combinatorially defined and are basically unknown quantities.

\subsection{}\label{1.5}

In \cite {J4} we initiated a new approach to the description of $B_J(\infty)$.  This was based on the duality operation on $B(\infty)$, due to Kashiwara for $\mathfrak g$ symmetrizable and extended to all $\mathfrak g$ Kac-Moody in \cite [Sect. 2] {J3}.  For this one must construct a set $Z_t$ of ``dual Kashiwara functions'' whose maximum must be ``almost'' invariant under the action of the Kashiwara operators. In particular invariance under the Kashiwara operators assigned to a fixed $s \in I$, resulted in the notion of an $S$-set of linear functions on $B_J$ (identified with $\mathbb N^{|J|}$). Then $Z_t$ is to be a disjoint union of $S$-sets for every $s \in I$ with one exceptional element $z^1_t$ called the initial driving function associated to $t$.  (For a more precise statement see \ref {8.1}.)

We call such a set $Z_t$, if it exists - a giant $S$-set.  Its structure must be particularly intricate and a challenge to describe.

Given a giant $S$-set $Z_t$, its maximum on $B_J$ restricts to the function $\varepsilon^\star_t$ on $B_J(\infty)$ described in \cite {Ka2}. We call $\varepsilon^\star_t:t \in I$ the dual Kashiwara parameter associated with $t \in I$.  The set of dual Kashiwara parameters may be used to determine $B_J(\infty)$ as a precise polyhedral subset of $B_J$, for any choice of $J$.  Here we remark that the $\varepsilon^\star_t:t \in I$ are determined intrinsically through the crystal structure of $B_J(\infty)$, which is independent of $J$; but viewed as functions on $B_J(\infty)\subset B_J$, are dependent on $J$.

%We call such a set $Z_t$ if it exists - a giant $S$-set.

To us the interest in the present work lies in the description of $S$-sets, giant $S$-sets, trails and their properties and inter-relations rather than in just describing $B_J(\infty)$.

  Thus in \cite {J4}, $S$-sets were described as the vertices of certain $S$-graphs.   A canonical family of $S$-graphs was described \cite {J4}, \cite {JL} and it was shown that the functions they define exhibit some remarkable convexity properties \cite {J5}.   In the present work we concentrate on the much more subtle giant $S$-set.  It may again be viewed as the set of vertices of a graph called a giant $S$-graph whose vertices are obtained from those of its $S$-subgraphs. At first it was completely unclear for the moment why our canonical $S$-graphs are the correct ones to choose in describing a giant $S$-graph.  The answer can now be seen to be given by the results of Section \ref {7}, which relates $S$-sets to trails under the hypothesis that there are no ``false trails''.

\subsection{}\label{1.6}

One may read off from \cite {GP} (resp. \cite [Thm. 3.9] {BZ}) the dual Kashiwara functions in type $A$ (resp. in finite type).  In the latter case they are defined in terms of unknown trails.  In this we remark that a trail in the lowest weight module $V(-\varpi_t): t \in I$ is a sequence of \textit{non-zero} vectors $v_j \in V(-\varpi_t):j \in J$ with $v_1=v_{-s_t\varpi_t}$ and then defined inductively through $v_{j+1}=e^{n_j}_{\alpha_{i_j}}v_j$ for some $n_j \in \mathbb N$ satisfying the following boundary condition. Namely there exists $k \in \mathbb N^+$ such that $v_{j+1}$ is proportional to the extremal vector $v_{-w_j\varpi_t}$, for all $j\geq k$.  In particular these extremal vectors will have several different presentations.  This boundary condition implies that the linear function on $B_J$ defined by a given trail is locally finite.  It can happen that $v_j$ is always an extremal vector for all $j \in J$ (though not necessarily proportional to $v_{-w_{j-1}\varpi_t}$).  Such trails are easy to describe and of course exhaust all trails if $V(-\varpi_t)$ is minuscule.  Otherwise the boundary condition makes it extremely difficult to write down a trail.

One should of course be aware that the set of trails in $V(-\varpi_t)$ depends on the choice of $J$. Here we shall always take $J$ to be fixed and not mention this dependence further.

\subsection{}\label{1.7}

The present work is a by-product of attempting to reconcile our work on $S$-graphs with the description in \cite {BZ} of dual Kashiwara functions.  In this we conjecture the absence of ``false trails''.  This involves a natural suggestion about relations between monomials expressions in simple root vectors for Demazure modules extending the relations in $U(\mathfrak n)$ itself. Here we remark that the latter are known at least when $\mathfrak g$ is symmetizable to be given by the Chevalley-Serre relations, a result \cite {GK}  considered deep, though ultimately not too difficult.

Given the truth of this conjecture we are also able to give a simple algorithm to compute the proportionality factors between the extremal vectors given by the various trails.  These proportionality factors are binomial coefficients, hence positive as would be expected from the positivity result of Lusztig in the simply-laced case mentioned in \ref {1.2}.

The main result of this paper is that the absence of ``false trails'' in $V(-\varpi_t)$ is equivalent to the existence of the giant $S$-set envelope $\mathscr K_t$ associated to $t\in I$. In principle the desired giant set $S$-set $Z_t$ can be recovered as the set of extremal points of $\mathscr K_t$, whilst the latter should be the $\mathbb Z$-convex hull of $Z_t$. For the moment these refinements are open questions and in any case not of any importance to the computation of the dual Kashiwara parameter $\varepsilon^\star_t$ - see Theorem \ref {8.7}.

In general one may compute $Z_t$ by induction on $j \in J$, using our knowledge of $S$-sets; but in general it is not assured that this procedure can result in a giant $S$-set unless ``false trails'' are absent.  Nevertheless when $V(-\varpi_t)$ is minuscule one can easily show that this inductive procedure does lead to a giant $S$-set, though this computation will not be presented here.  Again the complexity of $Z_t$ depends on the choice of $J$ and it is often possible to find choices for which the $Z_t:t \in I$ are rather simple. We find these cases rather boring and are more concerned with ``nasty'' choices of $J$ for which $Z_t$ is very complicated.

 It is always true that  $Z_t=\{z_t^1\}$, if $i_1=t$, and indeed this combined with the behaviour of $B_J(\infty)$ under changes of $J$ led to the above mentioned result in \cite [Thm. 3.9] {BZ}).  This procedure is restricted to finite type since otherwise it is impossible to alter a given $J$ so that $j_1$ becomes a chosen element of $I$.  Moreover handling the resulting computations based on tropical fields and cluster algebras is by no means easy.

 It should be noted that using trails in $V(-\varpi_t)$ leads to a giant $S$-set not for $\mathfrak g$, but for the Langlands dual of $\mathfrak g$.  This is no surprise considering that \cite {BZ} considered fundamental modules for the Langlands dual in the first place.  To obtain a giant $S$-set for $\mathfrak g$ one has only (to remember!) to interchange roots and co-roots in the expressions for the functions they define.

%Again the truth of our conjecture implies that the set of all trails in $V(-\varpi_t)$ describe a set of dual Kashiwara functions of type $t\in I$.

%Their maximum gives the function on $\varepsilon^\star_t$ on $B_J(\infty)$ described in \cite {Ka2}.  Then\ the $\varepsilon^\star_t:t \in I$ may be used to determine $B_J(\infty)$ as a precise polyhedral subset of $B_J$, for any choice of $J$.  Here we remark that the $\varepsilon^\star_t$ are determined intrinsically through the crystal structure of $B_J(\infty)$, which is independent of $J$; but viewed as functions on $B_J(\infty)\subset B_J$ with $B_J$ identified with $\mathbb N^{|J|}$, are dependent on $J$.

\subsection{}\label{1.8} In the absence of false trails, $\mathscr K_t$ indexes all trails in $V(-\varpi_t)$ which can then be described purely combinatorially.   Then by the results of \cite {JZ} the maximum of the resulting functions on $B_J(\infty)$ gives dual Kashiwara parameter $\varepsilon^\star_t$.  Up to the absence of false trails, this extends the result of Berenstein and Zelevinsky from the finite to the general case and since trails are determined is a much more explicit result.    In the sense that we only need to compute maxima, the set of all trails is highly redundant and the description of these maxima by a giant $S$-set much more economic. At first sight this redundancy was a great surprise.  It is now understood through the present work.

 However perhaps the most fascinating part of this work is that $S$-graphs which were invented to describe invariance under the Kashiwara operators made a second appearance as providing extremal points of certain convex sets \cite {J5}, which is not so surprising, but now also make a third appearance by virtue of the fact that the integer points of these convex sets are determined by  relations in Demazure submodules of fundamental modules.

 \textbf{Acknowledgements.}  I would like to thank S. Zelikson for explaining to me the results in \cite {BZ} and his understanding of the relation between trails and functions on $B_J(\infty)$, without which this paper would not have been brought into existence.  He also provided me with extremely useful computer outputs of the description of these functions in types $F_4$ and $D_5$ for some (nasty) special choices of $J$.  These are far to long to reproduce here and it seems practically impossible to illustrate diagrammatically the giant $S$-graphs in these cases.

\section{Trails}\label{2}

\subsection{}\label{2.1}

The notion of an $\textbf{i}$ trail, which we call simply, a trail, was introduced by Berenstein and Zelevinsky \cite {BZ}, in the case when $W$ is finite to describe $B_J(\infty)$ as a polyhedral subset of $B_J$.  In this they were also able to describe the dual Kashiwara parameters \cite [Thm. 3.9]{BZ}.   One may remark that in \cite {BZ} trails are defined in fundamental modules for the Langlands dual of $\mathfrak g$.  Here we prefer to stay with $\mathfrak g$ itself and eventually just interchange roots and co-roots in the resulting functions (see \ref {2.4}).

 When $W$ is finite, a trail is a sequence of vectors in a fundamental module from the almost lowest weight vector of weight $-s_t\varpi_t$ to its highest weight vector.  In the general case we replace the latter end-point by a boundary condition.  Although this is an obvious thing to do, we stress that this boundary condition is rather subtle and leads to a surprising and valuable conclusion (Lemma \ref {4.3}).

\subsection{Definitions}\label{2.2}

Define $I,J$ as in \ref {1.3} and recall the notation introduced there. As a set $B_J$ identifies with $\mathbb N^{|J|}$, that is to say with $|J|$-tuples of natural numbers in which almost all entries are zero. If $W$ is finite, let $w_0$ denote its unique longest element and identify $J$ with $\{1,2,\ldots,\ell(w_0)\}$, where $\ell(\cdot)$ denotes reduced length.
 %It is convenient to set $\hat{J}=\{1,2,\ldots,\ell(w_o)+1\}$ (resp. $J=\hat{J}$) if $W$ is finite (resp. infinite).
 Otherwise $J$ is a countable set and may be identified with $\mathbb N^+$.

Fix $t \in I$.

Let $V(-\varpi_t)$ denote the lowest weight module of lowest weight $-\varpi_t$ and $\Omega(V(-\varpi_t))$ its set of weights.

It is convenient to relabel $J$ in the following fashion.

For all $s \in I, k \in \mathbb N^+$ let $(s,k)$ denote the element of $J$ for which $i_j=s$ for the $k^{th}$ time counting from the right.  Observe that the natural order on $J$ defines a total order on this set of pairs.  We may also use these two parameterizations simultaneously and in particular $(s,k)+1$ means $j+1$, when $(s,k)=j$. Note also that if $u=(s,k)$, then $i_u=s$ by definition.

A trail $K$ associated to $t$ is a sequence of non-zero vectors $v_{\gamma^K_j} \in V(-\varpi_t):j \in \hat{J}$ of weight $\gamma^K_j$ satisfying the following rules

\

$(T)$. \textbf{Trails}.  For all $j \in J$, there exists $n_j \in \mathbb N$ such that  $e_{i_j}^{n_j}v_{\gamma^K_j}=v_{\gamma^K_{j+1}}$.

\

(B).  \textbf{The Boundary Conditions.}

\

(i). $\gamma^K_1=-s_t\varpi_t$.

(ii). For every trail $K$ there exists $\varphi(K) \in \mathbb N^+$ such that $\gamma^K_{j+1}=-w_{j}\varpi_t$, for all $j\in J|j\geq  \varphi(K)$.

\

We shall say that a trail $K$ trivializes at $w_j:j \in J$, or simply at $j \in J$, if $j \geq \varphi(K)$.  From then on the trail is just the appropriate sequence of extremal vectors.

From this definition it follows that the subset of $\mathscr K^{BZ}_t$ of trails which trivialize at $j \in J$ is increasing in $j$.

When $W$ is finite we take $\varphi(K) =\ell(w_0)$, for all $K \in \mathscr K$.  Notice that this means that $\gamma^K_{\ell(w_0)+1}=-w_0\varpi_t$.

\

The initial driving trail $K^1_t$ associated to $t$ is defined as follows.

Set $u=(t,1)$, so then $i_u=t$.  Set

$$\gamma^{K_t^1}_j:=\left\{
                       \begin{array}{ll}
                         -s_t\varpi_t, & \hbox {if $j\leq u+1$,} \\
                         -s_{i_{j-1}}\ldots s_{i_{u+1}}s_t\varpi_t, & \hbox {if $j>u+1$.}
                       \end{array}
                     \right. \eqno {(1)}$$

The driving trail $K_t^1$ consists just of extremal vectors; but it only trivializes at
$(t,1)$, because $\gamma^{K_t^1}_{j+1}=-w_j\varpi_t$, if and only if $j \geq u=(t,1)$.

\

$(P)$. \textbf{The Positivity Condition.}

\

 $\gamma^K_j \in \gamma_j^{K_t^1} +\mathbb N\pi$, for all $j \in J$.

\

 We use $\mathscr K_t^{BZ}$ to denote the set of all (Berenstein-Zelevinsky) trails associated to $t$, that is to say all trails $K$  satisfying $(T),(B), (P)$  above.

 We write $v_{\gamma^K_j}$ simply as $v^K_j$. We remark that the sequence of weights $\gamma^K_j: j \in J$ \textit{uniquely determines} the trail $K$ \textit{if} it exists.

Notice our insistence on $v_{\gamma^K_j}$ being non-zero. When we drop this hypothesis we shall call the resulting object a potential trail.  The hard part of the present work is to determine when a potential trail is a trail.

\begin {lemma}  Let $j \in J$ be minimal such that there exists a trail $K\in \mathscr K^{BZ}_t$ which trivializes at $w_j$. Then $j=(t,1)$ and $K=K_t^1$.
\end {lemma}

\begin {proof}
 Since $K_t^1$ trivializes at $w_{(t,1)}$, it follows from the choice of $j \in J$ and  $K \in \mathscr K^{BZ}_t$, that $\gamma^K_{j+1}=-w_j\varpi_t$, if $j \geq (t,1)$. In particular $\gamma^K_{(t,1)+1}=-w_{(t,1)}\varpi_t=-s_t\varpi_t$. %Now suppose $K$ trivializes at some $j<(t,1)$.  We can write $j=(s,k)$, for some $s \in I\setminus \{t\}, k \in \mathbb N^+$.
 Since $-s_t\varpi_t$ is the unique minimal weight of $\Omega(V(-\varpi_t))\setminus\{-\varpi_t\}$, it follows from $(P)$ that $\gamma_j^K=-s_t\varpi_t$, for all $j \leq (t,1)$.  Hence $K=K_t^1$.
\end {proof}

\subsection{Functions}\label{2.3}

To a trail $K$ we may associate a linear function $z^K$ on $B_J$ as follows.  Let $m_j$ denote the $j^{th}$ co-ordinate function on $B_J$.  Set $\delta_j^K=\frac{1}{2}(\gamma_j^K+\gamma_{j+1}^K)$ and

$$z^K:= \sum_{j \in J}\alpha^\vee_{i_j}(\delta_j^K)m_j.  \eqno {(2)}$$

For example if $K=K_t^1$, then $z^K$ is, up to interchanging $\mathfrak g$ by its Langlands dual, just the initial driving function $z_t^1$ associated to $t$ introduced in \cite {JZ}.

By the second boundary condition, $z^K$ is locally finite on $B_J$, for all $K \in \mathscr K^{BZ}$.

When we write $j \in J$ as $(s,k) \in I\times \mathbb N^+$ as prescribed in \ref {2.2}, we set $m_j=m_s^k$.

Again we may define the Kashiwara functions $r_s^k: s \in I, k \in \mathbb N^+$ which in this paper we shall always take to be for the Langlands dual of $\mathfrak g$, that is to say by interchanging roots and co-roots in \cite [2.3.2]{J3}.  Explicitly one has
$$r_s^k=m_s^k + \sum_{j>(s,k)}\alpha^\vee_{i_j}(\alpha_s)m_j, \forall (s,k) \in I \times \mathbb N^+.$$

In particular these functions are linear on $B_J$ (though \textit{not} locally finite).

By the definition of $B_J$, given for example in \cite [2.4.2]{J3}, there exists for all $b \in B_J$ an element $j_b \in \mathbb N^+$ such that $m_j(b)=0$, for all $j >j_b$.  Thus $r_s^k(b)$ is a finite sum (of integers) for all $b \in B_J$.

Through the first boundary condition we can reconstruct the set of weights $\gamma_j^K:j \in J$ from $z^K$ and hence obtain a potential trail $K$.  Of course it is completely unclear when $K$ is a trail, that is when $K \in \mathscr K^{BZ}_t$.

Observe that if a trail $K$ trivializes at $n \in \mathbb N^+$, then the coefficient of $m_j$ in $z^K$ is zero for $j>n$, that is to say the sum in $(2)$ stops at $j=n$.  In the finite case a trail trail trivializes at $\ell(w_0)$ and then stops. Corresponding the sum in $(2)$ stops at $j=\ell(w_0)$.

For example the driving trail $K_t^1$ trivializes at $(t,1)$ and the sum in $(2)$ is up to $j=(t,1)$.

\begin {lemma}  For all $K \in \mathscr K^{BZ}_t, j \in J$ one has
$$s_{i_j}\gamma_{j+1}^K=\gamma_j^K-\alpha_{i_j}^\vee(\delta^K_j)\alpha_{i_j}.$$
In particular the right hand side is a weight of $V(-\varpi_t)$.

\end {lemma}

\begin {proof}  Set $\alpha_{i_j}=\alpha$ and $s=s_{i_j}$ and omit the superscript $K$.  There exists $q_{j+1}\in \mathbb Q$ and $ \varpi \in P(\pi)$ such that $\gamma_{j+1}=q_{j+1}\alpha+\varpi$ and $\alpha^\vee(\varpi)=0$. By $(T)$ we obtain $\gamma_{j}=q_j\alpha+\varpi$, with $q_{j+1}+q_j=\alpha^\vee(\delta_j)$.  Yet $s\gamma_{j+1}=\gamma_j-(q_{j+1}+q_j)\alpha$, so the assertion obtains.
\end {proof}

\subsection{Faces}\label{2.4}

 Given $(s,k) \in I\times \mathbb N^+$, with $k>1$, we define the closed face $F_s^k$ of type $s$ as follows.  Set $u=(s,k), v=(s,k-1)$ viewed as elements of $J$ (as above). Recall that $i_u=i_v=s$. Set
   $$\gamma_j^{F_s^k}:=\left\{
                       \begin{array}{ll}
                         \alpha_s, & \hbox {if $v<j\leq u$,} \\
                         0, & \hbox {otherwise.}
                       \end{array}
                     \right. \eqno {(3)}$$

 Here we regard $F_t^1:=K_t^1$ as the unique open face associated to $t$.

 This terminology and what follows is motivated by wiring diagrams in type $A$.  In this, faces are bounded by paths and adjoining a face corresponds to altering the original path to go around the added face.  Gleizer and Postnikov \cite {GP} found a way to exclude certain paths and then the set of all allowed paths going from one end point of the unique open face for $t$ to the other end point, identifies with $Z_t$.  Wiring diagrams can be defined for all $\mathfrak g$, though are easier to visualize if the Dynkin diagram is linear.  However it is not clear what should be the allowed paths and some paths need to go around a face more than once.

 Given a trail $K$, a basic question is the following.  When can the face $F_s^k$ be adjoined to $K$? More precisely given a trail $K$ when is there a trail $K+F_k^s$ determined by the condition
 $$\gamma^{K+F_s^k}_j=\gamma_j^{F_s^k}+\gamma^K_j, \forall j \in J. \eqno {(4)}$$

It is clear that we may define a function $z^{F_s^k}$ on $B^J$ by replacing $K$ by $F_s^k$ in $(2)$ above.  It satisfies $(4)$, with $\gamma_j$ is replaced by $z$.  It is instructive to compute $z^{F_s^k}$.  Recall that $k >1$ and set $u=(s,k), v=(s,k-1)$. One checks that
$$ z^{F_s^k}=  m_u + \sum_{j=v+1}^{u-1}\alpha^\vee_{i_j}(\alpha_s) m_j +  m_v=r_s^{k-1}-r_s^k, \quad \forall (s,k-1) \in I \times \mathbb N^+.\eqno {(5)}$$

One may define $r_s^0=-\sum_{j \in J}\alpha_{i_j}^\vee(\alpha_s)m_j$ and check that $$z_s^1=r_s^0-r_s^1=\sum_{j=1}^{(s,1)}(\alpha_{i_j}^\vee(\alpha_s)-\delta_{i_j,s})m_j=
m_s^1+\sum_{j=1}^{(s,1)-1}\alpha^\vee_{i_j}(\alpha_s)m_j,$$
where $\delta$ is the Kronecker delta.  Thus the successive differences $r_s^{k-1}-r_s^k:(s,k) \in I\times \mathbb N^+$ of the Kashiwara functions are locally finite and are related to the co-ordinate functions of $B_J$ by a triangular matrix with ones on the diagonal.

For all $t \in I$, set
$$X_t:=z_t^1+\sum_{(s,k) \in I \times \mathbb N^+}\mathbb N(r_s^k-r_s^{k+1}).$$
In this sums are viewed as being finite.  Observe that there is a partial order $\leq$ on $X_t$ defined by $z\leq z'$ if $z'-z \in \sum_{(s,k) \in I \times \mathbb N^+}\mathbb N(r_s^k-r_s^{k+1})$.

%Interchanging roots and co-roots, the right hand side of $(5)$ becomes the difference $r_s^{k-1}-r_s^k$ of successive Kashiwara functions of type $s$.

The basic idea formulated in \cite {JZ} is that the giant $S$-set $Z_t$ can be constructed as a subset of $X_t$ inductively, starting from the initial driving function $z_t^1$ associated to $t$, by successively adding non-negative multiples of the $r_s^k-r_s^{k+1}: k \in \mathbb N^+$ given by an appropriate $S$-set of type $s$ specified by a driving function of type $s$ lying in the subset of $Z_t$ previously obtained.  In this we remark that $z_t^1$ is not a driving function of type $t$.  It is a driving function of type $s$ for all $s \in I\setminus \{t\}$; leading to $S$-sets strictly larger than $\{z_t^1\}$ if and only $\alpha^\vee_s(\alpha_t) \neq 0$ and $(s,1) <(t,1)$.  If this condition is not satisfied for all $s \in I\setminus\{t\}$, then $Z_t$ is reduced to $z_t^1$, which in turn equals $m_t^1$ in this case.
 %$r_s^k-r_s^{k+1}: k \in \mathbb N^+$ according to the rules one can obtain the set of dual Kashiwara functions of type $t$ by adding suitable non-negative multiples of successive differences of the dual Kashiwara functions for $\mathfrak g$ to the initial driving function $z_t^1$ for $\mathfrak g$ and to subsequent driving functions of type $s$, with the particular choice of sums to be given by the appropriate $S$-set associated to $s$.

 Even if the $S$-sets in question can be defined at each step, it is still not clear that $Z_t$ will have the remarkable property  of being a disjoint union of $S$-sets for each choice of $s \in I$.

 The main goal of this paper is to show that this construction is compatible with adjoining the faces obtained from an appropriate $S$-set to give new trails.  Then the fact that trails are associated to a fundamental $\mathfrak g$ module should enable us to recover this remarkable decomposition property, in analogy with the simultaneous decomposition of an integrable module as a direct sum of simple $\mathfrak {sl}(2)$ modules for every simple root.
 %will allow us to establish the existence of a giant $S$-set of type $t$ required to have the remarkable property of being a disjoint union of $S$ sets for each choice of $s \in I$.

 %Finally to return to $\mathfrak g$ itself we simply interchange roots and co-roots in the expression for these functions.  Then the maximum of these functions on $B_J(\infty)$ will be $\varepsilon_t^\star$.

 %In this paper we can consider that we are constructing a giant $S$-set for the Langlands dual of $\mathfrak g$.  This is entirely consistent \cite {BZ} who consider trails in the fundamental modules for the Langlands dual in the first place.

 Unfortunately to carry out this procedure we are forced in the present paper to assume that there are no ``false trails'' in $V(-\varpi_t)$. Under this hypothesis the envelope $\mathscr K_t$ of the sought after giant $S$-set $Z_t$ will be constructed and moreover we show that it identifies with the set $\mathscr K^{BZ}_t$ of all Berenstein-Zelevinsky trails (Corollary \ref {8.3}). Conversely we show (Theorem \ref {9.1}) that the possibility to construct $\mathscr K^{BZ}_t$ by adjoining faces excludes the presence of false trails in $V(-\varpi_t)$.  This is a weaker condition than saying that our construction gives a giant $S$-set.  We expect (though at present cannot prove) that $Z_t$ to be just the set of extremal elements of $\mathscr K_t$, but this refinement is of no consequence in computing the dual Kashiwara parameter $\varepsilon^\star_t$  since it is defined by taking a maximum.

 %The notion of faces and of adjoining faces comes from the geometry of wiring diagrams.

\subsection{}\label{2.5}

Fix $s \in I$. The driving trail $K_t^1$ is rather special and sometimes has to be excluded.  For this we introduce the condition

\

$(H)$.   $s\neq t$ or $K \neq K_t^1$.

\

\begin {lemma}  Take $K \in \mathscr K^{BZ}_t$ and assume $(H)$.  Then

\

(i) $\alpha_s^\vee(\delta^K_{(s,1)})\leq 0$.

\

(ii) $\gamma^K_{(s,1)}-\alpha_s \notin \Omega(V(-\varpi_t))$.
%Again $\alpha_t^\vee(\delta^K_{(t,1)})\leq 0$, if $K \neq F_t^1$.
\end {lemma}

\begin {proof}  Observe that (i) follows from (ii), $(T)$ and the definition of $\delta^K_{(s,1)}$.

Suppose $s \neq t$.

%By \cite [$(1)$]{JZ}, it is immediate that $\gamma^{F_t^1}_{(s,1)}$ is an extremal weight and for all $s \in I \setminus \{t\}$ one has $\alpha_s^\vee(\gamma^{F_t^1}_{(s,1)}) \leq 0$. This gives (ii) for $K=F_t^1$.

By $(B)$ every trail $K$ starts at the weight $\gamma_1^K:=-s_t\varpi_t$.  Then by $(T)$ we obtain
$\gamma^K_{(s,1)} \in -s_t\varpi_t+\sum_{s'\in I|(s',1)<(s,1)}\mathbb N\alpha_{s'}\subset -\varpi_t+\mathbb N(\pi\setminus\{\alpha_s\})$.  Yet every weight of $V(-\varpi_t)$ lies in $-\varpi_t +\mathbb N\pi$ and so (ii) holds for all $K \in \mathscr K_t^{BZ}$.

Suppose $s=t$.

As above we obtain $\gamma^K_{(t,1)} =-s_t\varpi_t + \mathbb N(\pi\setminus\{\alpha_t\})$.  Then $\gamma^K_{(t,1)}-\alpha_t \in \Omega(V(-\varpi_t))$ implies that $\gamma^K_{(t,1)}-\alpha_t=-\varpi_t$ and so $\gamma^K_{(t,1)} =-s_t\varpi_t$.  Then by $(T)$ we obtain $\gamma^K_j=-s_t\varpi_t$, for all $j \leq (t,1)+1$. Let us show inductively that this implies that $\gamma_j^K=\gamma_j^{F_t^1}$ for all $j \in J$.

We may assume $j\geq (t,1)+1$. Then by the induction hypothesis $\gamma_j^K-\alpha_{i_j}=\gamma_j^{F_t^1}-\alpha_{i_j}$ is not a weight of $V(-\varpi_t)$ and so by $(T)$ we obtain $\gamma_{j+1}^K=\gamma_j^K+n_j\alpha_{i_j}$, where $n_j\leq -\alpha_{i_j}^\vee(\gamma^K_j)$. Then by the induction hypothesis and the definition of $\gamma_j^{F_t^1}$ we obtain $\gamma_{j+1}^K \in \gamma_{j+1}^{F_t^1}-\mathbb N \alpha_{i_j}$ and we conclude by $(P)$.

Thus $K=F_t^1$ which is excluded by $(H)$.  We conclude that (ii) holds for $s=t$ also.
\end {proof}

%\subsection{Faces}\label{2.6}
%
%The above lemma is needed to speak of minimal trails (Definition \ref {4.6} with respect to a fixed $s \in I$.  Condition $(H)$ is automatically satisfied since the driving trail $K_t^1$ is never considered to be minimal with respect to $t$.  If there is no

 \section{An $\mathfrak {sl}(2)$ computation}\label{3}

 \subsection{}\label{3.1}

 In order to study trails we need the following standard (or at least not far from standard) $\mathfrak {sl}(2)$ computation.  Here we take the standard basis $(e,h,f)$ for $\mathfrak {sl}(2)$ given by the relations $[e,f]=h,[h,e]=2e,[h,f]=-2f$.

  \subsection{}\label{3.2}

 Fix a positive integer $n$ and non-negative integers $a_i,k_i,\ell_i:i=1,2,\ldots,n$.  Let $e_{-a_i}$ be a lowest weight vector of $h$ eigenvalue $-a_i$ and let $V(-a_i)$ denote the simple finite-dimensional $\mathfrak {sl}(2)$ module generated by the action of $e$ on $e_{-a_i}$ viewed as a lowest weight vector.  Let $v_\textbf{k}$ denote the vector $e^{k_n}v_{-a_n}\otimes e^{k_{n-1}}v_{-a_{n-1}}\otimes \cdots \otimes e^{k_1}v_{-a_1}$ in the $n$-fold tensor product $V(-a_n)\otimes V(-a_{n-1})\otimes \cdots \otimes V(-a_1)$ given the ``diagonal'' action of $\mathfrak {sl}(2)$, that is to say that given by the Leibnitz rule or in modern terminology by the coproduct on $U(\mathfrak {sl}(2))$.

 For all $j=1,2,\ldots,n$, set $b_i=k_i-\ell_i$ and

 $$a^{(j)}=\sum_{i=1}^ja_i, \quad b^{(j)}=\sum_{i=1}^jb_i\quad k^{(j)}=\sum_{i=1}^jk_i,\quad \ell^{(j)}=\sum_{i=1}^j\ell_i.\eqno {(6)}$$

 Set $\textbf{a}=\{a_i\}_{i=1}^n, \textbf{b}=\{b_i\}_{i=1}^n, \textbf{k}=\{k_i\}_{i=1}^n,\textbf{l}=\{\ell_i\}_{i=1}^n, $.

  Let $A_\textbf{b}(\textbf{k},\textbf{l})$ denote the coefficient of $v_\textbf{l}$ in $f^bv_\textbf{k}$. This can be non-zero only if $b_i:=k_i-\ell_i\geq 0$, for all $i=1,2,\ldots,n$ and that $b=\sum_{i=1}^nb_i$, which we shall assume. Note that $b^{(n)}=b$.

 %Indeed let us set
% $$B_n(\textbf{k},\textbf{m}):=(-1)^nn!\prod_{i=1}^r {k_i \choose m_i}, \quad  C_\textbf{n}(\textbf{k},\textbf{m})=\prod_{j=1}^r\prod_{i=1}^{n_j}(i+k^{j-1}+m^j-1-a^j). \eqno {(11)}$$

 \begin {lemma}  One has $$A_\textbf{b}(\textbf{k},\textbf{l})=b!\prod_{i=1}^n {k_i \choose \ell_i}\prod_{j=1}^n\prod_{i=1}^{b_j}(a^{(j)}+1-i-k^{(j-1)}-\ell^{(j)}). \eqno {(7)}$$
 \end {lemma}

 \begin {proof}
  The proof by induction on $b$ and on $n$.

 Set
$$B_\textbf{b}(\textbf{k},\textbf{l}):=(-1)^bb!\prod_{i=1}^n {k_i \choose \ell_i}, \quad  C_\textbf{b}\textbf{}(\textbf{k},\textbf{l})=\prod_{j=1}^r\prod_{i=1}^{b_j}(i+k^{(j-1)}+\ell^{(j)}-1-a^{(j)}). \eqno {(8)}$$

 We need to show that
 $$A_\textbf{b}(\textbf{k},\textbf{l})=B_\textbf{b}(\textbf{k},\textbf{l})C_\textbf{b}(\textbf{k},\textbf{l}). \eqno {(9)}$$

 Let $\textbf{k}-\delta_t: t \in \{1,2,\ldots,n\}$ designate that $k_t$ has been reduced by $1$ leaving the remaining $k_i$ unchanged.  Give a similar meaning to $\textbf{b}-\delta_t$.

 Recall that for all $u \in \mathbb N^+$ one has $[f,e^u]=-e^{u-1}u(h+u-1)$.  From this we obtain the recurrence relation
 $$A_\textbf{b}(\textbf{k},\textbf{l})= - \sum_{t=1}^n k_t(k_t+2k^{(t-1)}-a^{(t)}-1)A_{\textbf{b}-\delta_t}(\textbf{k}-\delta_t,\textbf{l}).\eqno {(10)}$$

 The right hand side of $(10)$ may be computed through induction on $b$.

 Observe that the replacement of $\textbf{k}$ by $\textbf{k}-\delta_t$ replaces $k^{(j)}$ by $k^{(j)}-1$ if $j\geq t$.

 Consider $C_{\textbf{b}-\delta_t}(\textbf{k}-\delta_t,\textbf{l})$.  For $j<t$, we remove from  $\prod_{i=1}^{b_j}(i+k^{(j-1)}+\ell^{(j)}-1-a^{(j)})$ the factor corresponding to $i=b_j$, noting that $b_j+k^{(j-1)}+\ell^{(j)}=k^{(j)}+\ell^{(j-1)}$.  For $j>t$ we replace the dummy index $i$ by $i+1$ and remove the factor corresponding to $i=0$.  For $j=t$, we drop the factor for which $i=b_j$ entirely.

  Finally we note that ${k_t-1 \choose \ell_t}=\frac {k_t-\ell_t}{k_t}{k_t \choose \ell_t}$.

 We thus conclude that each term in the sum in the right hand side of $(10)$ has the factor
 $$\frac{1}{b^n}B_\textbf{b}(\textbf{k},\textbf{l})\prod_{j=1}^n\prod_{i=1}^{b_j-1}(i+k^{(j-1)}+
 \ell^{(j)}-1-a^{(j)}),\eqno {(11)}$$
 and up to this common factor the $t^{th}$ term in the sum in $(10)$ takes the form
 $$(k_t-\ell_t)(k_t+2k^{(t-1)}-a^{(t)}-1)\prod_{j=1}^{t-1}(k^{(j)}+\ell^{(j-1)}-1-a^{(j)})
 \prod_{j=t+1}^{n}(k^{(j-1)}+\ell^{(j)}-1-a^{(j)}).\eqno {(12)}$$

 %In this we used the identity $n_j+k^{j-1}+m^j=k^j+m^{j-1}$ having obtained the first product by evaluating at $i=n_j$ and the second by replacing the dummy index $i$ by $i+1$ and then evaluating at $i=0$.
 Now for $t<n$, the terms in $(12)$ have the common factor $(k^{(n-1)}+\ell^{(n)}-1-a^{(n)})$.  Up to this common factor, they are just the terms we obtain from the right hand side of $(10)$ by reducing $n$ by $1$, up to the common factor given by $(11)$ with $n$ is reduced by $1$.

 By the induction hypothesis on $n$, this sum (with the said common factor) must equal the left hand side of $(12)$ in which $n$ is reduced by $1$.  We conclude that the sum of the first $n-1$ terms in $(12)$, must supply the missing factor for which $i=b_j$ in $(11)$ (with $n$ reduced by $1$).  Hence this sum  is equal to
 $$b^{(n-1)}(k^{(n-1)}+\ell^{(n})-1-a^{(n)})\prod_{j=1}^{n-1}(k^{(j)}+\ell^{(j-1)}-1-a^{(j)}). \eqno {(13)}$$
On the other hand the last term in $(12)$ is just
 $$(k_n-\ell_n)(k_n+2k^{(n-1)}-1-a^{(n)})\prod_{j=1}^{n-1}(k^{(j)}+\ell^{(j-1)}-1-a^{(j)}). \eqno {(14)}$$

 The expressions in $(13),(14)$ share the product as a common factor and so we have only to note that the sum of the quadratic expressions which proceeds them
is just $b^{(n)}(k^{(n)}+\ell^{(n-1)}-1-a^{(n)})$.  It follows that the sum of the terms in $(13)$ and $(14)$ equals $b\prod_{j=1}^n(b_j+k^{(j-1)}+\ell^{(j)}-1-a^{(j)})$.  Multiplication by the common factor given in $(11)$ recovers the asserted expression for $A_\textbf{b}(\textbf{k},\textbf{l})$.
  \end {proof}

   \subsection{}\label{3.3}

   From now on we view $n$ above is a fixed element of $\mathbb N^+$ and we set $N=\{1,2,\ldots,n-1\}, \hat{N}=\{1,2,\ldots,n\}$.

 \section{The Minimax Theorem}\label{4}

 \subsection{}\label{4.1}

 Recall the notation and hypotheses of \ref {2.4}, particularly the integers $u,v$ and recall $(3)$.  In order to adjoin the face $F_s^k$ to a trail $K$ to obtain a new trail $K+F_s^k$, we must replace $v^K_{(s,k-1)+1}$ by $e_sv^K_{(s,k-1)+1}$ and make the consequent replacements in the $v_j:v<j\leq u$.  This meets three obstructions.  First we must be sure that these vectors are all non-zero.  Second we need to know that $\gamma^K_{(s,k)}+\alpha_s$ does not exceed $\gamma^K_{(s,k)+1}$. This implies that $v^K_{(s,k)+1}=e_s^bv^K_{(s,k)}$, with $b\in \mathbb N^+$.  Finally we need to show, up to a non-zero scalar, that $v^K_{(s,k)+1}$ is also given by applying $e^{b-1}_s$ to the vector which replaces $v^K_{(s,k)}$.  We call this the matching condition at $(s,k)$.

 \subsection{}\label{4.2}

 At first sight it might seem that the matching condition could only be satisfied by a miracle. Actually it follows from the second boundary condition and an easy fact about Demazure modules.  Let us examine this in more detail.

 Consider a trail $K$ which trivializes at some $j\in J$.  Then we may write
 $$v_{-w_j\varpi_t}=v^K_{j+1}=e^{n_j}_{i_j}e^{n_{j-1}}_{i_{j-1}}\cdots e^{n_1}_{i_1}v_{-\varpi_t}. \eqno {(15)}$$

 Set $w=w_j$ and $s=i_j$.

 We introduce some further notation which will be compatible with that of Section \ref {3}.

 We write $j=(s,n)$ and identify $N$ with the set $\{(s,i)\}_{i=1}^n$.  The simple weight vectors $e_j:(s,i) < j < (s,i+1)$
%Define the finite subset $J_s=\{k \in J|1\leq k \leq j: i_k=s\}$ of $J$.  It inherits a total order from $J$.  Identify $J_s$ with the set $j,j-1,\ldots,1$.  Between successive elements $i,i-1:j\geq i>1$ of $J_s$, the simple weight vectors
  are distinct from $e_s$ and their product may be written as a weight vector $e_{\mu_i}$ of $U(\mathfrak n)$ of weight $\mu_i$.  This product commutes with $f_s$.

  In $(15)$ we must recall that a trail is defined to start at $v_{-s_t\varpi_t}$ and not at $v_{-\varpi_t}$.  Because of this we must assume that $(H)$ holds. Then by Lemma \ref {2.5}(ii) we have $f_sv^K_{(s,1)}=0$, so we can take $v_{\mu_1}=v^K_{(s,1)}$.  
  %When $s\neq t$ this is equivalent to taking $v_{\mu_1}$ to be the product of the simple weight vectors to the right of the smallest element of $J_s$ applied to $v_{-\varpi_t}$.

   Let $\textbf{e}$ denote the $n$-tuple consisting to the products $\{e_{\mu_i}\}_{i=2}^n$ and the weight vector $v_{\mu_1} \in V(-\varpi_t)$.  Fixing $\textbf{e}$, we may consider the set of all \textit{non-zero} vectors described by $(16)$, where only the set of exponents $\textbf{k}:= (k_n,k_{n-1},\ldots,k_1)$ are varied, as corresponding to a subset $T_s(\textbf{e})$ of trails trivializing at $w$.

 %Again the product of the simple weight vectors to the right of the smallest element of $J_s$ applied to $v_{-\varpi_t}$ is an element $v_{\mu_1}$ of $V(-\varpi_t)$ annihilated by $f_s$ of weight $\mu_1$.
%
% \textit{However} this argument is not quite correct if $s=t$, since a trail is defined to start at $v_{-s_t\varpi_t}$ and not at $v_{-\varpi_t}$.  Yet when $(H)$ holds, it follows from Lemma \ref {2.5}(ii) that $f_sv^K_{(s,1)}=0$, so we can take $v_{\mu_1}=v^K_{(s,1)}$.  When $s\neq t$, this holds automatically from definitions.

  Set $a_i:=-\alpha^\vee_s(\mu_i)$ which is non-negative integer.  Let $\textbf{a}$ denote the  $r$-tuple \newline $(a_1,a_2,\ldots,a_r)$ of non-negative integers.

 Set $h_s:=[e_s,f_s]$.  Then $\{e_s,h_s,f_s\}$ is a basis for an $\mathfrak {sl}(2)$ subalgebra $\mathfrak s_s$ of $\mathfrak g$, which we assume to be of standard form (cf. \ref {3.1}).   We shall often drop the $s$ subscript.

 We write $e_{\mu_i}$ simply as $e_{-a_i}$, for $r\geq i>1$ and $v_{\mu_1}$ simply as $v_{-a_1}$.  Then we may write $v_{-w\varpi_t}$ as $\overline{v}_{\textbf{k}}$, where
 $$\overline{v}_{\textbf{k}}:=e^{k_n}e_{-a_n}e^{k_{n-1}}e_{-a_{n-1}}\cdots e^{k_1}v_{-a_1}. \eqno {(16)}$$

% Let $\textbf{e}$ denote the $n$-tuple consisting to the products $\{e_{\mu_i}\}_{i=2}^n$ and the weight vector $v_{\mu_1} \in V(-\varpi_t)$.  Fixing $\textbf{e}$, we may consider the set of all \textit{non-zero} vectors described by $(16)$, where only the set of exponents $\textbf{k}:= (k_n,k_{n-1},\ldots,k_1)$ are varied, as corresponding to a subset $T_s(\textbf{e})$ of trails (trivializing at $w$).

 %\
%
% \textbf{Convention}.   Up to the end of Section \ref {7} only the values given by $\textbf{a}$ are of any importance and so we replace $\textbf{e}$ by $\textbf{a}$.
%
% \

 %Now fix $\textbf{a}=\{a_r,a_{r-1},\ldots,a_1\}$, or more precisely fix the vector $v_{\mu_1}$ and the products $\{e_{\mu_i}\}_{i=2}^r$.  Then we may regard  the set of all \textit{non-zero} vectors described by $(16)$, where only the set of exponents $\textbf{k}:= (k_r,k_{r-1},\ldots,k_1)$ are varied, as corresponding to a subset $T_s(\textbf{a})$ of trails (trivializing at $w$).

 The vector $v_\textbf{k}$ in the n-fold tensor product $V(-a_n)\otimes V(-a_{n-1})\otimes \cdots \otimes V(-a_1)$ defined in \ref {3.2} is non-zero if and only if $k_i\leq a_i$, for all $i \in \hat{N}$.

 Observe that there is an $\mathfrak s$ module map $\varphi$ of $V(-a_n)\otimes V(-a_{n-1})\otimes \cdots \otimes V(-a_1)$ into $V(-\varpi_t)$ taking $v_\textbf{k}$ to $\overline{v}_\textbf{k}$.

The\textit{ hard part} of the present work is to determine when $\overline{v}_\textbf{k}$ is non-zero.

 We may consider $v_\textbf{k}$ to be the inverse image of $\overline{v}_\textbf{k}\in V(-\varpi_t)$ in the tensor product even when the latter is zero, since its definition is just a matter of specifying $\textbf{k}$.

 We shall not always distinguish between the elements of $T_s(\textbf{e})$ and the vectors $\overline{v}_\textbf{k}$ they define.

  Let $M_s(\textbf{e})$, or simply $M$, denote their linear span of the elements of $T_s(\textbf{e})$.  It is clearly an $\mathfrak s$ module.

 % \
%
%
% \textbf{Convention}.   Up to the end of Section \ref {7} only the values given by $\textbf{a}$ are of any importance and so we replace $\textbf{e}$ by $\textbf{a}$.
%
%\

  By definition $M_s(\textbf{e})= \varphi(V(-a_n)\otimes V(-a_{n-1})\otimes \cdots \otimes V(-a_1))$.

  %We have the following remarkable fact, which depends crucially on the second boundary condition $(B)(ii)$ which imposes that eventually $\gamma^K_{j+1}$ is an extremal vector $-w\varpi_t$ with $w\in W$ of length $j$ - at which point we are saying that $K$ trivializes at $w$).  This forces a trail to move efficiently (that is without taking to many steps).  Of course since a fundamental weight has a large stabilizer in $W$, there are nevertheless many possible trails trivializing at a given element of $W$ and these need not only pass through extremal weights.

   %Moreover the map taking $v_{\textbf{k}}\in V(-a_r)\otimes V(-a_{r-1})\otimes \cdots \otimes V(-a_1)$ to $\overline{v}_{\textbf{k}}\in X(\textbf{a})$ extends linearly to a surjection of $V(-a_r)\otimes V(-a_{r-1})\otimes \cdots \otimes V(-a_1)$ onto $M_s(\textbf{a})$.  Quite remarkably

 \begin {lemma}  $M_s(\textbf{e})$ is a simple $\mathfrak s$ module.
 \end {lemma}

 \begin {proof}  Let $\mathfrak n^-$ be the opposed algebra to $\mathfrak n$ in $\mathfrak g$.  Recall that $w \in W$ has reduced decomposition $w=s_{i_j}s_{i_{j-1}}\cdots s_{i_1}$.  Set $V_w(-\varpi_t):=U(\mathfrak n^-)v_{-w\varpi_t}$, which is a so-called Demazure module.  It is a classical fact (see \cite [4.4.6(i)]{J2}) that $V_w(-\varpi_t)=k[e_{i_j}]k[e_{i_{j-1}}]\cdots k[e_{i_1}]v_{-\varpi_t}$.  Thus from our construction $v \in M$ belongs to $U(\mathfrak b^-)v_{-w\varpi_t}$.  Moreover if $v$ is a weight vector then its weight must lie in $-w\varpi_t+\mathbb Z \alpha_s$.  Yet the only weight vectors of weight $\mathbb Z \alpha_s$ lying in $U(\mathfrak n^-)$ are the powers of $f_s$.  Thus $v=f_s^nv_{-w\varpi_t}$ for some $n \in \mathbb N$, from which the simplicity of $M$ follows.
 \end {proof}

 \textbf{Remark}.  This is a remarkable fact, which depends crucially on the second boundary condition $B(ii)$ which imposes that eventually $\gamma^K_{j+1}$ is an extremal vector $-w\varpi_t$ with $w\in W$ of length $j$ - at which point we are saying that $K$ trivializes at $w$).  This forces a trail to move efficiently (that is without taking to many steps).  Of course since a fundamental weight has a large stabilizer in $W$, there are nevertheless many possible trails trivializing at a given element of $W$ and these need not only pass through extremal weights.  The simplest example of the latter occurs in type $B_2$ for its five dimensional fundamental module.

 This fact will play a significant role in our work (for example see Corollary \ref {6.2}) as well as providing an important guideline.  Indeed if vectors $\overline{v}_\textbf{k},\overline{v}_{\textbf{k}'}\in M_s(\textbf{a})$ have the same weight, then they must be proportional and it is natural to try to compute these proportionality factors.  This was one motivation for the present work.  We provide the solution when there are no false trails - see in particular Lemma \ref {5.2} and section \ref {7}.

 \subsection{}\label{4.3}

 It is clear that $v_{- w \varpi_t} \in M_s(\textbf{e})$ (resp. $v_{-s_\alpha w \varpi_t} \in M_s(\textbf{e})$) and is its unique up to scalars highest (resp. lowest) weight vector.  Any choice of $\textbf{k}$ for which $\overline{v}_{\textbf{k}}= v_{-s_\alpha w \varpi_t}$, up to a non-zero scalar, has the property that $k_n=0$; but the converse is false.

 We denote by $T_s^+(\textbf{e})$, (resp. $T_s^-(\textbf{e})$)  the subset of $T_s(\textbf{e})$, of all vectors $\overline{v}_\textbf{k}$ proportional to the highest (resp. lowest) weight vector of $M_s(\textbf{e})$.

 \

 \textbf{Convention}.   Up to the end of Section \ref {7} only the values given by $\textbf{a}$ are of any importance and so we replace $\textbf{e}$ by $\textbf{a}$.

\

 %We extend the convention of \ref {4.2} by writing $T_s^+(\textbf{e})$ as $T_s^+(\textbf{a})$ (resp. $T_s^-(\textbf{e})$ as $T_s^-(\textbf{a})$) up to the end of Section \ref {7}.

 We may and do regard $T_s^+(\textbf{a})$  (resp. $T_s^-(\textbf{a})$) as the set of all trails defined by the pair $(s,\textbf{a})$ trivializing at $v_{-w\varpi_t}$ (resp. $v_{-s_\alpha w\varpi_t}$).  Any element of $T_s(\textbf{a})$ becomes an element of $T_s^+(\textbf{a})$ by multiplying on the left by a suitable power of $e_s$.   In this sense we may view $T_s(\textbf{a})$ and $T_s^-(\textbf{a})$ as subsets of $T^+_s(\textbf{a})$.

 Recall the notation \ref {4.1} and let us describe how to adjoin (resp. remove)  the face $F_s^\ell:\ell>1$ from the trail $K$.  Indeed by $(B)$ we may assume that $K$ trivializes at $j \in J$ and is hence represented by $(16)$.  Then adjoining (resp. removing) the face $F_s^\ell$ to (resp. from) $K$ is simply a matter of taking a factor of $e_s$ through $e_{-a_\ell}$ from left to right (resp. right to left).  Of course to be able to do this we need to know that $k_\ell >1$ (resp. $k_{\ell-1}>0$), which is the condition we already met.  Another ``minor matter'' is to show that the resulting expression is non-zero which we referred to in \ref {4.1} as our first obstruction.  However given this, the matching condition at $(s,\ell)$ is automatically satisfied since non-zero vectors in $M_s(\textbf{a})$ of the same weight, which are just the elements of $T_s(\textbf{a})$, are proportional.  In particular all the elements of $T_s^+(\textbf{a})$ viewed as elements of $M_s(\textbf{a})$ are proportional.

  \subsection{}\label{4.4}

  The remainder of this paper is dedicated to trying to resolve the ``minor matter'' brought up in \ref {4.3}.  By the remark in \ref {4.2} and an induction argument on the length of the element $w\in W$ at which a given trail trivializes, it is equivalent to determining which expressions in $(15)$ (and hence in $(16)$ for every $s \in I$) are non-zero.

  The solution to the above problem we have in mind is that the non-zero expressions in $(16)$ are just the integer points of the convex hull of the $S$-set determined by $T_s(\textbf{a})$.  The next section addresses how this $S$-set should be described.

   \subsection{}\label{4.5}

   An $S$-set is determined by a set of parameters $c_1,c_2,\ldots,c_{n-1} \in \mathbb N$.  The latter are computed from its minimal element.  Therefore what we need is a partial order on $T^+_s(\textbf{a})$ compatible with adjunction of faces of the given type $s$, with the property that it admits a unique minimal element.

   Observe that there an obvious total ordering on $T^+_s(\textbf{a})$
   compatible with adjunction of faces of type $s$.  It is just that induced by the lexicographic order on $\mathbb N^n$, when we view an element $\overline{v}_\textbf{k}$ as being determined by an element $\textbf{k}\in \mathbb N^n$.

   \

   \textbf{Definition}.  The unique trail in $T^+_s(\textbf{a})$ which is minimal (resp. maximal) for the lexicographic ordering is called the $\ell$-minimal (resp. $\ell$-maximal) trail $K_{\ell \min}$ (resp. $K_{\ell \max}$) .

  \subsection{}\label{4.6}

  For the moment we shall consider a conceptually simpler definition of a minimal (resp. maximal) trail).

  Fix $s \in I$.

  \

  \textbf{ Definition}.   A trail $K \in \mathscr K^{BZ}$ is said to be minimal if  $f_sv^K_{(s,k)}=0$, for all $k \in \mathbb N^+$.

  \

  In terms of $(16)$, this condition is just that
 $$f_se_{-a_i}e_s^{k_{i-1}}e_{-a_{i-1}}\cdots e_s^{k_1}v_{-a_1}=0, \forall i=1,2\ldots,n. \eqno{(17)}$$

 \

 One may remark by $(1)$ and Lemma \ref {2.3}(ii) that $K_t^1$ is a minimal trail for all $s\neq t$.

 A difficulty is that $T^+_s(\textbf{a})$ does not obviously admit a minimal trail. However if it does then it is clearly $\ell$-minimal.  A related difficulty is that a minimal trail is not combinatorially defined that is to say it cannot be read off from the knowledge of $T^+_s(\textbf{a})$.

%One may remark by $(1)$ and Lemma \ref {2.3}(ii) that $K_t^1$ is a minimal trail for all $s\neq t$.
 \

 \textbf{ Definition}  A trail $K \in \mathscr K^{BZ}$ is said to be maximal if  $e_sv^K_{(s,k)+1}=0$, for all $k \in \mathbb N^+$.

 \

  In terms of $(16)$, this condition is just that
 $$e_se_s^{k_i}e_{-a_i}e_s^{k_{i-1}}e_{-a_{i-1}}\cdots e_s^{k_1}v_{-a_1}=0, \forall i=1,2\ldots,n. \eqno{(18)}$$

 \

A difficulty is that $T^+_s(\textbf{a})$ does not obviously admit a maximal trail.  However if it does then it is clearly  $\ell$-maximal.

 Take $\overline{v}_{\textbf{k}} \in T^+_s(\textbf{a})$. Then $(17)$ always holds for $i=1$.  Here we remark that condition $(H)$ can always be assumed since the driving trail $K_t^1$ is never considered to be an element of $T_s(\textbf{a})$ for $s=t$.  Again $(18)$ holds for $i=r$ since $\overline{v}_\textbf{k}$ is a highest weight vector of $M_s(\textbf{a})$.

 When $K_t^1$ is not minimal with respect to any $s \in I$, then $z_t^1=m_{t,1}$.  In this case $Z_t$ is reduced to $\{z_t^1\}$.  This holds in particular when $(t,1)=1$. When $W$ is finite, one can always find a reduced decomposition of $w_0$ so that this latter condition holds.  However the change of parametrization of $B(\infty)$ is extremely complicated.  Yet it can be handled by the tropical calculus and this was the basis behind \cite [Thm. 3.9] {BZ} to describe $Z_t$ using trails.   When $W$ is not finite, it is not at all clear why trails should give a correct description of $Z_t$, nor indeed even for $W$ finite, that trails should form a set which is so highly redundant.

\subsection{}\label{4.7}

The minimax theorem (Theorem \ref {4.7.4}) states that if $T^+_s(\textbf{a})$ admits a maximal element, then it admits a minimal element, whilst the converse holds by Lemma \ref {4.7.5}. In this the parameters which describe the minimal trail define an $S$-set (cf \ref {4.5}). In addition the maximal trail exactly corresponds to the unique maximal element of the $S$-set.  The various parts of this theorem will be handled separately in the subsections below.

The minimax theorem will not be used after section \ref {4.7} and so some may want to skip this section.  It could be useful if an $\ell$-minimal trail could be shown to be minimal.   We believe that it makes a useful preamble.

In the proof of the minimax theorem, Lemma \ref {4.2} will not be used, nor the full force of Lemma \ref {3.2}.  Again we shall not use the Chevalley-Serre relations.  All these will be needed when minimal trails are replaced by $\ell$-minimal trails, so presumably the latter theory is deeper and indeed goes much further.

Throughout we fix $s \in I$ and assume $(H)$.

\subsubsection{}\label{4.7.1}  Given integers $i\leq j$, set $[i,j]=\{i,i+1,\ldots,j\}$.

%Assume $(H)$ holds.  Fix $s \in I$ and let $\mathfrak s$ denote the $\mathfrak {sl}(2)$ subalgebra of $\mathfrak g$ generated by $e_s,f_s$ with $h_s=[e_s,f_s]$.

\begin {lemma}   Take $K \in \mathscr K^{BZ}_t$ and suppose there exists $k \in \mathbb N$ such that  $e_sv^K_{(s,k+1)+1} =0$ and $f_sv^K_{(s,k)}=0$. If $f_sv^K_{(s,k+1)} \neq 0$, then $k>0$ and one may remove a copy of the face $F_s^{k+1}$ from $K$ to obtain $K'\in \mathscr K^{BZ}_t$ by setting
$$ v_j^{K'}=\left\{
                       \begin{array}{ll}
                         f_sv_j^K, & \hbox{if $j \in [(s,k)+1,(s,k+1)]$,} \\
                         v_j^K, & \hbox{otherwise.}
                       \end{array}
                     \right. $$
Moreover this process may be repeated till $f_sv^{K'}_{(s,k+1)} =0$, for the new trail $K'$ obtained from $K$ by removing finitely many copies of $F_s^{k+1}$.
\end {lemma}

\begin {proof}  That $k>1$ follows from Lemma \ref {2.5}.  By definition of a trail there exists $b \in \mathbb N$ and a product $e_{-a_2}$ of the simple root vectors distinct from $e_s$ of $h_s$ eigenvalue $-a_2$ such that $v^K_{(s,k+1)}=e_{-a_2}v^K_{(s,k)+1}$ and $v^K_{(s,k)+1}=e^b_sv^K_{(s,k)}$.  Since $e_{-a_2}$ commutes with $f_s$, one has $a_2 \in \mathbb N$.

Choose $p \in \mathbb N$ maximal such that $f_s^pv^K_{(s,k+1)}=f^p_se_{-a_2}e_s^bv^K_{(s,k)}\neq 0$. Yet $f_sv^K_{(s,k)}=0$, so $1\leq p \leq b$.  Since $v_{(s,k)}$ is a lowest weight vector, it follows from $\mathfrak {sl}(2)$ theory that $f_s^pv^K_{(s,k+1)}=e_{-a_2}f_s^pv^K_{(s,k+1)}$ is a non-zero multiple of $e_{-a_2}e_s^{b-p}v^K_{(s,k)}$.  Set $v_{-a_1}=e_s^{n-p}v^K_{(s,k)}$.  One has $f_se_{-a_2}v_{-a_1}=0$, so then $a_1+a_2 \in \mathbb N$.

We conclude that the $v_{(k_2,k_1)}:=e^{k_2}_se_{-a_2}e_s^{k_1}v_{-a_1}:k_2,k_1 \in \mathbb N$ span an $\mathfrak s$ module not necessarily simple but having a unique up to scalars vector of the lowest possible weight, namely $v=v_{(0,0)}$. Since $v_{-a_1}$ is non-zero of weight $-a_1$ and $e_{-a_2}v_{-a_1}$ is non-zero of weight $-a_2-a_1$, the non-vanishing of $v_{(k_2,k_1)}$ implies that $k_1 \leq a_1$ and $k_1+k_2 \leq a_1+a_2$.  The converse, which is slightly less obvious, obtains by examining the zeros in the right hand side of $(7)$, taking $n=2,\textbf{l}=\textbf{0}$.

On the other hand there exists $q \in \mathbb N$ maximal such that $e_s^qv^K_{(s,k+1)}\neq 0$.  By the hypothesis $e_sv^K_{(s,k+1)+1} =0$, the resulting vector can be taken to be $v^K_{(s,k+1)+1}$.

When $k_1\leq a_1$ and $k_1+k_2=a_1+a_2$ the above vectors are all non-zero. All have weight $a_1+a_2$, which is the negative of the weight of $v$. Then since $v$ is unique up to scalars and of lowest possible weight, it follows that these vectors are all proportional.  One of these vectors namely $v_{(q,p)}$ is just $v^K_{(s,k+1)+1}$.  Then the non-vanishing and proportionality of the vectors $v_{(q+i,p-i)}: i =0,1,\ldots,p$, establishes the lemma.

\end {proof}

\subsubsection{}\label{4.7.2}

%Fix $s \in I$ and assume $(H)$.

Assume that $T^+_s(\textbf{a})$ admits a minimal trail $K_{\min}$ (necessarily unique).

 Set $c^{K_{\min}}_k(s):=-\alpha^\vee_s(\delta^{K_{\min}}_{(s,k)}):k \in \mathbb N^+$.  Observe that $c^{K_{\min}}_k(s)=0$ is non-zero, if $k\geq n$.

\begin {lemma} For all $k \in \mathbb N^+$, one has
$c^{K_{\min}}_k(s) \geq 0$.

\end {lemma}

\begin {proof}  This follows from the hypothesis that $f_sv^{K_{\min}}_{(s,k)}=0$, $(T)$ and $\mathfrak {sl}(2)$ theory.
   \end {proof}

   \subsubsection{}\label{4.7.3}

  Assume that $T^+_s(\textbf{a})$ admits a maximal trail $K_{\max}$ (necessarily unique).
  Set $d^{K_{\max}}_k(s):=\alpha^\vee_s(\delta^{K_{\max}}_{(s,k)}):k \in \mathbb N^+$.

   Observe that $d^{K_{\max}}_k(s)=0$ is non-zero, if $k > n$.

 \begin {lemma} For all $k \in \mathbb N^+$, one has

    \

    (i)  $d^{K_{\max}}_k(s) \geq 0$.

    \

    (ii)  $d^{K_{\max}}_1(s)=0$.

    \

    (iii) $f_sv^{K_{\max}}_{(s,1)}=0$.
   \end {lemma}

   \begin {proof}  (i) follows from the hypothesis that $e_sv^{K_{\max}}_{(s,k)+1}=0$, $(T)$ and $\mathfrak {sl}(2)$ theory.  Then (ii) (resp. (iii)) follow from (i) (resp. (ii)) of Lemma \ref {2.5}.
   \end {proof}

\subsubsection{}\label{4.7.4}
%
%
%
%We need the following useful
%
%\begin {lemma}  For all $K \in \mathscr K^{BZ}_t, j \in J$ one has
%$$s_{i_j}\gamma_{j+1}^K=\gamma_j^K-\alpha_{i_j}^\vee(\delta^K_j)\alpha_{i_j}.$$
%In particular the right hand side is a weight of $V(-\varpi_t)$.
%
%\end {lemma}
%
%\begin {proof}  Set $\alpha_{i_j}=\alpha$ and $s=s_{i_j}$ and omit the superscript $K$.  There exists $q_{j+1}\in \mathbb Q$ and $ \varpi \in P(\pi)$ such that $\gamma_{j+1}=q_{j+1}\alpha+\varpi$ and $\alpha^\vee(\varpi)=0$. By $(C)$ we obtain $\gamma_{j}=q_j\alpha+\varpi$, with $q_{j+1}+q_j=\alpha^\vee(\delta_j)$.  Yet $s\gamma_{j+1}=\gamma_j-(q_{j+1}+q_j)\alpha$, so the assertion obtains.
%\end {proof}

\begin {thm}  Assume that $T^+_s(\textbf{a})$ admits a maximal trail $K_{\max}$.

One may successively subtract with $k \in \mathbb N^+$ increasing, faces of $F_s^{k+1}$ from $K_{\max}$ to obtain a minimal trail $K_{\min}$ in $T^+_s(\textbf{a})$.

Moreover %setting $c^{K_{\min}}_k(s):=-\alpha^\vee_s(\delta^{K_{\min}}_{(s,k)}):k \in \mathbb N^+$ with $c^{K_{\min}}_0(s)=0$ one has

\

   (i)  $d^{K_{\max}}_{k+1}(s)=c^{K_{\min}}_k(s)$, for all $k \in \mathbb N$.

   \

   (ii) $K_{\max}=K_{\min}+\sum_{k=1}^\infty c_k^{K_{\min}}(s)F_s^{k+1}$, with the sum being finite.

\end {thm}

\begin {proof}  Write $d^{K_{\max}}_k(s)$ simply as $d_k$. Recall that $e_sv_{(s,k)+1}^{K_{\max}}=0$.

Set $K_1=K_{\max}$. For all $\ell> 1$ we show inductively that one may subtract $d_\ell$ copies of $F_s^\ell$ from $K_{\ell-1}$ and setting $K_\ell =K_{\ell-1} - d_\ell F_s^{\ell}$, that
$$f_sv^{K_\ell}_{(s,k)}=0, \forall \ell\geq k. \eqno{(19)}$$

$$v_j^{K_\ell}=v_j^{K_{\max}}, \forall j>(s,\ell).\eqno {(20)}$$

$$ \alpha_s^\vee(\delta^{K_\ell}_{(s,k)})=\left\{
                       \begin{array}{ll}
                         -d_{k+1}, & \hbox{if $k<\ell$,} \\
                         0, & \hbox{if $k=\ell$,} \\
                         d_{k}, & \hbox{if $k>\ell$.}
                       \end{array}
                     \right. \eqno {(21)}$$

   Take $\ell=1$ and recall that $K_1=K_{max}$. Then $(19)$ holds by Lemma \ref {4.7.3}(iii), whilst $(20)$ holds by definition of $K_1$.  Finally $(21)$ holds by Lemma \ref {4.7.3}(ii) and the definition of $d_k$.
%, where $\ell_{\ell-1} \in \mathbb N$ is maximal such that $K_\ell \in \mathscr K^{BZ}_t$.

  Since $s=i_{(s,\ell)}$, we obtain $s_{\alpha_s}\gamma^{K_{\ell}}_{(s,\ell+1)+1}=\gamma^{K_{\ell}}_{(s,\ell+1)}-d_{\ell+1}\alpha_s$, by Lemma \ref {2.3} and $(21))$.  On the other hand $e_sv_{(s,\ell+1)+1}^{K_\ell}=0$ by $(20)$ and the definition of $K_{\max}$. Thus $d_{\ell+1}$ is exactly the largest value of $b \in \mathbb N$ such that $f_s^bv_{(s,\ell+1)}^{K_{\ell}}\neq 0$, the latter being a vector of weight $s_{\alpha_s}\gamma^{K_{\ell}}_{(s,\ell+1)+1}$. Yet by $(19)$ one has $f_sv^{K_\ell}_{(s,\ell)} =0$.  Thus by Lemma \ref {4.7.1}, we may remove $d_{\ell+1}$ copies of $F_s^{\ell+1}$ from $K_\ell$ to obtain $K_{\ell+1}$.

  In this $v_{(s,\ell+1)}^{K_{\ell+1}}=f_s^{d_\ell}v_{(s,\ell+1)}^{K_{\ell}}$, so then $f_sv_{(s,\ell+1)}^{K_{\ell+1}}=0$, whilst $v_{(s,p)}^{K_{\ell+1}}=v_{(s,p)}^{K_{\ell}}$, for all $p \leq \ell$. This gives $(19)$ at the next induction step.  Again $v_j^{K_{\ell+1}}=v_j^{K_{\ell}}$, for all $j \notin [(s,\ell)+1, (s,\ell+1)]$.  This gives $(20)$ at the next induction step and the first and third lines of $(21)$. Again $\alpha_s^\vee(\delta^{K_{\ell+1}}_{(s,\ell+1)})=\alpha_s^\vee(\delta^{K_{\ell}}_{(s,\ell+1)})-d_{\ell+1}\alpha_s^\vee
  (\frac{1}{2}(0+\alpha_s))=0$, by the construction and the induction hypothesis. This is the second line of $(21)$.

  Finally by the observation in the first part of \ref {4.7.3}, $d_\ell=0$, for all $\ell>n$. Then $K_\ell=K_{n+1}$, for all $\ell>n$.  It follows from $(i)$ that $K_{\min}:=K_{n+1}$ is a minimal trail in $T_s(\textbf{a})$.  By $(21)$ it satisfies (i) and furthermore (ii) holds by the construction.

  \end {proof}

  \subsubsection{}\label{4.7.5}

  To describe a converse to Lemma \ref {4.7.1} a little more notation is needed.  Recall as in \ref {4.7.1} that by definition of a trail there exists a product $e_{-a_2}$ of the simple root vectors distinct from $e_s$ of $h_s$ eigenvalue $-a_2$ such that $v^K_{(s,k+1)}=e_{-a_2}v^K_{(s,k)+1}$.  Moreover for each $j \in [(s,k)+1,(s,k+1)]$ there is a sub-product $e[j]$ of $e_{-a_2}$ obtained by deleting simple root vectors from the \textit{left} such that $v^K_j=e[j]v^K_{(s,k)+1}$.  In this $e[(s,k)+1]=1,e[(s,k+1)]=e_{-a_2}$.

\begin {lemma}   Take $K \in \mathscr K^{BZ}_t$ and suppose there exists $k \in \mathbb N$ such that  $e_sv^K_{(s,k+1)+1} =0$ and $f_sv^K_{(s,k)}=0$.  Let $u \in \mathbb N$ be maximal such that $e_s^uv^K_{(s,k)+1}\neq 0$. Then for all $\ell=1,2, \ldots, u$, one may adjoin $\ell$ copies of the face $F_s^{k+1}$ to $K$ to obtain $K_\ell \in \mathscr K^{BZ}_t$ by setting
$$ v_j^{K_\ell}=\left\{
                       \begin{array}{ll}
                         e[j]e^\ell_sv_{(s,k)+1}^K, & \hbox{if $j \in [(s,k)+1,(s,k+1)]$,} \\
                         v_j^K, & \hbox{otherwise.}
                       \end{array}
                     \right. $$

Furthermore by construction $e_sv^{K_u}_{(s,k)+1}=0$.

\end {lemma}

\begin {proof}  Since the hypotheses of the lemma carry over to $K_1$ with $u$ reduced by $1$, it is enough to prove this assertion for $\ell =1$.  Again to show that $e[j]e_sv_{(s,k)+1}^K\neq 0$, it is enough to show that $e_{-a_2}e_sv_{(s,k)+1}^K\neq 0$.

The remainder of the proof follows closely that of Lemma \ref {4.7.1}.  As there we can write  $v^K_{(s,k+1)}=e_{-a_2}v^K_{(s,k)+1},v^K_{(s,k)+1}=e^b_sv^K_{(s,k)},v^K_{(s,k+1)+1}=
e_s^qv^K_{(s,k)+1}$, for some $b,q \in \mathbb N$. Again we can choose $p \in \mathbb N$ maximal such that $f_s^pv^K_{(s,k+1)}=f^p_se_{-a_2}e_s^bv^K_{(s,k)}\neq 0$. Then  $1\leq p \leq b$ and $f_s^pv^K_{(s,k+1)}=e_{-a_2}f_s^pv^K_{(s,k+1)}$ is a non-zero multiple of $e_{-a_2}e_s^{b-p}v^K_{(s,k)}$.  Set $v_{-a_1}=e_s^{n-p}v^K_{(s,k)}$.  Then $e_{-a_2}v_{-a_1}\neq 0$. Again $f_se_{-a_2}v_{-a_1}=0$, so  $a_1+a_2 \in \mathbb N$.

To show that  $e_{-a_2}e_sv_{(s,k)+1}^K\neq 0$ and the matching condition at $(s,k+1)$ holds, it suffices to show that $e_s^{q-1}e_{-a_2}e_s^{p+1}v_{-a_1}$ is non-zero and proportional to
 $e_s^{q}e_{-a_2}e_s^{p}v_{-a_1}$. The latter is the highest weight vector $v^K_{(s,k)+1}$.    Again $e_sv^K_{(s,k)+1}\neq 0$ by hypothesis, so $p+1\leq a_1$.  Then the required assertions follow as in \ref {4.7.1} using $(7)$.

%Choose $n \in \mathbb N$ maximal such that $v_{-a_1}:=f_s^nv^K_{(s,k)+1}\neq 0$.  By the hypothesis of the lemma, $v_{-a_1}$ is a non-zero multiple of $v^K_{(s,k)}$.
%
%One has
%$$f_s^{n+1}e_sv^K_{(s,k)+1}=-(n+1)f_s^n(h_s-n)v^K_{(s,k)+1}.\eqno {(22)}$$
%
%Now $v^K_{(s,k)+1}$ is an $h_s$ eigenvector.  We claim that its eigenvalue is distinct from $n$.  Otherwise $f_s^nv^K_{(s,k)+1}$ is a lowest weight vector of lowest possible weight in the $\mathfrak s$ module generated by $v^K_{(s,k)+1}$ and of $h_s$ eigenvalue $-n$.  This forces $v^K_{(s,k)+1}$ to be a highest weight vector contradicting the hypothesis that $e_sv_{(s,k)+1}^K\neq 0$.  (This argument would fail were it not for $f_s^nv^K_{(s,k)+1}$ being a lowest weight vector.) Thus the left hand side of $(22)$ is a non-zero multiple of $v_{-a_1}$.

\end {proof}

\subsubsection{}\label{4.7.6}

The proof of the minimax theorem is now complete.  If we regard $K_{\min}$ as corresponding to the minimal element of an $S$-set, then $K_{\max}$ corresponds to the maximal element of this $S$-set by \cite [5.7]{J4} and Theorem \ref {4.7.4}(i). 

 One may remark that it is possible to obtain further trails corresponding to elements of the $S$-set by inspection of the proof of Theorem \ref {4.7.4}.  Moreover from this we see an origin of their great redundancy in describing a maximum of the functions they define.  Indeed in the notation of Lemma \ref {4.7.4} we would expect the functions defined by the trails for which $\ell < u$ to be in the convex hull of those for which $\ell=0,u$ and hence be redundant.  This is hinted in subsection \ref {4.8} and will become much clearer following Section \ref {7}.

\subsection{}\label{4.8}

We shall take the minimax theorem as a paradigm of what we want from a theory of trails; but now replacing minimality by $\ell$-minimality.  A first step will be to obtain an analogue of Lemma \ref {4.7.2}.  For this we shall need the Chevalley-Serre relations.  To further proceed we shall need the full force of Lemma \ref {3.2} (and not just the case $r=2$).   Finally Lemma \ref {4.2} suggests that we should be able to compute the proportionality factors between the vectors lying in $T^+_s(\textbf{a})$.  For this we need to make a rather natural conjecture which we call the absence of ``false trails'' in $T^-_s(\textbf{a})$.  Then we shall need Lemma \ref {4.2}, for example in the proof of Corollary \ref {6.2}. Finally we shall show up to this conjecture that $T^+_s(\textbf{a})$ is the $\mathbb Z$ convex hull of the $S$-set it defines, which in turn through \cite [Thm. 1.4]{J5} form the set of its extremal elements.
%this will allow us to show that $Z_t$ can be taken to be the set of all trails (of type $t$).

We remark $T^-_s(\textbf{a})$ can be considered known by a previous induction step.  For the moment we have been unable to use this as the start of an induction procedure to prove our conjecture concerning false trails - see Section \ref {9}.

%It is an open question as to whether an $\ell$-minimal (resp. $\ell$-maximal) trail is a minimal (resp. maximal) trail.

\section{The Chevalley-Serre Relations and Their Consequences}\label{5}

 \subsection{}\label{5.1}

 Fix $s \in I$ and write $e_s$ simply as $e$. Consider a product $e_\mu$ of the simple root vectors $e_{s'}:s' \in I \setminus \{s\}$.  Set $a=-\alpha^\vee(\mu)$ which is a non-negative integer.  The Chevalley-Serre relations imply that
 $$(\ad e)^{a+1}e_{-a}: = \sum_{i=0}^{a+1}{a+1 \choose i}(-1)^{a+1-i}e^ie_{-a}e^{a+1-i}=0. \eqno {(22)} $$

 This relation can in principle be used to deduce that the $\overline{v}_\textbf{k} \in T_s(\textbf{a})$ are proportional in accordance with Lemma \ref {4.2}.  Our absence of false trails conjecture will eventually show (Lemma \ref {7.5.6}) that \textit{all} the scalar relations between the vectors defined by $T_s(\textbf{a})$ can be deduced from $(22)$.  Since $s \in I$ is arbitrary this determines all the scalar relations on extremal weight vectors defined by trails and in particular, up to our conjecture, \textit{describes all trails} - see Theorem \ref {7.6}.

 \subsection{}\label{5.2}

 Retain the above notation and hypotheses.    Let us write $e^{u_2}e_{-a}e^{u_1}:u_2,u_1 \in \mathbb N$ briefly as $<u_2,u_1>$.  We can view $<u_2,u_1>$ as the sub-expression in $\overline{v}_\textbf{k}$ at position $i+1\in \{1,2,\ldots,n-1\}$, by taking $a=a_{i+1},u_2=k_{i+1},u_1=k_i$.
 We write $<u_2,u_1>\approx 0$ to mean that the presence of this sub-expression in $\overline{v}_\textbf{k}$ makes the latter zero.

 Now given two such sub-expressions
 $<u_2,u_1>$ and $<u'_2,u'_1>$ at position $i+1$ with $u_2+u_1=u'_2+u'_1$, we write $c<u_2,u_1> \approx c'<u_2',u_1'>$ to mean that  $c\overline{v}_\textbf{k}=c'\overline{v}_{\textbf{k}'}$, given that the remaining exponents are unchanged, that is to say when $k'_j=k_j \in \mathbb N$, for all $j \in \{1,2,\ldots,r\}\setminus \{i+1,i\}$.  In this we say that these sub-expressions are quasi-equal.  We may also write the above as $c<u_2,u_1> - c'<u_2',u_1'> \approx 0$ and extend our convention to arbitrary finite sums.

 Let $p_2,p_1$ be non-negative integers such that $q:=p_2+p_1-a\geq 0$.

 \begin {lemma}  Fix $i \in \{1,2,\ldots,n-1\}$.  Assume $a\geq p_1$, so then $p_2-q \geq 0$. Suppose that $<u_2,u_1> \approx 0$ unless $u_2=p_2-v,u_1=p_1+v$, for $v \in \{0,1,2,\ldots,q\}$.

 Then the $\frac{<p_2-v,p_1+v>}{(p_2-v)!(p_1+v)! {q \choose v}}:v \in \{0,1,2,\ldots,q\}$, are quasi-equal.
 \end {lemma}

 \begin {proof} From $(22)$ we obtain on left multiplication by $e^u$ and right multiplication by $e^{q-u-1}$ that
 $$ \sum_{i=0}^{a+1}{a+1 \choose i}(-1)^{a+1-i}e^{i+u}e_{-a}e^{a-i+q-u}=0, \forall u \in \{0,1,\ldots,q-1\}. \eqno {(23)}$$

 Through the hypothesis of the lemma, this gives
 $$\sum_{v=0}^q \frac{(-1)^v}{(p_2-v-u)!(p_1-q+v+u+1)!} <p_2-v,p_1+v> \approx 0,  \forall u \in \{0,1,\ldots,q-1\}. \eqno {(24)}$$

 This is  system of a $q$ linear equations in $q+1$ unknowns which should properly be taken to be the appropriate $\overline{v}_{\textbf{k}'}$; but which we can simply take to be the \newline $\{<p_2-v,p_1+v>\}_{v=0}^q$ and replace $\approx$ by $=$.  Of course in every such term $\overline{v}_{\textbf{k}'}$, we require that $k'_j=k_j \in \mathbb N$, for all $j \in \{1,2,\ldots,r\}\setminus \{i+1,i\}$, as above.  

 This is close to a Toeplitz system, so probably its solution is known.  It can be shown to be unique up to scalars by showing that the determinant of the matrix with entries $\{\frac{(-1)^v}{(p_2-v-u)!(p_1-q+v+u+1)!}\}_{u,v=0}^{q-1}$ is non-zero.  As this paper is already long we leave this as an exercise to the reader.

 To show that our proposed solution is the correct one, it suffices to verify that
 $$\sum_{v=0}^q {q \choose v} \frac{(-1)^{q-v} (p_2-v)!(p_1+v)!}{(p_2-v-u)!(p_1-q+v+u+1)!}=0,  \forall u \in \{0,1,\ldots,q-1\}. \eqno {(25)}$$

 The proof of $(25)$ is standard.  Set $d_x=\partial/\partial x, d_y=\partial/\partial y$.  Then compute $d_x^{q-u-1}d_y^u$ on $x^{p_1}(x-y)^qy^{p_2-q}$ and evaluate at $x=y=1$.  This has an overall factor of $x-y$, so vanishes.   On the other hand it equals the left hand side of $(25)$.
 \end {proof}

 \textbf{Remark 1}. Consider $\overline{v}_\textbf{k}\in T_s(\textbf{a})$ and take $p_2=k_{i+1},a=a_{i+1},p_1=k_i$.  Then under the hypotheses of the lemma we can exactly shift up to $q_{i+1}:=k_{i+1}+k_i-a_{i+1}$ powers of $e$ to the right through $e_{-a_{i+1}}$ to obtain altogether $q_{i+1}$ elements of $T_s(\textbf{a})$.  Moreover the scale factors corresponding to the resulting elements $\overline{v}_\textbf{k} \in M_s(\textbf{a})$ are given by the conclusion of the lemma.

 \

 \textbf{Remark 2}.  As above consider a given sub-expression $e^{k_{i+1}}a_{i+1}e^{k_i}$ in a non-zero element $\overline{v}_\textbf{k}$ of $M_s(\textbf{a})$  and set $q_{i+1}=k_{i+1}+k_i-a_{i+1}$.  Let $E_i$ (resp. $E_{i+1}$) be a possible choice of the first (resp. second) exponent of $e$ for which this sub-expression again defines a non-zero element of $M_s(\textbf{a})$, with $E_i+E_{i+1}=k_i+k_{i+1}$.  Suppose that $k_i$ is a smallest possible value of $E_i$.    Then we may increase $E_i$ from $k_i$ to at least $k_i+q_{i+1}$ simultaneously decreasing $E_{i+1}$ from $k_{i+1}$ by the same amount.  Indeed if we could not increase $E_i$ strictly beyond $k_i+q_{i+1}$, then the lemma applies and the asserted conclusion also obtains.

 \subsection{}\label{5.3}

 Recall \ref {4.5} and let $\overline{v}_\textbf{l}$ denote the vector in $M_s(\textbf{a})$ corresponding to the unique $\ell$-minimal trail in $T_s^+(\textbf{a})$, more properly in $T_s^+(\textbf{e})$ - see  Convention \ref {4.2}.  Obviously this trail is determined by $\textbf{e}$, but not by $\textbf{a}$.  Both $\textbf{a}$ and $\textbf{l}$ are determined by $\textbf{e}$ and it is just $\textbf{a}$ and $\textbf{l}$ which will be important until the end of section \ref {7}.

%\textbf{ Notation.}  Fix a vector $\overline{v}_\textbf{k}$ of the form given in $(16)$.

 \subsubsection{}\label{5.3.1}

 Recall the notation of \ref {3.3}.

 \begin {lemma} One has $\ell_j\leq a_{j+1}$, for all $j\in N$.
 \end {lemma}

 \begin {proof}  %The proof is by induction on $j$.
  If not there exists $j\in N$ minimal such that $\ell_j > a_{j+1}$. Since $v_\textbf{l}$ is the unique $\ell$-minimal trail, one has $v_{\textbf{l}'}=0$, whenever $\ell'_i=\ell_i:i<j$ and $\ell'_j<\ell_j$.

 Now take $a=a_{j+1}$ in $(22)$ and multiply $(22)$ on the right by $e^{\ell_j-a_{j+1}-1}$.  Substitution into $v_\textbf{l}$ and using the vanishing of all the $v_{\textbf{l}'}$ above, implies that $v_\textbf{l}=0$, which is a contradiction.

 \end {proof}

  \subsubsection{}\label{5.3.2}

  Recall the notation of \ref {3.2}.

  \begin {lemma}  For all $j=1,2,\ldots,n$ one has $c_j:=a^{(j)}-\ell^{(j)}-\ell^{(j-1)}\geq 0$.
  \end {lemma}

  \begin {proof} The proof is by induction on $j$.  For $j=1$, the assertion is that $a_1\geq \ell_1$, which holds since $v_{-a_1}$ is a lowest weight vector.

  We remark that $e_{-a_2}e^{\ell_1}v_{-a_1}$ is not obviously a lowest weight vector, so this argument \textit{cannot} be repeated.  Indeed such a proof of the lemma would require knowing that $K_{\ell \min}=K_{\min}$.

  The continuation of the argument is rather complicated so as an illustration we consider first the case $j=2$.

  %One has $c_2-c_1=a_2-\ell_2-\ell_1$.  Thus if $c_2<0$, then
  Set $q_2:=\ell_2+\ell_1-a_2$. By Lemma \ref {5.3.1} we have $\ell_1\leq a_2$.  Taking $a=a_2, p_1=\ell_1$ recovers the first hypothesis of Lemma \ref {5.2}.

  Retain the notation of Remark 2 of \ref {5.2} and consider $E_1$.  It equals $\ell_1$.
    Since $K_{\ell \min}$ is minimal, $E_1$ cannot be decreased.  Then by Remark 2 of \ref {5.2}, we can increase $E_1$ by at least $q_2$. Yet since $v_{-a_1}$ is a lowest weight vector, non-vanishing forces $q_2+\ell_1 \leq a_1$, that is $c_1-q_2 \geq 0$. Yet  $c_2=c_1-q_2$, so the required assertion obtains.

  The general case obtains by repetition of this argument; but with a little extra difficulty.

  Assume $i \geq 2$. Recall that $c_1\geq 0$ and let $i\in N$ be minimal such that let $c_{i+1}<0$.   With this choice $q_{i+1}:=\ell_{i+1}+\ell_i-a_{i+1} = c_i-c_{i+1}>0$. By Lemma \ref {5.3.1} we have $\ell_i\leq a_{i+1}$.  Taking  $a=a_{i+1},p_1=\ell_i$ recovers the first hypothesis of Lemma \ref {5.2}.

  Consider $E_{i}$. It equals $\ell_i$. By the minimality of $K_{\ell \min}$ it cannot be decreased.  Then by Remark 2 of \ref {5.2} we can increase $E_i$ to at least $q_{i+1}+\ell_i=2\ell_i+\ell_{i+1}-a_{i+1}$.

  We claim that for all positive integer $j\leq i$ we can take $E_j$ to be at least $b_j:=2\ell_j+2\ell_{j+1}+\cdots+ 2\ell_i +\ell_{i+1}-(a_{j+1}+\cdots+a_{i+1})$, by successively increasing $\ell_i,\ell_{i-1},\ldots,\ell_j$ \textit{in that order}. This claim is proved by decreasing induction on $j$.  It has already been established for $j=i$.

  Suppose $b_j+\ell_{j-1}-a_j>0$.   Since $E_{j'}:j'\leq j-1$ has not yet been altered and so still equals $\ell_{j'}$, the minimality of $K_{\ell \min}$ implies that $E_{j-1}$ cannot be decreased. Then by Remark 2 of \ref {5.2}, we may increase $E_{j-1}=\ell_{j-1}$ to at least $(b_j+\ell_{j-1}-a_j)+\ell_{j-1}=b_{j-1}$, as required.

  Suppose that $b_j+\ell_{j-1}-a_j\leq 0$. Then $b_{j-1}\leq \ell_{j-1}$ and there is nothing to prove. This establishes the claim.

    Finally since $v_{-a_1}$ is a lowest weight vector, we obtain $b_1\leq a_1$, whilst $c_{i+1}=a_1-b_1$, so the required assertion obtains.

  \end {proof}

   \subsubsection{}\label{5.3.3}

  In analogy with \ref {4.7.2} we set $c^{K_{\ell \min}}_j(s):=-\alpha^\vee_s(\delta^{J_{\ell \min}}_{(s,j)}):j \in \mathbb N^+$. It is immediate from definitions that
  $$c^{K_{\ell \min}}_j(s)=a^{(j)}-(\ell^{(j)}+\ell^{(j-1)}). \eqno{(26)}$$

  \begin {cor} $c^{K_{\ell \min}}_j(s)\geq 0$, for all $j \in \mathbb N^+$.
  \end {cor}

  \begin {proof}

  It is immediate from definitions that $c^{K_{\ell \min}}_k(s)=0$, for all $k>n$.  For $k\leq n$, the assertion follows from Lemma \ref {5.3.2}.
  \end {proof}

  \subsubsection{}\label{5.3.4}

   Let $\overline{v}_\textbf{l}$, be the vector corresponding to the $\ell$-minimal trail.  Define $\overline{\textbf{l}}$ by setting $\ell_n=0$.
   % Obviously $\overline{\textbf{l}}$.

   In analogy with the first observation in \ref {4.7.2} we prove that

  \begin {lemma}

  \

  (i) $\overline{v}_{\overline{\textbf{l}}}$ is a non-zero multiple of $v_{-s_\alpha w\varpi_t}$.

  \

  (ii) $c^{K_{\ell \min}}_n(s)=0$.

  \end {lemma}

  \begin {proof}
 Choose $b \in \mathbb N$ maximal such that $f^b\overline{v}_{\overline{\textbf{l}}}\neq 0$. It is a lowest weight vector in $M_s(\textbf{a})$, hence proportional to $v_{-s_\alpha w\varpi_t}$. We can write it as $\overline{v}_{\textbf{l}'}$. Set  $\ell'_n=\alpha^\vee(-w\varpi_t)$ and $\textbf{k}=(\ell',\textbf{l}')$. Then $\overline{v}_\textbf{k}$ is a non-zero multiple of $v_{- w\varpi_t}$ and hence a non-zero multiple of $\overline{v}_\textbf{l}$.  On the other hand unless $b=0$, we obtain $\textbf{k}<\textbf{l}$ for the lexicographic order.  Hence (i).  Finally (ii) follows from (i).
  \end {proof}

  \textbf{Remark}.  Note that $\textbf{l}$ and $\overline{\textbf{l}}$ only differ in the last factor which is $\ell_n$ in the first case and $0$ in the second case. It is the action of $e^{l_n}$ which takes the lowest weight vector of $M_s(\textbf{a})$ to its highest weight vector.

  \subsubsection{}\label{5.3.5}

  Thus we have finally proved with the help of the Chevalley-Serre relations that the unique $\ell$-minimal trail in $T_s(\textbf{a})$ provides a suitable set of parameters for an $S$-set.

  We set $\textbf{c}:=\{c^{K_{\ell \min}}_k(s)\}_{k=1}^{n-1}$.  By Lemma \ref {5.3.2}, $\textbf{c}$ is determined by $\textbf{a}$ and by $\textbf{l}$, hence by $\textbf{e}$, via \ref {5.3}.

   \subsection{}\label{5.4}

   The $S$-set of type $s$ defined by the coefficients $\{c^{K_{\ell \min}}_i(s)\}_{i=1}^n$ must be shown to give rise to new trails $K \in T_s(\textbf{a})$.  If such a trail $K$ exists, then it is completely determined by its function $z^K$ defined in $(2)$.

   %To describe this we write $c^{K_{\ell \min}}_i(s)$, simply as, $c_i$.

   If $K$ is the trail defined by $\overline{v}_\textbf{k}$, one obtains
    $$c_i^K(s):=-\alpha_s^\vee(\delta^K_{(s,i)})=a^{(i)}-(k^{(i)}+k^{(i-1)}), \eqno{(27)}$$
    in the notation of \ref {3.2}.  It is immediate that $c^K_i(s)=0$, for all $i >n$.

   On the other hand an $S$-set defines a set of parameters $\{c'^K_i(s)\}_{i=1}^n$ which are certain specific linear combinations of the $\{c^{K_{\ell \min}}_i(s)\}_{i=1}^n$.  Then one may define
   $$z^K= \sum_{i=1}^n c'^K_i(s)(r_s^i-r_s^{i+1})-\sum_{i=1}^n c^{K_{\ell \min}}_i(s)m_s^i. \eqno {(28)}$$

   In this recall that
   %where for all $i \in \mathbb N^+$, the $m_s^i$ are the co-ordinate functions on $B_J$  and $r_s^i$ are the Kashiwara functions pertaining to $s \in I$, defined in \ref {2.3}.  One has
   $$r_s^i-r_s^{i+1}=m_s^i+m_s^{i+1}  \mod \sum_{s' \in I\setminus \{s\}, i \in \mathbb N^+} \mathbb Zm_{s'}^i. \eqno {(29)}$$

   Thus in order to recover $(27)$ from $(28),(29)$, we must take
   $$c'^K_i(s)=k^{(i)}-\ell^{(i)}, \forall i=1,2,\ldots,n, \eqno {(30)}$$

More generally we view $(30)$ as defining a set of parameters $c'^K_i(s):i \in \hat{N}$, or simply $c'_i: i \in \hat{N}$, relating a vector $\overline{v}_\textbf{k}\in M_s(\textbf{a})$ to the vector $\overline{v}_\textbf{l}$ associated to the $\ell$-minimal trail.

  \section{The Rigid Case}\label{6}

  Fix $s \in I$ and retain the notation of \ref {5.3}.  Recall in particular the definition of $\textbf{l}$ and of $\overline{\textbf{l}}$.

  Suppose that
  $$a_{j+1}-\ell_{j+1}-\ell_j\geq 0, \forall j =1,2,\ldots,n-2, \eqno {(31)}$$
  equivalently that the $c_j^{K_\ell \min}(s):j=1,2,\ldots,n-1$ are increasing.

  In this case, Lemma \ref {5.2} does not lead to any new elements of lowest weight in $M_s(\textbf{a})$ of the form $\overline{v}_{\overline{\textbf{k}}}$.  For this reason we say that $(31)$ defines the rigid case.

  Here we shall analyze the rigid case under the \textit{conjecture} that $\overline{v}_{\overline{\textbf{l}}}$ is indeed the unique vector obtained from a trail in $T^-_s(\textbf{a})$.
  %, that is to say a lowest weight vector of  $M_s(\textbf{a})$.
  However our first result described in \ref {6.1} will \textit{not} need this hypothesis.

   \subsection{}\label{6.1}

   Let $v_\textbf{k},v_\textbf{l}$ denote elements in the $n$-fold tensor product $V(-a_n)\otimes V(-a_{n-1})\otimes \cdots \otimes V(-a_1)$ defined as in \ref {3.1}.

    As in \ref {3.2} we assume that $b_j:=k_j-\ell_j\geq 0$, for all $i=1,2,\ldots,n$ and let $b$ denote their sum.  Set $\ell^{(0)}=0$.

   \begin {lemma}

   \

   (i) $a^{(j)}-k^{(j)}-\ell^{(j-1)} \geq 0$, for all $j=1,2,\ldots,n$ implies that the coefficient of $v_\textbf{l}$ in $f^bv_\textbf{k}$ is non-zero.
   % $\overline{v}_\textbf{k}\neq 0$.

   \

   Suppose that $(31)$ holds.  Suppose further that the coefficient of $v_\textbf{l}$ in $f^bv_\textbf{k}$ is non-zero.
   %so in particular $\overline{v}_\textbf{k}\neq 0$.
   Then

   \

   (ii) $a^{(j)}-k^{(j)}-\ell^{(j-1)} \geq 0$, for all $j=1,2,\ldots,n$.

   \

   (iii) Each linear factor in $A_\textbf{b}(\textbf{k},\textbf{l})$ is positive. In particular $A_\textbf{b}(\textbf{k},\textbf{l})$ is positive.

   \end {lemma}

   \begin {proof} Recall that $b_j+k^{(j-1)}+\ell^{(j)}=k^{(j)}+\ell^{(j-1)}$.  Thus the hypothesis of (i) implies that every linear factor in $A_\textbf{b}(\textbf{k},\textbf{l})$ is positive.  Consequently by Lemma \ref {3.2} it follows that the coefficient of $v_\textbf{l}$ in $f^bv_\textbf{k}$ is positive.  Hence (i).

   For (ii) we first show that
   $$a^{(j)}-k^{(j)}-\ell^{(j-1)} \geq 0, \eqno {(32)}$$
   by induction on $j$.

   For $j=1$ it is just the assertion that $a_1\geq k_1$ which follows from the fact that $\overline{v}_\textbf{k}\neq 0$ and that $v_{a_1}$ is a lowest weight vector.

   Combining  $(31)$ and $(32)$ we obtain
   $$a^{(j+1)}-k^{(j)}-\ell^{(j+1)} \geq 0, \eqno {(33)}$$

   On the other hand the factors $a^{(j+1)}+1-i-k^{(j)}-\ell^{(j+1)}$ occurring in $A_\textbf{b}(\textbf{k},\textbf{l})$ decrease in $i$, are integer and for $i=1$ is non-negative by $(33)$.  Then by the hypothesis of (ii) they must be positive for all $i \in [1,b_{j+1}]$.   Taking $i=b_{j+1}$ recovers $(32)$ with $j$ increased by $1$.  This gives (ii).  Then (iii) obtains from (ii).

    \end {proof}

   %\textbf{Remark}.  (iii) can fail in the non-rigid case. Amazingly this can be compensated by the fact that $f^n \overline{v}_{\textbf{k}}$ is in general a sum of terms in the $\overline{v}_{\textbf{l}'} \in T_s(\textbf{a})$ which are all proportional to $\overline{v}_\textbf{l}$, so the overall sum may be positive.

   \subsection{}\label{6.2}

  Continue to assume that $(31)$ holds.  Assume our conjecture that this implies $T^-_s(\textbf{a})$ to be reduced to just one trail, necessarily the $\ell$-minimal one by Lemma \ref {5.3.4}, and let $\overline{v}_{\overline{\textbf{l}}}$ denote the corresponding lowest weight vector in $M_s(\textbf{a})$.  Define $b \in \mathbb N$ as in \ref {6.1}.

   \begin {cor}
   Under these assumptions,
$\overline{v}_\textbf{k}\neq 0$ if and only if $a^{(j)}-k^{(j)}-\ell^{(j-1)} \geq 0$, for all $j=1,2,\ldots,n$ and $k_j\geq \ell_j$, for all $j=1,2,\ldots,n-1$.
   \end {cor}

   \begin {proof} The simplicity of $M_s(\textbf{a})$ implies that $\overline{v}_{\overline{\textbf{l}}}$ is the unique up to scalars lowest weight vector of $M_s(\textbf{a})$. Hence $\overline{v}_\textbf{k}\neq 0$ if and only if $f^b\overline{v}_\textbf{k}$ is a non-zero multiple of $\overline{v}_{\overline{\textbf{l}}}$.  Yet $f^bv_\textbf{k}$ is a sum of terms of the form $v_{\overline{\textbf{l}}'}$ and by our conjecture their images $\overline{v}_{\overline{\textbf{l}}'}$ in $M_s(\textbf{a})$ are zero unless $\textbf{l}'=\textbf{l}$. Thus  $\overline{v}_\textbf{k}\neq 0$ if and only if the coefficient of $v_{\overline{l}}$ in $f^bv_\textbf{k}$ is non-zero. Then apply Lemma \ref {6.1} with $\textbf{l}$ replaced by $\overline{\textbf{l}}$, recalling that $\overline{\ell}_n=0$.
   \end {proof}

 \subsection{}\label{6.3}

 %Our conjecture that in the rigid case (equivalently increasing coefficient case) can be expressed by saying that $T_s^-(\textbf{a})$ is reduced to a single element. Then
In the rigid case and under the assumption $|T^-_s(\textbf{a})|=1$, Corollary \ref {6.2} describes all the trails in $T_s(\textbf{a})$.

 In this next subsection we show the above set of trails is exactly the $\mathbb Z$ convex hull of the elements of the $S$-set defined by $\{c_j^{K_{\ell \min}}(s)\}_{j=1}^n$.

  \subsection{}\label{6.4}

  Observe that
  $$c^{K_{\ell \min}}_j(s)-c'^K_j(s)=a^{(j)}-k^{(j)}-\ell^{(j-1)},\forall j=1,2,\ldots,n. \eqno {(34)}$$
  $$ c'^K_j(s)-c'^K_{j-1}(s)=k_j-\ell_j, \forall j=1,2,\ldots,n.\eqno{(35)}$$

  Thus the conclusion of Corollary \ref {6.2} can be expressed by saying that the $c_j'^K(s)$ are increasing and $c^{K_{\ell \min}}_j(s)\geq c'^K_j(s)$ for all $j=1,2,\ldots,n$.  On the other hand the $c'^K_j(s)$ are clearly integers. In view of \cite [Thm. 1.4] {J5} we obtain

  \begin {thm} Assume that the $\{c_j^{K_{\ell \min}}(s)\}_{j=1}^n$ are increasing and that  $T^-_s(\textbf{a})$ is reduced to just one trail.  Then the $z^K:K\in T_s(\textbf{a})$ consists of the integer points of the convex set whose extremal elements form the $S$-set defined by $\{c_j^{K_{\ell \min}}(s)\}_{j=1}^n$.
  \end {thm}

 \subsection{}\label{6.5}

 The unique up to a scalar solution to $(23)$ given by Lemma \ref {5.2} was not obtained by using Cramer's rule, nor by an inspired guess and least of all by a literature search.  Indeed it was obtained by a disarmingly simple method which will be useful in the sequel.

 As in \ref {4.2}, let $v_{-a_3},v_{-a_2},v_{-a_1}$ be lowest weight vectors for the $\mathfrak {sl}(2)$ subalgebra $\mathfrak s$ of $\mathfrak g$ defined in \ref {4.2}.

 Let $\varphi:v_{-a_3}\otimes e^{u_2}v_{-a_2}\otimes e^{u_1} v_{-a_1}\mapsto e_{-a_3}e^{u_2}e_{-a_2}e^{u_1} v_{-a_1}$ of $V_{-a_3}\otimes V(-a_2)\otimes V(-a_1)$ in $M_s(\textbf{a})$ be the $\mathfrak s$ module map noted in \ref {4.3} (for $r=3$).

 % As in \ref {4.2}, take the adjoint action of $\mathfrak s$ on $e_{-a_2},e_{-a_2}$  and consider $e_{-a_3}e^{u_2}e_{-a_2}e^{u_1} v_{-a_1}$ as an image under an $\mathfrak s$ module map $\varphi$ of $v_{-a_3}\otimes e^{u_2}v_{-a_2}\otimes e^{u_1} v_{-a_1}\in V_{-a_3}\otimes V(-a_2)\otimes V(-a_1)$.

 Fix $p_2,p_1 \in \mathbb N$ and assume $q:=p_2+p_1-a_2\geq 0$.

 Now take $u_2,u_1 \in \mathbb N$ with $u_2 \geq p_2-q,u_1 \geq p_1$.
 %Assume that $(u_2-(p_2-q))+(u_1-p_1)$ is fixed.
 By Lemma \ref {3.2} we may calculate the coefficient $A(q,v,p_2,p_1,a_1)$ of $v_{-a_3} \otimes e^{p_2-q}v_{-a_2}\otimes e^{p_1} v_{-a_1}$ in $f^q v_{-a_3}\otimes e^{p_2-v}v_{-a_2}\otimes e^{p_1+v} v_{-a_1}:v \in [0,q]$, by taking $n=3,k_3=\ell_3=0,b=q,k_2=p_2-v,k_1=p_1+v,\ell_2=p_2-q,\ell_1=p_1$ in $(7)$.

 Noting that $a_2+a_1-k_1-(\ell_2+\ell_1)=a_1-p_1-v$, the products for $i=1,2$ in $(7)$ coalesce into a single product and we obtain

$$A(q,v,p_2,p_1,a_1)=q!{p_2-v \choose p_2-q}{p_1+v \choose p_1}\prod_{i=1}^q(a_1+1-p_1-i).$$

Thus up to a common factor we can simply write
 $$A(q,v,p_2,p_1,a_1)=(p_2-v)!(p_1+v)!{q \choose v}, \forall v \in [0,q].\eqno {(36)}$$

 Now the conclusion of Lemma \ref {5.2} shows that these are exactly the proportionality factors between the $\varphi(v_{-a_3}\otimes e^{p_2-v}v_{-a_2}\otimes e^{p_1+v} v_{-a_1})=e_{-a_3}e^{p_2-v}e_{-a_2}e^{p_1+v}v_{-a_1}:v \in [0,q]$ when we use the Chevalley-Serre relations and make the additional hypothesis (of Lemma \ref {5.2}) that these expressions vanish for $v \notin [0,q]$.  Moreover these terms must either be all non-zero or all zero.  Since $v_{-a_1}$ is a lowest weight vector, the former forces $a_1 \geq q+p_1$, which in turn implies that the common factor we have eliminated above to be non-zero.

Actually this last fact can be better seen as follows.  Eventually we will identify $e_{-a_3}e^{p_2}e_{-a_2}e^{p_1}v_{-a_1}$ with the lowest weight vector $\overline{v}_{\overline{\textbf{l}}}$.  In this, $(26)$ gives $c_1=a_1-p_1,c_1-c_2=q$.  Then by Lemma \ref {5.3.2} we obtain $a_1-p_1\geq q$, which forces this common factor to be non-zero.

Notice we are \textit{not} claiming that this gives a further proof of Lemma \ref {5.2}.  Yet it has already proved its usefulness in guessing a solution to $(24)$ and will further prove invaluable in computing sums resulting from $f^b\overline{v}_\textbf{k}$.  Here we shall need a slight generalization of the above.  It is given in \ref {7.5.5}.

 \section{The General Case}\label{7}

 Here we drop the condition of rigidity. Recall the notation of \ref {4.2} and \ref {4.3}.  The aim of this section is to use the $\mathfrak s$ module structure of $M_s(\textbf{a})$ to deduce the structure of $T^+_s(\textbf{a})$ from that of $T^-_s(\textbf{a})$, as we have already done in the rigid case.  Ultimately a comparable result (Theorem \ref {7.6}) is obtained but the proof is far more technical.
\subsection{}\label{7.1}

Recall the notation of \ref {3.3} and \ref {5.3.3}. Write $c^{K_{\ell \min}}_j(s)$ simply as $c_j$, for all $j \in N$.  Recall that the $c_j:j \in N$ are non-negative integers (Lemma \ref {5.3.2}).

 As in \cite [1.4]{J5}, we lift the natural order on $\textbf{c}:=\{c_i\}_{i \in N}$ to a total order which we view as a total order $\prec$ on $N$.  Choose a labelling $u_i:i \in N$ such that $u_1\prec u_2 \prec  \ldots \prec u_{r-1}$.  Then for all $j\in N$ relabel the subset $N_j:=\{u_i\}_{i=1}^j$ of $N$, as $\{v_i\}_{i=1}^j$  with the $v_i$ increasing in the natural order.  Let $\theta$ be the bijection of $N$ defined by $u_j=v_{\theta(j)}$, for all $j \in N$.

 Define a convex set $K(\textbf{c}):=\{c'_j\}_{j\in N} \subset \mathbb Q^{r-1}$ by
 $$0\leq c_j'\leq c_j, \eqno{(37)}$$
 $$c'_{v_{\theta(j)+1}}-c'_{v_{\theta (j)}} \geq -(c_{v_{\theta (j)}}-c_{v_{\theta (j)+1}}), \quad c'_{v_{\theta(j)}}-c'_{v_{\theta (j)-1}} \geq 0, \eqno {(38)}$$
 for all $j \in N$.

 One may remark that if the $\{c_j\}_{j \in N}$ are increasing, then $\theta$ is the identity map and so $(38)$ just means that the $\{c'_j\}_{j \in N}$ are increasing.  This is the rigid case which is much easier to handle.

 A main result of \cite {J5} is that the canonical $S$-set $Z(\textbf{c})$  defined by $\textbf{c}$ (as the functions attached to the vertices of the corresponding canonical $S$-graph) consists of the extremal points of $K(\textbf{c})$.  Although one can compute these functions from \cite {J4}, \cite {JL}, their presentation as the extremal elements of $K(\textbf{c})$ is more convenient and will be that mainly used here.

 Now let $K_\mathbb Z(\textbf{c})$ denote the set of integer points of $K(\textbf{c})$, that is to say $K_\mathbb Z(\textbf{c}):=K(\textbf{c})\cap \mathbb Z^{n-1}$.

\subsection{}\label{7.2}

In order to prove the key proposition below we need to recall the properties of $S$-sets and the construction of the canonical $S$-sets.

\subsubsection{}\label{7.2.1}

Let $\textbf{c}:=\{c_i\}_{i \in N}$ be the set of coefficients viewed as being a totally ordered set by taking a lift of the natural order on the coefficients.  The results below will not depend on the lift.  This is because, although the graphs differ, the functions defined by their vertices still coincide \cite [5.8]{JL}.

An $S$-graph is a graph $\mathscr G$ with a number of properties.  These are listed in \cite [Sect. 6]{J4}. In this we note the following.

Let $V(\mathscr G)$ (resp. $E(\mathscr G)$) denote the set of vertices (resp. edges) of $\mathscr G$.

There is map from $E(\mathscr G) \rightarrow N$ (resp. $V(\mathscr G)\rightarrow\hat{N}$).  For all $i \in \hat{N}$, let $V^i(\mathscr G)$ denote its inverse image in $V(\mathscr G)$.

Given a vertex $v \in V^k(\mathscr G)$, we write $i_v=k$.  Given neighbouring vertices $v,v'$, we let $(v,v')$ be the (unique) edge joining them and let $i_{(v,v')}$ denote its image in $N$.  We write $c_{i_{(v,v')}}$ simply as $c_{(v,v')}$.

Recall the co-ordinate (resp. Kashiwara) functions $\{m^i_s\}_{s \in I,i \in \mathbb N^+}$ (resp. $\{r^i_s\}_{s \in I,i \in \mathbb N^+}$) on $B_J$ defined in \cite [2.3.2]{J3}.  They are related by $(29)$ which is all we need here. Let $Z(B_J)$ denote the free $\mathbb Z$ module generated by the $\{m^i_s\}_{s \in I,i \in \mathbb N^+}$. A given $S$-set can be viewed as being of fixed type $s \in I$ and we then drop the subscript $s$.

%For all $i \in \mathbb N^+,s \in I$, let $r^i_s$ be the linear functions on $B_J$ defined in \cite [2.3.2]{J3}.  There are certain linear combinations of the co-ordinate functions $m_s^i:i \in \mathbb N^+,s \in I$ on $B_J$. (Here the exponent does \textit{not} signify a power.)

%Let $R$ denote the linear space with basis $r^i: i \in \mathbb N^+$.  In this the exponent does \textit{not} signify a power. Eventually the $r^i$ are viewed as Kashiwara functions and a subscript $s \in I$ is added to indicate that they are of type $s \in I$.  Specifically these functions are the linear functions on $B_J$ defined in \cite [2.3.2]{J3}.

There is a map $v \mapsto f_v$ of $V(\mathscr G) \rightarrow Z(B_J)$.  It has the property that
$$f_v-f_{v'}=c_{v,v'}(r^{i_v}-r^{i_{v'}}), \eqno {(39)}$$
for any pair of adjacent vertices $v,v'$.

An $S$-graph admits a unique chain (called a pointed chain \cite [6.3]{J4}) of adjacent vertices $v_j\in V^j(\mathscr G):j \in \hat{N}$ such that $i_{(v_{j+1},v_j)}=j$, for all $j \in N$. %Let $v^j:j=1,2,\ldots,r$ be the unique vertex in this chain lying in $V^j(\mathscr G)$.

The function $f_{v_n}$
%, where $v_r$ is the unique element in the above chain with label $r$,
is called the driving function of the $S$-set.  We assign to it $z_{K_{\ell \min}}$, where $K_{\ell \min}$ is the unique $\ell$-minimal element of $T_s(\textbf{a})$.   We may write (using the notational convention of \ref {7.1})
$$f_{v_n}= -\sum_{k=1}^n c_k m_s^k  \mod \sum_{s' \in I\setminus \{s\}, k \in \mathbb N^+} \mathbb Zm_{s'}^k. \eqno {(40)}$$

Here the terms in $\sum_{s' \in I\setminus \{s\}, k \in \mathbb N^+} \mathbb Zm_{s'}^k$, obtained by an induction procedure on Weyl group element lengths, can be ignored and the $s$ subscript dropped.

Since an $S$-graph is connected, $(39),(40)$ determine the set of functions $\{f_v\}_{ v \in V(\mathscr G)}$.  It is called the $S$-set determined by the $S$-graph $\mathscr G$.

An $S$-graph is \textit{not} unique.  However for each $\textbf{c}$ as above, there is a specific $S$-graph $\mathscr G(\textbf{c})$ constructed by ``binary fusion'' \cite [7.2]{J4}.  It is canonical in the sense of \cite [5.6]{JL}.  We recall briefly this construction.

When $n=1$, we take $\mathscr G(\textbf{c})$ to be the graph with one vertex and having label $1$.

Let $u$ be the unique maximal element of $N$ for the order relation $\prec$.

 %and recall the process of binary fusion \cite [7.2]{J4}, the details of which we shall repeat briefly here.

  Set $\textbf{c}^-=\textbf{c}\setminus \{c_u\}$.  Then the graph $\mathscr G(\textbf{c}^-)$ involving one less coefficient can be assumed to be defined.

We define new graphs $\mathscr G^\pm$ isomorphic to $\mathscr G(\textbf{c}^-)$ as unlabelled graphs.

 Let $\mathscr G^+$ be the graph obtained from $\mathscr G(\textbf{c}^-)$ by leaving the labels in $[1,u-1]$ unchanged and increasing the labels in $[u,r-1]$ by $1$. Let the labelling on $\mathscr G^-$ be defined by required that the above defined unlabelled graph isomorphism $\varphi:\mathscr G^+ \iso \mathscr G^-$ fixes all labels with the following exception. If $v \in V(\mathscr G^+)$ has label $u+1$, then $\varphi(v) \in \mathscr G^-$ is assigned label $u$.

Notice that if $v \in V(\mathscr G(\textbf{c}^-)$ has label $i\in N$, then the corresponding vertex $v \in \mathscr G^\pm$ has label $\phi^\pm(i)\in \hat{N}$, where
 %$$ \phi^\pm(i)=\left\{
%                       \begin{array}{ll}
%                         i, & \hbox{if $i<u+\frac{1}{2}\mp \frac{1}{2}$} \\
%                         i+1, & \hbox{if $i\geq u+\frac{1}{2}\mp \frac{1}{2}$,}
%                       \end{array}
%                     \right. $$

$$ \phi^+(i)=\left\{
                       \begin{array}{ll}
                         i, & \hbox{if $i<u$} \\
                         i+1, & \hbox{if $i\geq u$,}
                       \end{array}
                     \right. $$

$$ \phi^-(i)=\left\{
                       \begin{array}{ll}
                         i, & \hbox{if $i\leq u$} \\
                         i+1, & \hbox{if $i> u$,}
                       \end{array}
                     \right. $$
for all $i \in N$.
 %One may observe that $V^u(\mathscr G^+)=\phi$ whilst $V^{u+1}(\mathscr G^-)=\phi$.

One may remark that $\phi^+(k)=\phi^-(k)$, for all $k \in N\setminus \{u\}$, whilst $\phi^+(u)=u+1,\phi^-(u)=u$.   As a consequence $V^u(\mathscr G^+)=\phi$ whilst $V^{u+1}(\mathscr G^-)=\phi$.

  Then $\mathscr G(\textbf{c})$ is defined as the union of $\mathscr G^+$ and $\mathscr G^-$ in which each vertex $v$ of $\mathscr G^+$ with label $u+1$ is joined to $\varphi(v)$, so having label $u$, with an edge having label $u$.

  Through $\varphi$ we may view $\phi^+$ (resp. $\phi^-$) as a labelled graph embedding of $\mathscr G^+$ (resp. $\mathscr G^-$) into $\mathscr G(\textbf{c})$.

  The graphs $\mathscr G(\textbf{c})$ have some extra properties noted below.

  \

 \textbf{Property 1}.  For all $v \in V^k(\mathscr G(\textbf{c}))$ the coefficient of $m^k$ in $f_v$ is zero \cite [Lemma 5.4 and $P_8$ of 6.7]{J4}. (A possibly more transparent proof obtains from \cite [4.3]{J5}.) If the $c_i>0:i \in N$, then the converse also holds and results from the manner in which the graphs degenerate as described in \cite [Sect. 5.8.2]{JL}.   Again if $c_{n-1}\neq 0$, then the converse holds for $k=n$, that is to say the coefficient of $m^n$ in $f_v$ is zero only if $v \in V^n(\mathscr G(\textbf{c}))$.  This follows by also using \cite [Lemma 5.8.1]{JL}.

 \

 %One can ask if this holds for all $S$-graphs (by combining $(29),(39),(40)$. For example this is easily checked if $v$ lies in the pointed chain.

 \textbf{Property 2.} The $S$-set $Z(\textbf{c})$ determined by $\mathscr G(\textbf{c})$, admits the convexity property described in \cite [Thm. 1.4]{J5} as noted in \ref {7.1}.  Here the difference of successive Kashiwara functions $r^i-r^{i+1}:i \in N$ is viewed as the $i^{th}$ co-ordinate function on $\mathbb Q^{|N|}$.  Moreover the driving function can be taken equal to zero, since convexity does not depend on a choice of origin.   This fact is crucial in proving that the dual Kashiwara parameter $\varepsilon^\star_t$ is given (in the absence of false trails) by the maximum of the functions defined by the elements of $\mathscr K^{BZ}_t$ (see Remarks following Theorem \ref {8.7}).

%A third property special to $\mathscr G(\textbf{c})$ is obtained by induction using binary fusion.  The proof is more tricky than one might have anticipated due to changes in labelling.
%
%Fix $j \in N$ and let $\mathscr G^j(\textbf{c})$ be the graph with vertices $V^j(\mathscr G(\textbf{c}))$ and edges joining pairs $v_1,v_2$ given $f_{v_1}-f_{v_2}=(c_{u_1}-c_{u_2})(r^{u_1+1}-r^{u_1})$, for some $u_1,u_2 \in N$ distinct and $j \notin \{u_1,u_1+1\}$.  When this holds we say that $v_1,v_2$ are neighbours in $\mathscr G^j(\textbf{c})$, through $u_1$.

\

\textbf{Remark}.  The functions $\{f_v\}_{v \in V(\mathscr G^+)}$ defined by the vertices of the subgraph $\mathscr G^+$ of $\mathscr G(\textbf{c})$ are independent of $c_u$.  Moreover they are obtained from the functions $\{f_v\}_{v \in V(\mathscr G(\textbf{c}^-))}$ by the relabelling prescribed by $\phi^+$. This means that the former are the \textit{same} functions as the latter defined with respect to the coefficient set $\{c_1,c_2,c_{u-1},c_{u+1},\ldots,c_{n-1}\}$ \textit{but} by viewing $r^{\phi^+(i)}-r^{\phi^+(i+1)}$ as the $i^{th}$ co-ordinate function.  This statement, which is a a trivial consequence of definitions, is embodied in \cite [Eq. $(12)$]{J5}.  Thus we may identify the functions in these two sets.

Similarly the functions $\{f_v\}_{v \in V(\mathscr G^-)}$ coincide with those given by $\{f_v\}_{v \in V(\mathscr G(\textbf{c}^-))}$ defined with respect to the coefficient set $\{c_1,c_2,c_{u-1},c_{u+1},\ldots,c_{n-1}\}$ \textit{but} by viewing $r^{\phi^-(i)}-r^{\phi^-(i+1)}$ as the $i^{th}$ co-ordinate function  \textit{and up to} adding the fixed term $(c_u-c_{u+1})(r^u-r^{u+1})$.  The first part, which is a trivial consequence of definitions, is embodied in \cite [Eq. $(13)$]{J5}.  The second part, which is less trivial, is just \cite [Eq. $(14)$]{J5}. Thus up translation by the fixed term $(c_u-c_{u+1})(r^u-r^{u+1})$ we may identify the functions in these two sets.

\subsubsection{\textbf{Property 3}}\label{7.2.2}

A third property special to $\mathscr G(\textbf{c})$ is described by the Lemma below.  It is obtained by induction using binary fusion.  The proof is more delicate than one might have anticipated due to changes in labelling.

\

\textbf{Definition}. Fix $j \in N$ and let $\mathscr G^j(\textbf{c})$ be the graph with vertex set $V^j(\mathscr G(\textbf{c}))$ whose edges are obtained by joining pairs $v_1,v_2$ satisfying $$f_{v_1}-f_{v_2}=(c_{u_1}-c_{u_2})(r^{u_1+1}-r^{u_1}), \text{with} \ u_1,u_2 \in N \text{ distinct and} j \notin \{u_1,u_1+1\}.$$
 When this holds we say that $v_1,v_2$ are neighbours in $\mathscr G^j(\textbf{c})$, through $u_1$.

\

Let $(\mathscr G^\pm)^j$ be the full subgraphs of $\mathscr G^j(\textbf{c})$ with vertices in $\mathscr G^\pm$.

\begin {lemma} For all $j \in N$, the graph $\mathscr G^j(\textbf{c})$ is connected.
\end {lemma}

\begin {proof} The proof is by induction on $|N|$.  When $|N|=1$, $\mathscr G(\textbf{c})$ has just two vertices with different labels so there is nothing to prove.

Take $u \in N$ to be the unique maximal element for $\prec$.

\

  Assume first that $j\in N\setminus \{u,u+1\}$.

\

 \textbf{ Step 1.}

  Consider $v\in V^j(\mathscr G^+)$.  By our assumption on $j$ one has $\varphi(v) \in V^j(\mathscr G^-)$.  By \cite [Lemma 7.8(ii)]{J4} there exists $u' \in R\setminus \{u\}$ such that $f_v-f_{\varphi(v)}=(c_u-c_{u'})(r^{u+1}-r^u)$. Thus $v,\varphi(v)$ are neighbours in $\mathscr G^j(\textbf{c})$, through $u$.

  \

  Through the induction hypothesis $\mathscr G^{j_-}(\textbf{c}^-)$ is connected, for all $j_- \in N\setminus \{n-1\}$.  Given $k_-\in \{1,2,\ldots,n-2\}\setminus \{u\}$, set
   $$ k=\left\{
                       \begin{array}{ll}
                         k_-, & \hbox{if $k_-<u$} \\
                         k_-+1, & \hbox{if $k_->u$.}
                       \end{array}
                     \right. $$

   In this notation we can obtain $j$ to be any element of $N\setminus \{u,u+1\}$ by taking $j_- \in N\setminus \{u,n-1\}$.   Furthermore in the notation of \ref {7.2.1} one has
   $$\phi^+(k_-)=\phi^-(k_-)=k: \ \text{unless} \quad k_-=u.$$

   \

  \textbf{ Step 2.}

   Suppose that $v_1,v_2$  are neighbours in $\mathscr G^{j_-}(\textbf{c}^-)$, through $u'_-\in N \setminus \{n-1\}$.  By definition $j_- \notin \{u'_-,u'_-+1\}$.

  Suppose that $u \notin \{u_-',u'_-+1\}$.  On passing from $\mathscr G(\textbf{c}^-)$ to $\mathscr G^\pm$, the pair $(u'_-,u'_-+1)$ as labels on vertices become $(u',u'+1)$.  Observe further that $j \notin \{u',u'+1\}$.  We conclude that  $v_1,v_2$ (resp. $\varphi(v_1),\varphi(v_2)$) are neighbours in $(\mathscr G^+)^j$ (resp. $(\mathscr G^-)^j$) through $u'$.

  Suppose that $u'_-=u-1$.  Since the pair $(u-1,u)$ goes to $(u-1,u+1)$ (resp. $(u-1,u)$) on passing from $\mathscr G(\textbf{c}^-)$ to $\mathscr G^+$ (resp. $\mathscr G^-$) as labels on vertices, it follows that $\varphi(v_1),\varphi(v_2)$ are neighbours in $(\mathscr G^-)^j$, through $u'_-$; but $v_1,v_2$ are \textit{not} neighbours in $(\mathscr G^+)^j$.

   Suppose that $u'_-=u$.  Since the pair $(u,u+1)$ goes to $(u+1,u+2)$ (resp. $(u,u+2)$) on passing from $\mathscr G(\textbf{c}^-)$ to $\mathscr G^+$ (resp. $\mathscr G^-$) as labels on vertices, it follows that $v_1,v_2$ are neighbours in $(\mathscr G^+)^j$, through $u'_-+1$.  On other other hand $\varphi(v_1),\varphi(v_2)$ are \textit{not} neighbours in $(\mathscr G^-)^j$.

   \

   We conclude by Step 1, for any pair of neighbours $v_1,v_2 \in \mathscr G^{j_-}(\textbf{c}^-)$, that  $v_1,\varphi(v_1)$ and $v_2,\varphi(v_2)$ are neighbours in $\mathscr G^j(\textbf{c})$ and by Step 2, either $v_1,v_2$ or $\varphi(v_1),\varphi(v_2)$, or both, are neighbours in $\mathscr G^j(\textbf{c})$. Thus $\mathscr G^j(\textbf{c})$ is connected if $j \notin \{u,u+1\}$.

   \

 \textbf{ Step 3.}

  Assume $j=u$.  In this case $V^j(\mathscr G^+)$ is empty.  Thus it suffices to show that $(\mathscr G^-)^j$ is connected.   Suppose that $v_1,v_2$  are neighbours in $\mathscr G^u(\textbf{c}^-)$, through $u'_-$, which means in particular that $u\notin \{u'_-,u'_-+1\}$. Then as labels on vertices the pair $(u'_-,u'_-+1)$ goes to the pair $(u',u'+1)$ on passing from  $\mathscr G(\textbf{c}^-)$ to $\mathscr G^-$, whilst $u$ does not change.  Moreover $u \notin \{u',u'+1\}$, for otherwise $u= u_-'+1$.  Thus $\varphi(v_1),\varphi(v_2)$ are neighbours in $(\mathscr G^-)^u$, through $u'$.

   Assume $j=u+1$.  In this case $V^{j}(\mathscr G^-)$ is empty.  Thus it suffices to show that $(\mathscr G^+)^{j}$ is connected.   Suppose that $v_1,v_2$  are neighbours in $\mathscr G^u(\textbf{c}^-)$, through $u'_-$, which means in particular that $u\notin \{u'_-,u'_-+1\}$.  Then the pair $(u'_-,u'_-+1)$ goes to the pair $(u',u'+1)$  and $u$ goes to $u+1$ on passing from  $\mathscr G(\textbf{c}^-)$ to $\mathscr G^+$ as labels on vertices and $u+1 \notin \{u',u'+1\}$, for otherwise $u=u'_-$.  Thus $\varphi(v_1),\varphi(v_2)$ are neighbours in $(\mathscr G^+)^{u+1}$, through $u'$.

   \end {proof}

   \textbf{Remark 1}.  Retain the above notation.  We may conclude that as a graph $\mathscr G^j(\textbf{c})$ is a hypercube with some edges missing. Indeed if $j \in \{u,u+1\}$, then $\mathscr G^j(\textbf{c})$ is isomorphic to $\mathscr G^j(\textbf{c}^-)$.  Otherwise a line joining a pair of neighbours $v_1,v_2$ in $\mathscr G^{j_-}(\textbf{c}^-)$ gives rise to a square in $\mathscr G^{j}(\textbf{c})$ with lines joining the pairs $v_i,\varphi(v_i):i=1,2$ and at most one of the lines joining $v_1,v_2$ and $\varphi(v_1),\varphi(v_2)$ deleted.

  \

  \textbf{Remark 2}.  The choice of which line to delete in the square described in Remark 1, does \textit{not} depend on $j\in N$, but rather by $u' \in N\setminus \{u\}$. If $u'=u-1$ (resp. $u'=u+1$), then the line in $(\mathscr G^+)^j$  (resp. $(\mathscr G^-)^j$) is deleted.  Otherwise neither are deleted. In terms of the notation of \ref {7.1} deletions are determined by the \textit{natural} order on the integers $u_i$ satisfying $u_1\prec u_2 \prec \ldots \prec u_{r-1}$, and so by the function $\theta$ (which is a permutation of $N$).

\subsubsection{}\label{7.2.3}

In the notation of \ref {7.2.1} define a function in $Z(B_J) \mod \sum_{s' \in I\setminus \{s\}, k \in \mathbb N^+} \mathbb Zm_{s'}^k$ by
$$f=\sum_{j \in N}c'_j(r_s^j-r_s^{j+1})-\sum_{k\in \hat{N}} c_k m_s^k. \eqno {(41)}$$

  Here we recall (\ref {2.4}) that adding the difference $(r_s^{j}-r_s^{j+1})$ of successive Kashiwara functions of type $s$ is by $(3)$ implemented by adding the face $F_s^{j+1}$.
  %In terms of $\overline{v}_\textbf{k}$,
  %This means decreasing (resp. increasing) the exponent of $e_s$ in its $(j+1)^{th}$ (resp. $j^{th}$) place by $1$ in $\overline{v}_\textbf{k}$. The resulting element of $M_s(\textbf{a})$, if it is defined, will be proportional to $\overline{v}_\textbf{k}$, but not necessarily non-zero,

  In terms of the change in $z^K$, as defined in $(2)$, adding the above difference of successive Kashiwara functions becomes adding the function $z^{F_s^{j+1}}$.
  %Finally since $s \in I$ is fixed, it may be omitted.

  Thus $(41)$ may be rewritten as
$$f=\sum_{j \in N}c'_jz^{F_s^{j+1}}-\sum_{k\in \hat{N}} c_k m_s^k. \eqno {(42)}$$

Recalling that $c_n=0$, it follows from $(29)$ that the coefficient of $m_s^n$ in $f$ is zero if and only if $c'_{n-1}=0$.  Take $f=f_v$ for some $v \in V(\mathscr G(\textbf{c}))$.  Then by Property 1 of \ref {7.2.1} it follows that $c'_{n-1}=0$ if $i_v=n$ and only if given $c_{n-1}\neq 0$.

Recall the definition of the convex set $K(\textbf{c})$ given in \ref {7.1}.

Let $K^-(\textbf{c})$ (resp. $K^-_\mathbb Z(\textbf{c})$) denote the subset of $K(\textbf{c})$ (resp. $K_\mathbb Z(\textbf{c})$) defined by further imposing that $c'_{n-1}=0$. It is again a convex set.

 Suppose $c_{n-1}=0$, then $(37)$ forces $c'_{n-1}=0$, that is to say $K^-(\textbf{c})=K(\textbf{c})$.  Since we are trying to deduce properties of $K(\textbf{c})$ from those of $K^-(\textbf{c})$ we can assume that $c_{n-1}\neq 0$ without loss of generality. Thus from now on we adopt the

 \

 \textbf{Assumption.} $c_{n-1}\neq 0$.

  \

  It follows from this assumption and the above remarks that extremal elements of $K^-(\textbf{c})$ form the set $Z^-(\textbf{c}):=\{f_v\}_{v \in V^n(\mathscr G(\textbf{c}))}\subset K^-_\mathbb Z(\textbf{c})$.  One may remark that $K^-(\textbf{c})$ is reduced to the driving function (identified with $\{\textbf{0}\}$) if the $c_j:j \in N$ are increasing.

  Again by the definition of $T^-_s(\textbf{a})$ the coefficient in $z^K: K \in T^-_s(\textbf{a})$ of $m_s^n$ is zero.

%\textbf{Remark}.  In the generic case (that is when the coefficients $c_i$ are positive and pairwise distinct) the cardinality of $V(\mathscr G(\textbf{c}))$ is a power of $2$, namely $2^{n-1}$.  Again
%$V^j(\mathscr G(\textbf{c}))$ is a power of $2$.  This follows from the construction of binary fusion and the power can be calculated.   By contrast the cardinality of $Z^-(\textbf{c})$ is the Catalan number $C(n-1)$, the proof being more delicate \cite [6.7]{JL}.
\subsubsection{}\label{7.2.4}

Recall Lemma \ref {4.7.5} which was proved under the hypothesis that $K_{\ell \min}=K_{\min}$.  Here we prove a similar result but under a different hypothesis (stated in the first line of the lemma) and using the Chevalley-Serre relations via Lemma \ref {5.2}.

\begin {lemma}   Assume $\{z^K\}_{K\in T_s(\textbf{a})} \subset K_\mathbb Z(\textbf{c})$ and set $K^n=K_{\ell \min}$.
For all $j=n-1,n-2,\ldots,1$, adjoining $c_{j}$ copies of $F_s^{j+1}$ to $K^{j+1}$ gives a trail $K^{j}$.  Moreover $z^{K^j}=f_{v^j}$.
\end {lemma}

\begin {proof}   By definition $z^{K_n}$ is the driving function and by definition this is chosen to be $f_{v_n}$.  Again the vector in $M_s(\textbf{a})$ corresponding to $K_n$ is $v_\textbf{l}$.

Through the definition of the vertices $\{v_j\}_{j=1}^n$ of the pointed chain $(39)$ becomes using $(5)$ that
 $$f_{v_j}=f_{v_{j+1}}+c_jz^{F_s^{j+1}}, \forall j \in N. \eqno{(43)}$$

 By $(26)$ one has $c_j-c_{j-1}=a_j-\ell_j-\ell_{j-1}$.

 Now suppose that the trails $K^{j}:1<j\leq n$ have been constructed as described by the conclusion of the lemma and that $z^{K^{j}}=f_{v^{j}}$. Let $v_{\textbf{l}^{j}}\in M_s(\textbf{a})$ be the vector corresponding to $K^{j}$.  By definition, adjoining $c_j$ copies of $F_s^{j+1}$ to $K^{j+1}$ (to obtain $K^j$) means moving $c_j$ copies of $e$ to the right across $e_{a_{-(j+1)}}$. This gives in particular $\ell^j_j=\ell^{j+1}_j+c_j$ and $\ell^j_{j'}=\ell^{j+1}_{j'}$, for all $j'<j$.  As a consequence $\ell_j^j=\ell_j+c_j, \ell_{j-1}^j=\ell_{j-1}$.

 If $c_{j-1}=0$ we may set $K^{j-1}=K^{j}$.  Otherwise since $f_{v_{j}}$ is an extremal element of $K(\textbf{c})$, it follows from $(43)$ that $f_{v_{j}} - z^{F_s^{j}} \notin K(\textbf{c})$.  Then by the hypothesis of the lemma we cannot remove the face $F_s^{j}$ from $K^{j}$ to obtain a trail $K^{j-1}$.  This means that we cannot move a copy of $e$ to the left across $e_{-a_j}$ in $K^j$ to obtain a trail.

 Then applying Lemma \ref {5.2}, we can move $\ell_j^j+\ell_{j-1}^j-a_j=c_{j-1}$ copies of $e$ to the right across $e_{-a_j}$. This means that we can adjoin $c_{j-1}$ copies of $F_s^{j}$ to $K^j$ obtaining a new trail which we denote $K^{j-1}$.  Moreover in this $z^{K^{j-1}}=f_{v_j}+c_{j-1}z^{F_s^j}=f_{v_{j-1}}$, by $(43)$.

\end {proof}

 \subsubsection{}\label{7.2.5}  Both in \ref {7.2.4} and in Lemma \ref {4.7.5} the argument can be used to give more trails.  Indeed one may adjoin strictly less than $c_j$ faces $F^{j+1}$ at say the $j^{th}$  step.  Of course in this at the next step one is constrained to adjoin strictly less than $c_{j-1}$ faces $F^j$. In the rigid case one may easily check that the resulting trails exhaust $K_\mathbb Z(\textbf{c})$, that is to say the hypothesis $\{z^K\}_{K\in T_s(\textbf{a})} \subset K_\mathbb Z(\textbf{c})$ implies equality.  This argument fails to obtain equality in general.  The difficulty can be illustrated in the example $c_2> c_1 > c_3$.  In this case $V^4(\mathscr G(\textbf{c}))$ has four elements. They can be defined in terms of $\textbf{c}_j'=(c'_3,c'_2,c'_1)$, where $c_k'=c_k'^K(s)$ is given by $(30)$.  They are $(0,0,0),(0,c_2-c_3,0),(0,c_2-c_3,c_1-c_3),(0,c_1-c_3,c_1-c_3)$.  The first corresponds to the $\ell$-minimal trail.  Appealing to Lemma \ref {5.2} (and checking it applies!), the second is obtained by adjoining $(c_2-c_3)$ copies of $F^3$ to the first and then the third is obtained by adjoining $c_1-c_3$ copies of the face $F^2$ to the second. Finally for the fourth trail one ``doubles-back'' and \textit{removes} $c_2-c_1$ copies of $F^3$ from the third trail. Though one can first adjoin just $c_1-c_3$ copies of $F^3$ to the minimal trail, one \textit{cannot} then further adjoin $c_1-c_3$ copies of $F^2$ to obtain the fourth trail by appealing to Lemma \ref {5.2}.

 %This kind of acrobatics should allow one to use Lemma \ref {7.2.2} to convert the hypothesis of Lemma \ref {7.2.4} into an equality (as we already saw holds for the rigid case); but it because increasingly precarious as $r$ becomes large.  In the section we shall modify this approach.

 This example also illustrates Remark \ref {7.2.4}.  Here $u=2$ and $j=4$, whilst $\mathscr G^{j_-}(\textbf{c}^-)$ is a graph with two vertices $v_1,v_2$ joined by an edge with $u_-=1$. The corresponding functions are $f_{v_1}=(0,0)$ and $f_{v_2}=(0,c_1-c_2)=f_{v_1}+(c_1-c_2)(r^1-r^2)$.  In $\mathscr G^+$ (resp. $\mathscr G^-$), these functions become $f_{v_1}=(0,0,0),f_{v_2}=(0,c_1-c_3,c_1-c_3)=(0,0,0)+(c_1-c_3)(r^1-r^3)$ (resp. $f_{\varphi(v_1)}=(0,c_2-c_3,0), f_{\varphi(v_2)}(0,c_2-c_3,c_1-c_3)=(0,c_2-c_3,0)+(c_1-c_3)(r^1-r^2)$).  Thus in $(\mathscr G^+)^j \hookrightarrow\mathscr G^j(\textbf{c})$ the edge joining $v_1$ and $v_2$ is deleted.

 This kind of acrobatics should allow one to use Lemma \ref {7.2.2} to convert the hypothesis of Lemma \ref {7.2.4} into an equality (as we already saw holds for the rigid case); but it because increasingly precarious as $n$ becomes large.  The next subsections constitute a modification of modify this approach.  In some sense what we have done so far is a warm up exercise to the real McCoy.

 %This example illustrates another phenomenon which we would like to extend.  In this $c_2$ is the largest element of $\textbf{c}$ and We would like to suppress the largest element $c_u$ of $\textbf{c}$ as part of an induction procedure. This is what we did in \ref {7.2.2} for binary fusion. Yet concerning trails there is an additional obstruction. Indeed to do this for trails we need to replace the factor.

\subsubsection{}\label{7.2.6}

%We now prove a quite remarkable and unexpected property of $K_\mathbb Z(\textbf{c})$.

Recall the notation of \ref {3.2},\ref {4.2} using the convention of \ref {4.3}.  Consider an element $\overline{v}_\textbf{k}$ of the form given in $(16)$\textit{ not} for the moment assumed to be non-zero.  We identify $\textbf{k}$ with a trail $K \in K_\mathbb Z(\textbf{c})\subset \mathbb Z^n$.  In this $\textbf{a}$ is assumed fixed; but the components of $\textbf{k}$ may vary. Using this identification $\textbf{l}$ identifies with $K_{\ell \min}(s)$ and we take $\textbf{c}$ to be given by $(26)$ with $c_j=c_j^{K_{\ell \min}}(s):j \in N$, which we recall are non-negative integers.  Here $\overline{v}_\textbf{l}$ \textit{is} assumed to be non-zero.

Following $(30)$ and the above identification, we set $c_j'^{\textbf{k}}=k^j-\ell^j:j\in \hat{N}$, or simply $c_j'$, if $\textbf{k}$ is fixed. One has $c_n'=0$. The possible values of $\{c_j'^{\textbf{k}}\}_{j\in N}$, for $\textbf{k} \in K_\mathbb Z(\textbf{c})$  are determined by $(37),(38)$ and of course the requirement that they be integer.

Take $u \in \hat{N}$.
% and $\textbf{k}\in K_\mathbb Z(\textbf{c})$.

A line of type $u$ in $K(\textbf{c})$ through $\textbf{k}\in K_\mathbb Z(\textbf{c})$ is defined to be the set $\{\textbf{k}(v): v \in \mathbb Q$, where
$$ k_j(v)=\left\{\begin{array}{ll}
                        k_j+v, & \hbox{if $j=u$,} \\
                        k_j-v, & \hbox{if $j=u+1$,}\\
                        k_j,& \hbox{otherwise}.
                       \end{array} \right. \eqno {(44)}$$

Recall that $k^n=\ell^n$ and view $\textbf{k}$ as being given by the $n-1$-tuple whose $j^{th}$ entry is $c'_j=k^j-\ell^j$.  In these new co-ordinates
 a line of type $j$ passing through $\textbf{k}\in K_\mathbb Z(\textbf{c})$ becomes an affine translate by multiples of the $j^{th}$ co-ordinate function. Indeed one checks from $(44)$ that $k^{i}(v)=k^i:i \neq j$, for all $v \in \mathbb Q$.

 Notice by the discussion following $(41)$ that, in terms of this $(n-1)$-tuple, the $j^{th}$ co-ordinate function is also represented as the successive difference $r^j-r^{j+1}$ of Kashiwara functions.

Let $\prec$ denote the order relation on $N$ introduced in \ref {7.1}.

\begin {lemma} Let $u$ be the unique maximal element of $N$ for $\prec$.  A line of type $u$ in $\mathbb Z^r$ through $\textbf{k}\in K_\mathbb Z(\textbf{c})$ has length $\max (k_{u+1}+k_u-a_{u+1},0)$ in $K_\mathbb Z(\textbf{c})$ and consists of  $\max (1+k_{u+1}+k_u-a_{u+1},1)$ elements of $K_\mathbb Z(\textbf{c})$.
\end {lemma}

\begin {proof} Combining $(26)$ and $(30)$ we obtain
$$k_{u+1}+k_u-a_{u+1}=(c'_{u+1}+c_u-c_{u+1}) -c'_{u-1}. \eqno {(45)}$$

%Recall that $k^r=\ell^r$ and view $\textbf{k}$ as being given by the $r-1$-tuple whose $j^{th}$ entry %is the successive difference $r_s^j-r_s^{j+1}$ of Kashiwara functions (or of the $z^{F^{j+1}_s}$).
%$k^j-\ell^j$.
%
%\
%
% \textbf{N.B.} A line of type $j$ passing through $\textbf{k}\in K_\mathbb Z(\textbf{c})$ becomes in these new co-ordinates, an affine translate by multiples of the $j^{th}$ co-ordinate function (because in particular $k^{j+1}$ does not change).
%
% \

Recall binary fusion and the notation of \ref {7.2.2}.

 A typical element of $Z(\textbf{c}^-)$, that is to say a function assigned to a vertex of $\mathscr G(\textbf{c}^-)$ takes the form $\textbf{c}'=(c'_{n-2},c'_{n-3},\ldots,c'_1)$.  Then (\cite [Eq. $(12)$]{J5}) viewed as the function assigned to the corresponding vertex of $V(\mathscr G^+)$, when $\mathscr G^+$ is viewed as a subgraph of $\mathscr G(\textbf{c})$ it becomes
 $$\textbf{c}'_+:=(c_{n-1}',c_{n-2}', \ldots, c_{u+1}', c_{u-1}',c_{u-1}',c'_{u-2},c_{u-3}'\ldots,c_1'). \eqno {(46)}$$

  On the other hand (\cite [Eq. $(17)$]{J5}) as the function assigned to the corresponding vertex of $V(\mathscr G^-)$, when $\mathscr G^-$ is viewed as a subgraph of $\mathscr G(\textbf{c})$ it takes the form
 $$\textbf{c}'_-:=(c_{n-1}',c_{n-2}', \ldots, c_{u+1}', c_{u+1}'+(c_u-c_{u+1}),c_{u-1}',c'_{u-2},c_{u-3}'\ldots,c_1').\eqno{(47)}$$
 % as the function assigned to the corresponding vertex of $V(\mathscr G^-)$, when $\mathscr G^-$ is viewed as a subgraph of $\mathscr G(\textbf{c})$.

%As noted in the proof of \cite [Prop. 3.5]{J5} a typical element $\textbf{c}'$ of $Z(\textbf{c}^-)$, that is to say a function assigned to a vertex of $\mathscr G(\textbf{c}^-)$ takes the form (\cite [Eq. $(12)$]{J5})
% $$\textbf{c}'':=(c_n',c_{n-1}', \ldots, c_{u+1}', c_{u-1}',c_{u-1}',c'_{u-2},c_{u-3}'\ldots,c_1'), \eqno {(45)}$$
% as the function assigned to the corresponding vertex of $V(\mathscr G^+)$, when $\mathscr G^+$ is viewed as a subgraph of $\mathscr G(\textbf{c})$. On the other hand it takes the form (\cite [Eq. $(17)$]{J5})
% $$\textbf{c}''':=(c_n',c_{n-1}', \ldots, c_{u+1}', c_{u+1}'+(c_u-c_{u+1}),c_{u-1}',c'_{u-2},c_{u-3}'\ldots,c_1'),\eqno{(46)}$$
%  as the function assigned to the corresponding vertex of $V(\mathscr G^-)$, when $\mathscr G^-$ is viewed as a subgraph of $\mathscr G(\textbf{c})$.

  It was also noted (\cite [line above Eq. (18)] {J5}) that the expression in $(45)$ is non-negative. By linearity the expressions in $(46),(47)$ and this last assertion remain valid on taking \textit{convex} linear combinations of elements of $Z(\textbf{c}^-)$, so in particular we can assume that $\textbf{c}' \in K^-_\mathbb Z(\textbf{c})$.

  By taking convex linear combinations of the elements defined by $(46),(47)$ we may obtain any element $\tilde{\textbf{c}}$ differing from $\textbf{c}'_+$ and $\textbf{c}'_-$ just in its $u^{th}$ co-ordinate $\tilde{c}_u$ and which lies between (the non-negative integers) $c'_{u-1},c'_{u+1}+c_u-c_{u+1}$.  Conversely if $\tilde{c}_u \notin c'_{u-1},c'_{u+1}+c_u-c_{u+1}$, then $\tilde{c} \notin K(\textbf{c})$ by \cite [Lemma 2.3]{J5}.  (These facts led to $K(\textbf{c})$ being the convex hull of $Z(\textbf{c})$, established in \cite [Thm. 1.4]{J5}.)

  Combined with $(45)$ and the remarks following it, the conclusion of the lemma obtains.

  \end {proof}

  \textbf{Remark}.  One has $k_{u+1}+k_u-a_{u+1}=(c'_{u+1}-c'_u)-(c_{u+1}-c_u) \geq 0$, by $(38)$ and the definition of $u$. Thus the statement of the lemma may be simplified.  However as stated it prompts the question as to whether the statement is true for all $j \in N$.   This is false, as can seen in the example $n=3,c_1=3,c_2=2$. In this case $u=1$.  The lines of type $1$ whose second co-ordinate is $0,1,2$ meets $K(\textbf{c})$ in $2,3,4$ points respectively, as required by the conclusion of the lemma.  However the lines of type $2$ whose first co-ordinate is $0,1,2,3$ meets $K(\textbf{c})$ in $3,3,2,1$ points respectively, so the lemma cannot be extended to this case. However it extends (in this case) when both end points lie in $Z(\textbf{c})$.

   \subsubsection{}\label{7.2.7}

    Following the above remark, we extend the previous lemma for lines of type $u_1 \in N$ joining elements of $Z(\textbf{c})$.  Recall Definition \ref {7.2.2}.

   \begin {lemma}  Take $j \in N$ and let $v_1,v_2$ be neighbours in $\mathscr G^j(\textbf{c})$ through some $u_1 \in N$.  Then on the line $L$ of type $u_1$ joining $v_1,v_2$ there are exactly $1+k_{u_1+1}+k_{u_1}-a_{u_1+1}$ elements of $K_\mathbb Z(\textbf{c})$, where $\textbf{k}=\{k_i\}_{i \in \hat{N}}$ is defined by either $v_1$ or $v_2$.
   \end {lemma}

   \begin {proof}  By definition one has $v_1,v_2 \in Z(\textbf{c})$.  By \cite [Thm. 1.4]{J5}, the elements of $Z(\textbf{c})$ are the extremal elements of $K(\textbf{c})$, which is the convex set they define.  Thus $L \cap K(\textbf{c})$ exactly consists of convex linear combinations of $v_1,v_2$.

   In the notation of Definition \ref {7.2.2}, the functions $f_{v_1},f_{v_2}$ only differ in the $u_1^{th}$ co-ordinate and by the integer $c_{u_1}-c_{u_2}$, which as we shall see must be non-negative. Thus the cardinality of $L\cap K_\mathbb Z(\textbf{c})$ is just $1+(c_{u_1}-c_{u_2})$.

   The proof is completed by induction on $n$.  If $n=1$, there is nothing to prove.

   If $u_1$ is the maximal element $u$ of $N$ for the order relation $\prec$, then the assertion holds by Lemma \ref {7.2.6}.  In this we note that by $(44)$ the value of $k_{u+1}-k_u-a_{u+1}$ is independent of the choice of $\textbf{k} \in L$, it is non-negative by Remark \ref {7.2.6} and equals $c_{u_1}-c_{u_2}$ by the previous paragraph.  (The latter can directly seen to be non-negative.  Indeed $u_2 \prec u_1=u$ means that $c_{u_1}\geq c_{u_2}$.)

   If $u_1$ is not the maximal element $u$ of $N$, then both $v_1,v_2$ can be assumed to lie in $V(\mathscr G^+)$ or in $V(\mathscr G^-)$. Yet by Remark \ref {7.2.1}, the functions defined by $V(\mathscr G^+)$ (resp. $V(\mathscr G^-)$), identify with those given by the functions defined by $V(\mathscr G(\textbf{c}^-))$  (resp. and up to an overall translation).  Thus we can take $v_1,v_2$ to lie in $V(\mathscr G(\textbf{c}^-))$ and apply the induction hypothesis.

     %(Of course the convex sets they define, differ from $K(\textbf{c})$, and this is why we need the first paragraph of the proof.)

   \end {proof}

 \subsubsection{}\label{7.2.8}

 \begin {cor}

\

(i) $\{z^K\}_{K\in T_s^-(\textbf{a})} \subset K^-_\mathbb Z(\textbf{c})$ implies
that $Z^-(\textbf{c}) \subset \{z^K\}_{K\in T_s^-(\textbf{a})}$.

\

(ii) $\{z^K\}_{K\in T_s(\textbf{a})} \subset K_\mathbb Z(\textbf{c})$ implies that $Z(\textbf{c}) \subset \{z^K\}_{K\in T_s(\textbf{a})}$.

\end {cor}

\begin {proof}  Recall (\ref {7.2.1}) the definition of $v_r \in V^r(\mathscr G(\textbf{c}))$. It is assigned to the unique $\ell$-minimal element of $T_s(\textbf{a})$ which in particular belongs to $T^-_s(\textbf{a})$. Then taking $j=n$ in Lemmas \ref {7.2.2}, \ref {7.2.7} we obtain (i) through Lemma \ref {5.2}.  Again by Lemma \ref {7.2.4} to each element $v_j\in V^j(\mathscr G(\textbf{c})):j=n,n-1,\ldots,1$ of the pointed chain, there is a trail in $T_s(\textbf{a})$ and so as above (ii) results through applying Lemma \ref {5.2} to Lemmas \ref {7.2.2}, \ref {7.2.7}.
\end {proof}

 \subsubsection{}\label{7.2.9}

   Assume that

   \

   $(*)$ \quad  $\{z^K\}_{K\in T_s(\textbf{a})} \subset K_\mathbb Z(\textbf{c})$.

   \

    Take $u$ to be the unique maximal element of $N$ for $\prec$.

   Suppose that $\textbf{k}$ is a trail, equivalently that $\overline{v}_\textbf{k}\neq 0$. Then by Lemma \ref {5.2} applied to the conclusion of Lemma \ref {7.2.6} we conclude that the points in $K_\mathbb Z(\textbf{c})$ lying on the line of type $u$ through $\textbf{k}$ are also trails.

   Recall that $Z(\textbf{c})$ is defined to be the set of functions obtained from the vertices of $\mathscr G(\textbf{c})$ and $K_\mathbb Z(\textbf{c})$ to be their $\mathbb Z$ convex hull.  Recall further that $\mathscr G(\textbf{c})$ is obtained by an appropriate joining of the graphs $\mathscr G^\pm$.
    %\begin{figure}
%  \centering
%  % Requires \usepackage{graphicx}
%  \includegraphics[width=]{}\\
%  \caption{}\label{}
%\end{figure}
%   $\mathscr G^\pm$.

   Now let $Z(\textbf{c})^\pm$ be the subset of functions obtained from the vertices of $\mathscr G^\pm$ and $K_\mathbb Z(\textbf{c})^\pm$ their $\mathbb Z$ convex hull\footnote{$Z(\textbf{c})^-$ and $K_\mathbb Z(\textbf{c})^- $
   are not to be confused with $Z^-(\textbf{c})$ and $K^-_\mathbb Z(\textbf{c})$ defined in \ref {7.2.3}.}.  If we could say that any element of $K_\mathbb Z(\textbf{c})^+$ or of $K_\mathbb Z(\textbf{c})^-$ is a trail, then we may conclude by the previous argument that any element of $K_\mathbb Z(\textbf{c})$ is a trail.
   At first sight this might seem easy as by induction any $\overline{v}_\textbf{k} \in K_\mathbb Z(\textbf{c}^-)$ is a trail.  However here there is a difficulty which one might not spot at first sight.  Indeed the change in the numbering on going from $\mathscr G(\textbf{c}^-)$ to $\mathscr G^\pm$ upsets the notion of a trail.  This is illustrated below.

   Notice that in the expression for $\textbf{c}'_+$ given in $(45)$ the doubling of $c'_{u-1}$ means that in terms of $\textbf{k}$ one has $k_u=\ell_u$.  Thus in $\overline{v}_\textbf{k}$ we have the sub-product $e_a:=e_{-a_{u+1}}e_s^{\ell_u}e_{-a_u}$.  In order to use the fact that the corresponding element coming from $K_\mathbb Z(\textbf{c}^-)$ is a trail, we would have to be able treat this subproduct as if it were a product of negative root vectors distinct from $e_s$ of $h_s$ eigenvector of eigenvalue $a:=2\ell_u-a_{u+1}-a_u$, so that in particular Lemma \ref {5.2} can be applied.  This is not obviously true since for example it is not annihilated by$(\ad e_s)^{a+1}$.

   %Now in applying Lemma \ref {5.2} to assert that $v_\textbf{k} \in K_\mathbb Z(\textbf{c}^-)$ is a trail. we would like to treat $e_a$ not only as having $h_s$ eigenvalue $a:=2\ell_u-a_{u+1}-a_u$, but also with respect to the action of $\ad e_s$ to behave as a product of simple positive root vectors distinct from $e_s$ and in particular to be annihilated by $(\ad e_s)^{a+1}$, which of course is false.

   %We show nevertheless that this sub-product does behave in the required fashion to recover Lemma \ref {5.2}.

   We shall overcome this difficulty by the same process we used to obtain Corollary \ref {7.2.8}.  The key point is that every square defined by Remark $1$ of \ref {7.2.2} at most one edge is missing (which had the consequence that $\mathscr G^j(\textbf{c})$ is connected (Lemma \ref {7.2.2})).
%
%  The key point is that the above difficulty is the same as that overcome in the proof of Lemma \ref {7.2.2}.  In particular combining this result with Lemmas \ref {7.2.4}, \ref {7.2.7} and Lemma \ref {5.2} we already obtain from Lemma \ref {5.2} that $(*)$ implies that $Z(\textbf{c}) \subset \{z^K\}_{K\in T_s(\textbf{a})}$.

%This example illustrates another phenomenon which we would like to extend.  In this $c_2$ is the largest element of $\textbf{c}$ and We would like to suppress the largest element $c_u$ of $\textbf{c}$ as part of an induction procedure. This is what we did in \ref {7.2.2} for binary fusion. Yet concerning trails there is an additional obstruction. Indeed to do this for trails we need to replace the factor.
%
%\end {proof}

\subsection{}\label{7.3}

 The key point in describing a set of dual Kashiwara functions of type $t\in I$, is to show that $\{z^K\}_{K \in T^+_s(\textbf{a})}=K_\mathbb Z(\textbf{c})$, where $\textbf{c}$ is defined by $T^+_s(\textbf{a})$ as in \ref {5.3.3}, and this for all choices of $\textbf{a}$ and all $s \in I$.

Notice that an element of $K_\mathbb Z(\textbf{c})$ defines a ``potential trail'' trivializing at the fixed element extremal vector of weight $-w_j\varpi_t$, that is to say a product of the form $\overline{v}_\textbf{k}$ (as defined in $(16)$).   Recall that knowing when it does not vanish is the hard part of the present story.  In this case we are calling it a trail.

We start with the hypothesis that $\{z^K\}_{K\in T_s^-(\textbf{a})} \subset K^-_\mathbb Z(\textbf{c})$.  Notice that $K_{\ell \min} \in T^-_s(\textbf{a})$.  If the $\{c_k\}_{k \in N}$ are increasing this just means that $T^-_s(\textbf{a})$ is reduced to $\{K_{\ell \min}\}$.  This was the hypothesis of Theorem \ref {6.4} and its conclusion is that {\newline} $\{z^K\}_{K \in T^+_s(\textbf{a})}=K_\mathbb Z(\textbf{c})$. The general case is surprisingly much more difficult.

A first step for the general case is provided by the

\begin {prop}

\

(i) $\{z^K\}_{K\in T_s^-(\textbf{a})} \subset K^-_\mathbb Z(\textbf{c})$ implies equality.
%that $Z^-(\textbf{c}) \subset \{z^K\}_{K\in T_s^-(\textbf{a})}$.

\

(ii) $\{z^K\}_{K\in T_s(\textbf{a})} \subset K_\mathbb Z(\textbf{c})$ implies equality.
%that $Z(\textbf{c}) \subset  \{z^K\}_{K\in T_s(\textbf{a})}$.

\end {prop}

\begin {proof}  Notice that equality in $(i)$ (resp. in $(ii)$) means in the above language that a potential trail defined by an element of $K^-_\mathbb Z(\textbf{c})$ (resp. $K_\mathbb Z(\textbf{c})$), is in fact a trail.

It is clear that the proof of (i) follows exactly that of (ii) since it simply means imposing $c'_{n-1}=0$ in $(37),(38)$, that is in terms of co-ordinates the first entry is taken to be zero.

The proof of (ii) obtains by induction on the cardinality of the index set $\hat{N}$, namely $n$.  The assertion is empty for $n=1$ and an easy consequence of Lemma \ref {5.2} for $n=2$, as obtained in the more general result noted in \ref {7.2.5}.

Now choose $u \in N$ as in \ref {7.2.2} and retain the notation used there.
%An element of $K_\mathbb Z(\textbf{c}^-)$ is a trail by the induction hypothesis.

Let $\phi^\pm$ be the labelled graph embedding of $\mathscr G^\pm$ into $\mathscr G(\textbf{c})$ defined in \ref {7.2.1}.  It implements a relabelling of vertices (given in \ref {7.2.1}) and a relabelling of coefficients, given by $\phi^\pm(\textbf{c}'):=\textbf{c}'_\pm$, the latter being given by $(46),(47)$.

 %Given $\textbf{c}' \in Z(\textbf{c}^-)$, one has  $\textbf{c}'_\pm=\theta_\pm(\textbf{c}')$, where the latter are given by $(46),(47)$. %In terms of the co-ordinates for which the $j^{th}$ entry is $k^j-\ell^j$, this new element is given by $(44)$, roughly speaking by doubling the $(u-1)^{th}$ entry.  This just means setting $k_u=\ell_u$.

 As in \ref {7.2.9}, let $K_\mathbb Z(\textbf{c})^\pm$ denote the subset of $K_\mathbb Z(\textbf{c})$ which lies in the convex hull of $\im \phi^\pm$.

 Recall (\ref {7.1}) that by definition $u_{i+1}$ is the successor of $u_i$ for $\prec$ and that $u_{n-1}=u$.

 %As noted in $(44)$, the elements in $\im (\mathscr G^+ \rightarrow\mathscr G(\textbf{c})$, an element of $\textbf{c}' \in  Z(\textbf{c}^-)$ becomes the element of $\textbf{c}'' \in Z(\textbf{c})$ given by doubling the $(u-1)^{th}$ entry (in the co-ordinates for which $j^{th}$ entry is $k^j-\ell^j$)  and this means putting $k_u=\ell_u$.  Let $\overline{K}_\mathbb Z(\textbf{c})$ denote the subset of $K_\mathbb Z(\textbf{c})$

 %Suppose that $u-1 \prec u+1 \prec u$ (resp. $u+1 \prec u-1 \prec u$).   We claim that every element of $K_\mathbb Z(\textbf{c})^+$ (resp. $K_\mathbb Z(\textbf{c})^-$) is a trail.

 Suppose that $u_{n-2}>u_{n-1}$ (resp. $u_{n-2}< u_{n-1}$).

 \

\textbf{ Claim.}  Every element of $K_\mathbb Z(\textbf{c})^+$ (resp. $K_\mathbb Z(\textbf{c})^-$) is a trail.

\

  %An element of $K_\mathbb Z(\textbf{c}^-)$ is a trail by the induction hypothesis.  It follows by Remark 2 of \ref {7.2.2} that the elements of $(K_+)_\mathbb Z(\textbf{c})$ (resp. $(K_-)_\mathbb Z(\textbf{c})$) are trails if $u_{r-2} > u_{r-1}$ (resp $u_{r-2} < u_{r-1}$).

  It is clear (from the last part of the proof of Lemma \ref {7.2.6}) that an element $K$ of $K_\mathbb Z(\textbf{c})$ lies on a line of type $u$ joining an element of $K_\mathbb Z(\textbf{c})^+$ to an element of $K_\mathbb Z(\textbf{c})^-$ in $K_\mathbb Z(\textbf{c})$.  Through the claim and combining Lemmas \ref {5.2}, \ref {7.2.6} with the hypothesis of the proposition we deduce that $K$ is a trail.  This gives (i) up to the claim.

   To prove the claim suppose first that $u_{n-2}>u_{n-1}$. % and set $u':=u_{n-2}-1 \geq u$.
    Recall that $\mathscr G^+$ (resp. $\mathscr G^-$) decomposes as the union of graphs $\mathscr G^{+,+}$ and $\mathscr G^{+,-}$ (resp. $\mathscr G^{-,+}$ and $\mathscr G^{-,-}$) obtained by suppressing the edges with label $u_{n-2}$.

    Consider the subgraph $\mathscr G(\tilde{\textbf{c}})$ of $\mathscr G(\textbf{c})$ obtained as the binary fusion of $\mathscr G^{+,-}$ and $\mathscr G^{-,-})$) that is to say by joining the vertex $v \in V^{u+1}(\mathscr G^{+,-})$ to the corresponding vertex, namely $\varphi(v)$ of $V^{u}(\mathscr G^{-,-})$.  It is just the $S$-graph with coefficient set $\tilde{\textbf{c}}:=\textbf{c}\setminus \{c_{u_{r-2}}\}$.  Since  $\tilde{\textbf{c}}$ has cardinality $r-2$, it follows by the induction hypothesis that the $\mathbb Z$ convex set of functions defined by $\mathscr G(\tilde{\textbf{c}})$ are all trails.
     %(and the resulting set of trails identifies with $K_\mathbb Z(\textbf{c}_{u'})$).
     In particular the $\mathbb Z$ convex hull of the functions defined by $\mathscr G^{+,-}$ are all trails.

  Now in $\mathscr G^+$ the label $u_{n-1}$ does not appear on edges and so for $\prec$ the largest label on edges is $u_{n-2}$. Moreover by our assumption that $u_{n-2}>u_{n-1}$, it follows that that a vertex $v \in V^{u_{n-2}+1}(\mathscr G^{+,+})$ is joined to the corresponding vertex,
  namely $\varphi(v)$ in $V^{u_{n-2}}(\mathscr G^{+,-})$ in $\mathscr G^+$ by an edge with label $u_{n-2}$.  Thus $\mathscr G^+$ is the binary fusion of union of these two graphs. (This can fail for $\mathscr G^-$ if $u_{n-2}=u_{n-1}+1$, since then a vertex  $v \in V^{u_{n-2}+1}(\mathscr G^{-,+})$ is joined to the corresponding vertex,
  namely $\varphi(v)$ in $V^{u_{r-2}-1}(\mathscr G^{+,-})$ in $\mathscr G^+$ by an edge with label $u_{n-2}$.)  Then by Lemmas \ref {5.2}, \ref {7.2.6}, the hypothesis of the proposition and the result of the previous paragraph, the $\mathbb Z$ convex hull of the functions defined by $\mathscr G^+$ are all trails.

  The case $u_{n-2}<u_{n-1}$ is similar, except that we take  $\mathscr G(\tilde{\textbf{c}})$ to be the subgraph of $\mathscr G(\textbf{c})$ obtained by joining the vertex $v \in V^{u+1}(\mathscr G^{+,+})$ to the corresponding vertex, namely $\varphi(v)$ of $V^{u}(\mathscr G^{-,+})$. As in the first case, it follows that the $\mathbb Z$ convex hull of the functions defined by $\mathscr G^{-,+}$ are all trails. By our assumption that $u_{n-2}<u_{n-1}$, it follows that that a vertex $v \in V^{u_{n-2}+1}(\mathscr G^{-,+})$ is joined to the corresponding vertex,
  namely $\varphi(v)$ in $V^{u_{n-2}}(\mathscr G^{-,-})$ in $\mathscr G^-$, which is hence the binary fusion of these two graphs. (Again this can fail for $\mathscr G^+$ if $u_{n-2}=u_{n-1}-1$.) Then by Lemmas \ref {5.2}, \ref {7.2.6}, the hypothesis of the proposition and the result of the previous paragraph, the $\mathbb Z$ convex hull of the functions defined by $\mathscr G^-$ are all trails.

  This proves the claim and hence the proposition.

  %The labels in $\mathscr G_u^-$ are exactly the same as the labels in the subgraph $\mathscr G_u(\textbf{c}^-)$ of $\mathscr G(\textbf{c})^-)$ obtained by restricting all labels to be $\leq u$. Moreover the latter is just the $S$-graph $\mathscr G(\textbf{c}_{u-1})$ in which the coefficient set is reduced to $\textbf{c}_{u-1}:= \{c_i\}_{i=1}^{u-1}$.  Thus by the induction hypothesis the $\mathbb Z$ convex set of functions it defines are all trails (and the resulting set of trails identifies with $K_\mathbb Z(\textbf{c}_{u-1})$).  Then combining Lemmas \ref {5.2}, \ref {7.2.6} with the hypothesis of the proposition we deduce that the $\mathbb Z$ convex set of functions defined by $\mathscr G_u(\textbf{c})$ are all trails (and the resulting set of trails identifies with $K_\mathbb Z(\textbf{c}_{u})$).  In particular the $\mathbb Z$ convex hull of the functions defined by $\mathscr G_u^+$ are all trails.
\end {proof}

\textbf{Remark}.  Notice that Corollary \ref {7.2.8} has been made redundant. Nevertheless elements of the latter arise in the present proof.  Indeed the square with edges being the subgraphs $\mathscr G^{+,+},\mathscr G^{+,-},\mathscr G^{-,+},\mathscr G^{-,-}$, joined by an edge whenever the pair leads to a binary fusion (so that at most one edge is missing, the latter when $u_{r-2}=u_{r-1}\pm 1$) is the exact analogue of the square described in Remark $2$ of \ref {7.2.2}.  Also the proof of \ref {7.2.8} is somewhat easier, so makes a good preamble.

\subsection{The Lower Bounds}\label{7.4}  We would like to prove that the hypothesis of the first part of Proposition \ref {7.3} implies the hypothesis of the second part, as in the rigid case, using $\mathfrak {sl}(2)$ theory. For this we must show that $(37),(38)$ hold for any trail $K\in T_s(\textbf{a})$ given that they hold for any trail $K\in T^-_s(\textbf{a})$.  Here if a trail $K$ is given by a vector $\overline{v}_\textbf{k}$ of the form described in $(16)$, then the corresponding elements $c'_j:j \in N$ to be substituted into $(37),(38)$ are the $c'^K_j(s)$ (or simply $c'^K_j$) given by $(35)$, where in the latter the $\ell_j:j \in \hat{N}$ are given in terms of the unique $\ell$-minimal trail.

The key point which allows the use of $\mathfrak {sl}(2)$ theory is that the $\overline{v}_\textbf{k}$ belong to a simple module $M_s(\textbf{a})$ (Lemma \ref {4.3}).  In this  the elements of $T^-_s(\textbf{a})$ are proportional to its lowest weight vector.

We shall refer to the first part of $(37)$, namely $0\leq c_j':j \in N$ together with $(38)$ as the \textbf{lower bounds in $(37),(38)$}.  Here we may note that the right hand sides in $(38)$ are determined by the unique $\ell$-minimal trail and in particular are fixed.

In the rigid case the lower bounds are equivalent to the inequalities $k_j\geq \ell_j$, for all $j \in N$.   These are an obvious consequence of the hypothesis that $T^-_s(\textbf{a})$ is reduced to a single element being the unique minimal trail (and the simplicity of $M_s(\textbf{a})$). The general case is barely more difficult.

\begin {lemma}  Suppose the lower bounds in $(37),(38)$ hold for all $K\in T^-_s(\textbf{a})$, then they hold for all $K \in T_s(\textbf{a})$.
\end {lemma}

\begin {proof}  As in \ref {3.2} one may observe that the action of powers of the lowering operator $f_s$ on a trail $v_\textbf{k}$ gives a sum of expressions $v_{\textbf{k}'}$ in which one exponent of $e_s$ is reduced by one.  Moreover any such non-vanishing term corresponds to a trail in $T_s(\textbf{a})$.  Now since the only trails in $T_s(\textbf{a})$ which are non-zero multiples of the \textit{unique up to scalars} lowest weight vector of $M_s(\textbf{a})$ lie in  $T^-_s(\textbf{a})$, it follows that $k_j\geq k'_j:j \in \hat{N}$, for any $v_{\textbf{k}'}$ appearing in such a sum and so in particular for some trail $v_{\textbf{k}'}$ lying in $T^-_s(\textbf{a})$.  Now since $c'^K_j=k^{(j)}-\ell^{(j)}$ for any trail $K$, the lower bounds for $\textbf{k}'$ together with the inequalities $k_j\geq k'_j:j \in \hat{N}$ imply the lower bounds for $\textbf{k}$.   Hence the assertion of the lemma.
\end {proof}

\subsection{The Upper Bounds}\label{7.5}

We shall refer to the second part of $(37)$, namely $c_j'\leq c_j$ as the upper bounds (in $(37),(38)$).  Substitution in $(34),(35)$ for a given trail $K \in T_s(\textbf{a})$ shows that these upper bounds are equivalent to
$$a^{(j)}-k^{(j)}-\ell^{(j-1)}\geq 0, \forall j \in \hat{N}.\eqno {(47)}$$

The set of all upper and lower bounds for $T^-_s(\textbf{a})$ (resp. $T_s(\textbf{a})$)  is equivalent to the hypothesis of $(i)$ (resp. $(ii)$) of Lemma \ref {7.3}.

The upper bounds are much more difficult to establish and their proof will need most of the results of Section \ref {7.2}.
\subsubsection{}\label{7.5.1}

  Recall that $\textbf{l}$ defines the unique $\ell$-minimal trail and note that $\ell_r$ does not appear in $(47)$, so we do not have to worry about whether we are setting $\ell_r=0$ which makes the unique minimal trail an element of $T^-_s(\textbf{a})$ and the corresponding vector, a lowest weight vector for $M_s(\textbf{a})$.

  One may recall that $(47)$ is equivalent to all the factors in the left hand side of $(7)$ (or $(8)$) being strictly positive. Its proof was a delicate matter even in the rigid case (Lemma \ref {6.1}(ii)).  It cannot be expected to be easier in the general case.

  \subsubsection{}\label{7.5.2}

  Lemmas \ref {7.2.7}, \ref {5.2} combined with the induction argument in \ref {7.3}  (which overcomes the difficulty raised in \ref {7.2.9}) allows us to compute the proportionality coefficients between the vectors $\overline{v}_{\textbf{l}'}$ corresponding to the trails in $T_s^-(\textbf{a})$ - these all being proportional to the lowest weight vector of the simple module $M_s(\textbf{a})$.  These proportionality factors are all binomial coefficients and hence strictly positive.

   Assume that $\overline{v}_{\textbf{k}}$ is proportional to a weight vector of $M_s(\textbf{a})$. In principle one may use Lemma \ref {3.2} to compute the coefficient of $v_{\textbf{l}'}$ in $f^bv_{\textbf{k}}$, use the above proportionality factors to compute the coefficient of $f^b\overline{v}_{\textbf{k}}$ of the lowest weight vector in $M_s(\textbf{a})$ and then show that this expression is non-zero only if the upper bounds are satisfied. Though possible in principle, this calculation could be quite horrendous involving complicated sums with both positive and negative terms.

   We provide a short-cut to the above computation by using the observation in \ref {6.5}.  The method will first be illustrated in the case $n=3$.  Here we adopt the notation of \ref {6.5} and note that failure of rigidity means that $q:=\ell_2+\ell_1-a_2 >0$.

   For each trail $L^{i} \in T_s^-(\textbf{a})$, let $\overline{v}_{\textbf{l}^{(i)}}$  be the (non-zero) element in $M_s(\textbf{a})$ it defines.  These vectors are all proportional to the unique up to scalars lowest weight vector in $M_s(\textbf{a})$, which recalling Lemma \ref {5.3.4}, we can take to be $\overline{v}_{\overline {\textbf{l}}}$.  Thus there exist non-zero scalars $h_i$ such that $\overline{v}_{\textbf{l}^{(i)}}=h_i\overline{v}_{\overline {\textbf{l}}}$.  We assume that the conclusion of Proposition \ref {7.3}(i) holds. This provides the hypothesis of Lemma \ref {5.2} and by its conclusion we may compute the $h_i$ up to an overall non-zero scalar.

   Let $\overline{v}_\textbf{k}$ be an element of the simple module $M_s(\textbf{a})$ of the form given in $(16)$.
   %Yet $v_\textbf{k}\neq 0$ if and only if there exists $n \in \mathbb N$ such that $f^nv_\textbf{k}$ is a non-zero multiple of $v_{\overline {\textbf{l}}}$.
   We may lift $\overline{v}_\textbf{k}$ to the corresponding element of the tensor product as $v_\textbf{k}:=e^{k_3}v_{-a_3}\otimes e^{k_2}v_{-a_2}\otimes e^{k_1}v_{-a_1}$ and compute the coefficient say $g_i$ of $v_{\textbf{l}^{(i)}}$, the latter again viewed as an element of the tensor product, using Lemma \ref {3.2}.  This gives $f^b\overline{v}_\textbf{k}= (\sum g_ih_i)\overline{v}_{\overline {\textbf{l}}}$.

   Now set $\tilde{\textbf{l}}=(\tilde{\ell_3},\tilde{\ell_2},\tilde{\ell_1})$, with $\tilde{\ell_1}=\ell_1, \tilde{\ell_2}=\ell_2-q,\tilde{\ell_3}=0$. Then the observation of \ref {6.5} is that $h_i$ is also the coefficient of $f^qv_{\textbf{l}^{(i)}}$ in $v_{\tilde{\textbf{l}}}$ (up to a non-zero scalar independent of $i$).  Consequently, up to this non-zero scalar, $\sum g_ih_i$ is the coefficient of $v_{\tilde{\textbf{l}}}$ in $f^{q+b}v_\textbf{k}$, which can again be computed from Lemma \ref {3.2}.  We wish to know when it is non-zero.

   Now for a fixed choice of $\textbf{k}$ we may compare $\tilde{\textbf{c}},\tilde{\textbf{c}}'$ computed relative to $\tilde {\textbf{l}}$, with $\textbf{c},\textbf{c}'$ computed relative to $\textbf{l}$ by $(26),(30)$.  Since $\tilde{\ell_1}=\ell_1$, we obtain $\tilde{c_1}=c_1$ and  $\tilde{c_1}'=c_1'$.  On the other hand $\tilde{c_2}-c_2=\tilde{c_2}'-c_2'=q$. We conclude that $\tilde{c_i} - \tilde{c_i}'=c_i-c_i':i=1,2$.

   Finally observe that since $\tilde{c_2}=\tilde{c_1}$ by construction, the system they define is rigid.  Thus by Lemma \ref {6.1}(ii), using $(34)$, we conclude that the coefficient of $v_{\tilde{\textbf{l}}}$ in $f^{q+b}v_\textbf{k}$ is non-zero only if $\tilde{c_i} \geq \tilde{c_i}':i=1,2$.

    We conclude that $f^b\overline{v}_\textbf{k}$ that is non-zero multiple of $\overline{v}_{\overline {\textbf{l}}}$ only if $c_i\geq c_i':i=1,2$.  Yet by the simplicity of $M_s(\textbf{a})$, this implies that $\overline{v}_\textbf{k}\neq 0$ only if $c_i\geq c_i':i=1,2$.  This means that the upper bounds on the trail $K\in T_s(\textbf{a})$ corresponding to $v_\textbf{k}$ are satisfied,  as required.

   Notice that in this proof we used the fact that $T^-_s(\textbf{a})$ is given by the conclusion of Proposition \ref {7.3}(i) to obtain the proportionality factors between the corresponding vectors $\overline{v}_{\textbf{l}^{(i)}}$ via Lemma \ref {5.2}.  On the other hand we have \textit{not} used the hypothesis of Corollary \ref {6.2} with respect to the tilde variables, which we have no reason to assume true.

   \subsubsection{}\label{7.5.3}

   The proof of the upper bounds for arbitrary values of $n$ follows a similar path.

   First we construct a rigid system $\tilde{c}_i:i\in N$ from the given set $c_i:i\in N$ of coefficients.

   %First we define the $\tilde {l_i}$ inductively as follows.
   Set $\tilde{\ell_1}=\ell_1$ and define $\tilde{\ell_i}:i \in N$ inductively through $\tilde{\ell_i}= \min{(\ell_i,a_i-\tilde{\ell}_{i-1})}$, for $i >1$.  Notice that this formula also holds for $i=1$ if we take $\tilde{\ell_0}=0$, since $\ell_1 \leq a_1$ because $v_{-a_1}$ is a lowest weight vector.  Obviously this implies that $\tilde{\ell}_{i-1} \leq \ell_{i-1}$.  Then by Lemma \ref {5.3.1} we obtain $a_{i}-\tilde {\ell}_{i-1}\geq 0$.  Since $\ell_i \geq 0$, we obtain that $\tilde{\ell_i}\geq 0$.  In particular $v_{\tilde{\textbf{l}}}$, as an element of the tensor product, is defined.

   By construction $a_i \geq \tilde{\ell_i}+\tilde{\ell}_{i-1}$, for all $i \in N$.  Thus if we define $\tilde{c_i}:i \in N$ by $(26)$ using $\tilde{\ell}_i$ instead of $\ell_i$, it follows that the $\tilde{c_i}$ are increasing, that is form a rigid system.

   Similarly we define $\tilde{c_i}':i \in N$ by $(30)$ using $\tilde{\ell}_i$ instead of $\ell_i$.

   Now $(c_i-c_i')-(c_{i-1}-c'_{i-1})=a_i-k_i-k_{i-1}$, so is unchanged if we replace $\ell_i$ by $\tilde{\ell}_i$.   Consequently
   $$((\tilde{c}_i-\tilde{c}_{i-1})-(\tilde{c}_i'-\tilde{c}'_{i-1}))
   -((c_i-c'_{i-1})-(c_{i-1}-c'_{i-1}))=0, \forall i \in N\setminus \{1\} . \eqno {(48)}$$

   We conclude that
   $$\tilde{c}_i-\tilde{c}_i'=c_i-c_i', \forall i \in N. \eqno {(49)}$$

   Indeed since it holds trivially for $i=1$, it holds for all $i \in N$ through $(48)$ by induction on $i$.

   Finally set $\tilde{\ell}_n=0$ and $\tilde{\textbf{l}}=(\tilde{\ell}_n,\ldots,\tilde{\ell}_1)$.  We may suppose that $\ell_n=0$ by replacing $\textbf{l}$ by $\overline{\textbf{l}}$.

   \subsubsection{}\label{7.5.4}

   The second step is to generalize the non-vanishing of the common factor which we deleted to obtain $(36)$.  Set $b_j=\ell_j-\tilde{\ell}_j: j \in N$ and $b=\sum_{j \in N}b_j$.  This will eventually be proved (\ref {7.5.6}) through the following

   \begin {lemma}  The coefficient of $v_{\tilde{\textbf{l}}}$ in $f^bv_{\overline{\textbf{l}}}$ is non-zero (actually positive).
   \end {lemma}

   \begin {proof}  We have only to show that each factor in $C_\textbf{b} (\overline{\textbf{l}},\tilde{\textbf{l}})$ is positive.  For this it is enough that $a^{(j)}-\ell^{(j)}-\tilde{\ell}^{(j-1)}\geq 0$, for all $j\in N$.  Yet this expression is the sum of $c_j$ and $\ell^{(j-1)}-\tilde{\ell}^{(j-1)}$.  The first is non-negative by Lemma \ref {5.3.2}.  The second is non-negative by construction (see \ref {7.5.3}).
   \end {proof}

    \subsubsection{}\label{7.5.5}

    The third step is to generalize the comparison of coefficients described in $(36)$ when we modify factors in a sub-expression.  Specifically take a trail $L' \in T^-_s(\textbf{a})$ and let $v_{\textbf{l}'}$ be the corresponding lowest weight vector in $M_s(\textbf{a})$.  We know these to be all proportional; but we want to compute the proportionality factors by the method used in \ref {6.5}.

     As in \ref {5.2} consider modifying the exponents of a subexpression $<u_1,u_2>$ occurring at position $i+1\in N$ in $v_{\textbf{l}'}$, that is to say with $a_{i+1}=a$.  Specifically assume that $q:=\ell_{i+1}'+\ell_i'-a_{i+1}>0:\ell'_{i+1}-q\geq 0$ and for all $v \in [0,q]$, define the $n-1$-tuple $\textbf{l}(v)$ by setting $\ell(v)_{i+1}=\ell_{i+1}'-v,\ell(v)_{i}=\ell_i'+v$ and $\ell(v)_{j}=\ell'_{j}:j \in N\setminus \{i+1,i\}$.

     Using $(7)$ one may calculate the coefficient of $v_{\tilde{\textbf{l}}}$ in $f^bv_{\textbf{l}(v)}$ (with $b$ as in \ref {7.5.4}).  Using the relation $\ell'_{i+1}+\ell'_i-q=a_{i+1}$, the terms coming from $C_\textbf{b}(\textbf{l}(v),\tilde{\textbf{l}})$ with $j=i,i+1$ collapse into the single product.  Then we obtain this coefficient to be
$$b!{\ell'_{i+1}-v \choose \ell'_{i+1}-q}{\ell_i'+v \choose \ell_i'}\prod_{j=1}^q(a^{(i)}-j+1 -\tilde{\ell}^{(i)}-\ell(v)^{(i-1)}),$$
up to common factors given by $B_\textbf{b} (\textbf{l}(v),\tilde{\textbf{l}})$ (which are simply binomial coefficients and so positive) and by $C_\textbf{b} (\textbf{l}(v),\tilde{\textbf{l}})$, with $j\neq i,i+1$.

Thus up to an overall common factor we can take the coefficient $h_v$ of $v_{\tilde{\textbf{l}}}$ in $f^bv_{\textbf{l}(v)}$ to be given by
$$h_v=(\ell'_{i+1}-v)!(\ell'_{i}+v)! { q \choose v}. \eqno {(50)}$$

%For the moment we do not know that this common factor is non-zero, yet of course this condition is equivalent to $v_{\tilde{\textbf{l}}}$ occurring in $f^bv_{\textbf{l}'}$, with a non-zero coefficient.  By Lemma \ref {7.5.4}, this does hold if $\textbf{l}'=\overline{\textbf{l}}$.

For the moment we do not know that this common factor is non-zero, yet of course this condition is equivalent to $v_{\tilde{\textbf{l}}}$ occurring in $f^bv_{\textbf{l}'}$, with a non-zero coefficient.  By Lemma \ref {7.5.4}, this does hold if $\textbf{l}'=\overline{\textbf{l}}$.

The non-vanishing of these common factors will be established in the next section.

\subsubsection{}\label{7.5.6}
The fourth step is to show that the coefficients $h_i$ described by $(50)$ are those obtained from Lemma \ref {5.2} when the hypotheses of the latter are satisfied.  Notice here that if $q \leq 0$ then there is nothing to prove and furthermore the condition $\ell_{i+1} -q \geq 0$ (used in \ref {7.5.5}) was also assumed in Lemma \ref {5.2}.

Lemma \ref {7.2.7} sets up a framework in which the hypotheses of Lemma \ref {5.2} are satisfied.  It is not immediately applicable because of the difficulty in the induction step noted in \ref {7.2.9}. However this difficulty was overcome in the proof of Proposition \ref {7.3}.

Thus assume that the hypothesis of Proposition \ref {7.3}(i) holds.  Then by its conclusion (or rather by the construction given in its proof) we can find a sequence of trails $L^{(i)} \in T_s^-(\textbf{a})$, starting from the $\ell$-minimal trail $L^{(1)}$, such that any pair of successive trails in this sequence are extremal elements on a line of type $u \in N$ whose length exactly matches the conditions required to apply Lemma \ref {5.2}.  By Lemma \ref {5.2} the integer points of this line are trails and the corresponding vectors are related by the proportionality factors given by its conclusion.  These are exactly the proportionality factors given by $(50)$.  Furthermore \textit{every} trail in $T_s^-(\textbf{a})$ so obtains.  (Indeed this is what the conclusion of Proposition \ref {7.3}(i) states.)

Let $v_{\textbf{l}^{(i)}}$ be the lift in the tensor product corresponding to the vector $\overline{v_{\textbf{l}}}^{(i)}\in M_s(\textbf{a})$ defined by $L^{(i)}$.  Let us show that the coefficient of $v_{\tilde{\textbf{l}}}$ in $f^nv_{\overline{\textbf{l}}^{(i)}}$  is non-zero by induction. By Lemma \ref {7.5.4} this holds for $i=1$. Then as noted in the last paragraph of \ref {7.5.5} this means that the proportionality factor which we deleted on obtaining $(50)$ with respect to the pair $L^{(i)},L^{(i+1)}$ is non-zero, which in turn implies the coefficient of $v_{\tilde{\textbf{l}}}$ in $f^bv_{\overline{\textbf{l}}^{(i+1)}}$  is non-zero, as required.  Furthermore we also conclude that all these non-zero proportionality factors are equal.

We have thus proved the following
\begin {lemma}  For all trails $L^{(i)} \in T_s^-(\textbf{a})$ the coefficient $h_i$ of $v_{\tilde{\textbf{l}}}$ in $f^bv_{\overline{\textbf{l}}^{(i)}}$ satisfies $\overline{v}_{\textbf{l}^{(i)}}=h_i\overline{v}_{\overline {\textbf{l}}}$, up to a fixed non-zero scalar.
\end {lemma}

\subsubsection{}\label{7.5.7}

We may now conclude as in the special case described in \ref {7.5.2}.

   \begin {prop}  Suppose that the trails in $T_s^-(\textbf{a})$ satisfy the conclusion of Proposition \ref {7.3}(i).  Then the trails in $T_s(\textbf{a})$ satisfy the upper bounds.
   \end {prop}

   \begin {proof} %We fill in a few details.
   Take a vector $\overline{v}_\textbf{k}$ corresponding to a trail $K \in T_s(\textbf{a})$. Choose $b \in \mathbb N$ such that $f^b\overline{v}_\textbf{k}$ is a multiple of a lowest weight vector. The simplicity of $M_s(\textbf{a})$ implies that $\overline{v}_\textbf{k}\neq 0$ if and only if this multiple is non-zero.
   Set $q = \sum_{i\in N} \ell_i-\tilde{\ell}_i$.
   %Consider $f^nv_\textbf{k}$.
   Since the $\tilde{c_i}$ system is rigid, we conclude by Lemma \ref {6.1}(ii) that $f^{b+q}v_\textbf{k}$ is a non-zero multiple $d$ of $v_{\tilde{\textbf{l}}}$ only if $\tilde{c}_i'\leq \tilde{c}_i$, for all $i \in N$.

   %Let $\overline{v}_{\textbf{l}^{(i)}}$  be the non-zero element in $M_s(\textbf{a})$ defined by a trail $L^{(i)} \in T^-_s(\textbf{a})$. As in \ref {7.5.2} there exist non-zero scalars $h_i$ such that $\overline{v}_{\textbf{l}^{(i)}}=h_i\overline{v}_{\overline {\textbf{l}}}$.
%
%   We claim that $h_i$ is given, up to an overall non-zero scalar, by the coefficient of $v_{\tilde{\textbf{l}}}$ in $f^qv_{\textbf{l}^{(i)}}$.
%
%   By Proposition \ref {7.3}(i) one may recover every trail in $T_s^-(\textbf{a})$ as an element of $K^-_\mathbb Z(\textbf{c})$.  By Lemma \ref {7.2.6} this means the proportionality factors between the $\overline{v}_{\textbf{l}^{(i)}}$ are given by Lemma \ref {5.2}. Then by the observation in the last part of \ref {6.5} these proportionality factors are just those given by the claim.

   As in \ref {7.5.2} let $g_i$ denote the coefficient of $v_{\textbf{l}^{(i)}}$ in $f^bv_\textbf{k}$.  Then the coefficient of $\overline{v}_{\overline {\textbf{l}}}$ in $f^b\overline{v}_\textbf{k}$ is just $\sum_ig_ih_i$.  By Lemma \ref {7.5.6} this is just the coefficient $d$ of $v_{\tilde{\textbf{l}}}$ in $f^{b+q}v_\textbf{k}$, up to a non-zero scalar.  By the first paragraph above and $(49)$, this gives $c'_i\leq c_i$, for all $i \in N$, as required.

  % On the other hand the conclusion of Proposition \ref {7.3}(i) taken together with Lemma \ref {5.2} (which calculates the relative multiples $h_i$ of the $v_{\textbf{l}^{(i)}}$) and the observation in \ref {6.5} (which shows that the $h_i$ are given by the coefficient of $v_{\tilde{\textbf{l}}}$ in $f^qv_{\textbf{l}^{(i)}}$) implies exactly as in \ref {7.5.2} that $f^nv_\textbf{k}=dv_{\overline{\textbf{l}}}$, also.  Thus $v_\textbf{k}\neq 0$ implies $d \neq 0$, which by the first paragraph and $(49)$ gives the required upper bounds.
   \end {proof}

   \subsubsection{}\label{7.5.8}

   One can ask if there is a more direct way to obtain the upper bounds for trails in $T_s(\textbf{a})$ given those for trails in $T_s^-(\textbf{a})$.  This does work in the rigid case and ultimately gives a second proof of Theorem \ref {6.4}.

   Let $\overline{v}_{\textbf{k}}$ be a non-zero vector in $M_s(\textbf{a})$ corresponding to a trail $K \in T_s(\textbf{a})$ for which the upper bound is not satisfied.  We can assume that $\overline{v}_{\textbf{k}}$ has lowest possible weight with this property.  By the hypothesis $\overline{v}_{\textbf{k}}$ is not proportional to a lowest weight vector of the simple module $M_s(\textbf{a})$ and so there exists $b \in \mathbb N^+$ such that $f^b\overline{v}_{\textbf{k}}$ is a non-zero multiple of this lowest weight vector.

   Consider $f\overline{v}_{\textbf{k}}$. It is non-zero. The action of $f$ gives a sum of terms $\overline{v}_{\textbf{k}'}$ defined by reducing just one exponent of $e$ in $\overline{v}_{\textbf{k}}$ by $1$. At least one such term $\overline{v}_{\textbf{k}'}$ must be non-zero.  Moreover any such non-zero term  corresponds to a trail $K' \in T_s(\textbf{a})$ and satisfies the upper bound by the choice of $\overline{v}_{\textbf{k}}$.

   Suppose it is the $j'^{th}$ exponent of $e$, namely $k_{j'}$ of $\overline{v}_{\textbf{k}'}$, which  is one less than that of $\overline{v}_{\textbf{k}}$. Then
   $$ c_j'^{K'}=\left\{
                       \begin{array}{ll}
                         c_j'^{K}, & \hbox{if $j <j'$}, \\
                         c_j'^{K}-1, & \hbox{if $j\geq j'$.}
                       \end{array}
                     \right. $$

    From this it follows that the condition that $K'$ satisfies the upper bound, whilst $K$ does not, forces $c_j'^{K'}\leq c_j$ for all $j \in \hat{N}$ with equality for some $j \geq j'$. Yet for $j \geq j'$ we obtain from $(26),(30),(48)$ that
    %Now  for $j \geq j'$ one has
    %$c_j-c_j'^{K'}=a^{(j)}-\ell^{(j)}-\ell^{(j-1)}-(k^{(j)}-\ell^{(j)}-1)$ and so
    %there exists $j \in \hat{N}$ with $j \geq j'$ such that
   $$c_j-c_j'^{K'}=c_j-(c_j'^K -1)=a^{(j)}+1-k^{(j)}-\ell^{(j-1)}. \eqno {(51)}$$

   % We claim that for every choice $v_{\textbf{l}'}$ for which the corresponding trail $L'$ lies in $T_s^-(\textbf{a})$, the coefficient $A_n(\textbf{k},\textbf{l}')$ of $v_{\textbf{l}'}$ in $f^nv_\textbf{k}$, \textit{as computed in the tensor product}, is zero.  This will imply that $f^nv_\textbf{k}=0$, which is a contradiction.

   Recall that any $\overline{v}_{\textbf{l}^{(i)}}$ for which the corresponding trail $L^{(i)}$ lies in $T_s^-(\textbf{a})$, is a non-zero multiple $h_i$ of the lowest weight vector of $M_s(\textbf{a})$, which can be taken to be $\overline{v}_{\overline{\textbf{l}}}$, that is we may write $\overline{v}_{\textbf{l}^{(i)}}=h_i\overline{v}_{\overline{\textbf{l}}}$.

   We would like to show that the resulting coefficient of $\overline{v}_{\overline{\textbf{l}}}$    in $f^b\overline{v}_\textbf{k}$, namely $\sum_ig_ih_i$, where $g_i= A_\textbf{b} (\textbf{k},\textbf{l}^{(i)})$, is zero.   This would imply that $f^b\overline{v}_\textbf{k}=0$, which is a contradiction.  %This is what we did in \ref {7.5.7}, though it would be pleasant to be able to obtain this result from $(51)$.

    \

    Let us show that the above method works in the rigid case. In this case $\{\textbf{l}^{(i)}\}$ is reduced to $\textbf{l}$. This forces $k'_i \geq \ell_i$, for all $i \in \hat{N}$.

     Choose $j \geq j'$ minimal such that $c_j=c_j'^{K'}$ and suppose $k_j=\ell_j$. Since $k_{j'}=k'_{j'}+1$, we conclude that $j>j'$.

     Then $(c_j-c_j'^{K'})-(c_{j-1}-c_{j-1}'^{K'})= a_j-\ell_j-\ell_{j-1}-(k_j-\ell_j)$, which by rigidity is non-negative.  Thus    $c_{j-1}-c_{j-1}'^{K'}\leq 0$.  On the other hand $c_{j-1}\geq c_{j-1}'^{K'}$, so equality holds and this  contradicts the choice of $j$, proving that $k_j>\ell_j$, for this choice.

     Then inspection of $(7)$ shows that the right hand side of $(51)$ is a factor of $A_\textbf{b} (\textbf{k},\textbf{l})$ and its vanishing establishes what we require.

 \subsection{}\label{7.6}

 We extend Theorem \ref {6.4} to the general case.

 \

 \textbf{Definition.} We say that $T_s^-(\textbf{a})$ (resp. $T_s(\textbf{a})$) possess no false trails if the hypothesis of Corollary \ref {7.2.8}(i) (resp. (ii)) holds.  If the $c_j:j \in N$ are increasing, the former just means that $T_s^-(\textbf{a})$ is reduced to the $\ell$-minimal trail $L$.

 \begin {thm} Assume that $T^-_s(\textbf{a})$ admits no false trails.  Then $\{z^K\}_{K\in T_s(\textbf{a})}$ is the set $K_\mathbb Z(\textbf{c})$ of all integer points of the convex set whose extremal elements form the $S$-set $Z(\textbf{c})$ defined by $\textbf{c}:= \{c_j^{K_{\ell \min}}(s)\}_{j=1}^{n-1}$.  Moreover $T_s(\textbf{a})$ admits no false trails.
  \end {thm}

  \begin {proof}  Apply Lemmas \ref {7.4} and \ref {7.5.7} to Proposition \ref {7.3}(ii).
  \end {proof}

  \section{From the Absence of False Trails to Giant $S$-sets}\label{8}

  Fix a sequence $J$ of elements of $I$ as in \ref {1.3}.

  In this section an $S$-set of type $s$ will always mean the set of functions defined the vertices of the canonical $S$ graph $\mathscr G(\textbf{c})$, specified by a coefficient set $\textbf{c}$ deduced from an $\ell$-minimal trail following \ref {5.3.5}.

  \subsection{}\label{8.1}

  Fix $t \in I$. Recall (\ref {2.3}) the initial driving function $z_t^1$ associated to $t$. Recall (\ref {2.4}) the definition of the set $X_t$ of locally finite linear functions on $B_J$.

  \textbf{Definition.} A giant $S$-set (relative to $J$) associated to $t$ is a subset $Z_t$ of $X_t$  such that $Z_t\setminus \{z_t\}$ is a disjoint union of $S$-sets of type $t$ and for $s\in I\setminus \{t\}$ is a disjoint union of $S$-sets of type $s$.

 \subsection{}\label{8.2}

 %We cannot and do not specify when a trail in $\mathscr K^{BZ}_t$ is false.  Rather we say what it means for $\mathscr K^{BZ}_t$ to possess no false trails, or equivalently that false trails are absent from $\mathscr K^{BZ}_t$.    Recall the notation of \ref {4.3}.

We describe what it means for $\mathscr K^{BZ}_t$ to possess no false trails.    Recall the notation of \ref {4.3}.

 Fix $s \in I$ and recall the notation of \ref {4.2} and \ref {4.3}.
 %We say that $T_s(\textbf{a})$, more properly $T_s(\textbf{e})$ (resp. $T^-_s(\textbf{e})$, more properly $T^-_s(\textbf{e})$) possesses no false trails if the hypothesis of (i) (resp. (ii)) of Proposition \ref {7.3} holds.

 %The conclusion of Theorem \ref {7.6} implies that these two properties are equivalent, so we only have to verify that $T^-_s(\textbf{e})$ possesses no false trails.

 Let $\textbf{E}_s$ denote all choices of tuples $\textbf{e}$ of products of simple root vectors different from $e_s$, such that $T_s(\textbf{e})$ is non-empty, that is to say defined by all trails which trivialise at some $w_j:j \in J$, with $i_j=s$.  These sets are disjoint and their union is the set $\mathscr K^{BZ}_t$ of all trails in $V(-\varpi_t)$.  Each $T_s(\textbf{e})$ admits a unique $\ell$-minimal trail (which actually lies in $T^-_s(\textbf{e})$, equivalently trivializes at $w_{j-1}$) from which one may compute a tuple $\textbf{c}$ of non-negative coefficients (5.3.3) and hence an $S$-set $Z(\textbf{c})$.
 %, which may be more appropriate written as $Z_\textbf{e}$.

   Recall that $T_s(\textbf{e})$ (resp. $T^-_s(\textbf{e})$) possesses no false trails if the hypothesis of (i) (resp. (ii)) of Proposition \ref {7.3} holds.  (Here we recall that the convention of \ref {4.3} was used in \ref {7.3}.)

 The conclusion of Theorem \ref {7.6} implies that these two properties are equivalent, so we only have to verify that $T^-_s(\textbf{e})$ possesses no false trails.

  \

\textbf{Definition.}

 Fix $s \in I$. We say that $\mathscr K^{BZ}_t$ possesses no false trails relative to $s$ if $T_s(\textbf{e})$, equivalently $T^-_s(\textbf{e})$, possesses no false trails for all $\textbf{e} \in \textbf{E}_s$.

  We say that $\mathscr K^{BZ}_t$ possesses no false trails, if $V(-\varpi_t)$ possesses no false trails relative to $s$, for all $s \in I$.

  \

  To simplify presentation we shall in the remainder of this section identify a trail $K \in \mathscr K^{BZ}_t$ with the linear function $z^K$ it defines on $B_J$ given in \ref {2.3}.

  %\textbf{Remark}.  This in particular means that some knowledge of the trails which trivialize at $s_\alpha w$ determines those which trivialize at $w$.  Obviously we would like to show that the hypothesis of the above theorem obtains by induction.  This is not straightforward.  However its truth would imply that the set of all trails form the integer points of a convex set which for every $i \in I$ is the convex hull of a disjoint union of $S$-sets of type $i$.  It is also not quite obvious that the extremal points of this convex hull is the corresponding disjoint union of the $S$-sets of type $i$ (thereby proving the independence of this union on $i$). Yet even without this second point we may still conclude that the maximum of all functions determined by the trails in the fundamental module $V(-\varpi_i)$ determines the function $\varepsilon^\star_i$ on $B_J(\infty)$.

   %By our assumption on $v_{\textbf{k}}$, there must be exactly one place

 \subsection{}\label{8.3}

 The set $K_\mathbb Z(\textbf{c})$ (resp. $Z(\textbf{c})$) introduced in \ref {7.1} also depends on $s \in I$ and the choice of $\textbf{e} \in \textbf{E}_s$.  On the other hand as noted in \ref {5.3.5}, $\textbf{e}$ determines $\textbf{c}$ and so this set should be more properly written as
 $K_\mathbb Z(\textbf{e})$ (resp. $Z(\textbf{e})$).

  As in \ref {4.2}, let $\mathfrak s$ denote the $\mathfrak {sl}(2)$ subalgebra of $\mathfrak g$ spanned by the triple $(e_s,h_s,f_s)$.  We write the simple $\mathfrak s$ module $M_s(\textbf{e})$ introduced in \ref {4.3} simply as $M(\textbf{e})$.

 In the above the dependence on $s \in I$ is implicit in the fact that $\textbf{e} \in \textbf{E}_s$.

 These sets also depend on the choice of $t \in I$; but since $t$ is fixed throughout, this dependence is not indicated.  %Often we omit the superscript $s$, since it is implicit in the statement that $\textbf{e} \in  \textbf{E}_s$.

 When $s=t$ we augment $\textbf{E}_s$ by an element $\textbf{e}_\phi$ such that $Z_t(\textbf{e}_\phi)=\{z_t^1\}$, which is itself \textit{not} an $S$-set.

 From the discussion in \ref {8.2}, which basically sets up notation, we immediately obtain the following consequence of Proposition \ref {7.6}.

 \begin {cor}  Suppose that $\mathscr K^{BZ}_t$ possesses no false trails.  Then
 $\mathscr K_t(s):=\sqcup_{\textbf{e}\in \textbf{E}_s} K_\mathbb Z(\textbf{e})=\mathscr K^{BZ}_t$.
 In particular the left hand side is independent of $s \in I$.
 \end {cor}

 \textbf{Remark}.  In view of this independence we write the left hand side simply as $\mathscr K_t$.  It is called the giant $S$-set envelope associated to $t \in I$.  Under the hypothesis of the Corollary, one has $\mathscr K_t^{BZ}=\mathscr K_t$.

 \subsection{}\label{8.4}

 Under its hypothesis, Corollary \ref {8.3} gives an inductive procedure for computing $\mathscr K^{BZ}_t$, or equivalently the $T_s(\textbf{e}):s \in I, \textbf{e} \in \textbf{E}_s$.  This goes as follows and here we remark that it involves all the elements of $s \in I$.

 %To start the procedure, let $j \in J$ be minimal such that there exists a trail $K$ which trivializes at $w_j$.  Since $K_t^1$ trivializes at $w_{(t,1)}$ - see \ref {2.3} - it follows that $\gamma^K_{j+1}=-w_j\varpi_t$ if $j \geq (t,1)$. In particular $\gamma^K_{(t,1)+1}=-w_{(t,1)}\varpi_t=-s_t\varpi_t$. %Now suppose $K$ trivializes at some $j<(t,1)$.  We can write $j=(s,k)$, for some $s \in I\setminus \{t\}, k \in \mathbb N^+$.
% Since $-s_t\varpi_t$ is the unique minimal weight of $\Omega(V(-\varpi_t))\setminus\{-\varpi_t\}$, it follows that $\gamma_j^K=-s_t\varpi_t$, for all $j \leq (t,1)$.  Hence $K=K_t^1$.

 Suppose by way of induction we have constructed all the trails $K \in \mathscr K^{BZ}_t$ which trivialize at some $w_j:j\in J$.   By Lemma \ref {2.2} we may assume $j \geq (t,1)$.  Moreover when equality holds $K=K_t^1$, so the induction starts.

 Set $s=i_{j}$.  Then the trails which trivialize at $w_{j}$ are defined by a finite subset $\textbf{E}^j_s$ of $\textbf{E}_s$.
  %of tuples compatible with $J$ of products of the simple roots vectors different from $e_s$.
  For each tuple $\textbf{e}\in \textbf{E}^j_s$, there is defined an $\ell$-minimal trail which trivializes at $w_{j-1}$ given in particular by a tuple $\textbf{l}$, which can be computed from $\textbf{e}$.  Then a tuple $\textbf{c}$ of coefficients can be computed from $\textbf{e}$ using $(26)$ and the corresponding $S$-set $Z(\textbf{c})$ of type $s$ written down.

 Following \ref {8.2}, we write $Z(\textbf{c})$ as $Z(\textbf{e})$ and we write its driving function (defined by its unique $\ell$-minimal trail) by $z_\textbf{e}$.

 Assume the hypothesis of Corollary \ref {8.3}.  Then by Proposition \ref {7.6}, the $S$-set $Z(\textbf{e})$, determines the trails defined by $\textbf{e}$ as its $\mathbb Z$ convex hull $K_\mathbb Z(\textbf{e})$. As $\textbf{e}$ runs over $\textbf{E}_s^j$ this constructs all the trails which trivialize at $w_{j}$.

 This procedure is actually much easier to implement than the description would suggest.  This is because we only have to compute the functions $z^K:K\in \mathscr K^{BZ}_t$.  As we have already noted in \ref {2.3}, $K$ trivializes at $w_j$ if and only if the coefficient  of $m_{j'}$ in $z^K$ is zero for all $j'>j$, that is if $z^K$ has support in $[1,j]$ in the sense of $(2)$.

 With $s$ as above, a function $z^K$ attached to an $\ell$-minimal trail $K$ of type $s$ has coefficients $-c_k$ of the $m_{(s,k)}: k \in \mathbb N^+|(s,k)\leq j$. Let $\textbf{c}=(c_1,c_2,\ldots)$ be the resulting (finite) tuple of integers.  Notably the $c_k:k \in N$ are non-negative (by Lemma \ref {5.3.2}).

  It is slightly unfortunate that the converse of this last statement is false. In other words there may be a trail $K$ which is not $\ell$-minimal, for which the corresponding function $z^K$ has non-positive coefficients of the $m_{(s,k)}:k \in \mathbb N^+$. (An example was found with the help of S. Zelikson for a certain reduced decomposition in type $D_5$.)  In this case the resulting tuple $\textbf{c}$ should \textit{not} be used to construct an $S$-set (of type $s$).  Some tuples can be detected and then discarded by the following algorithm.

  \

  $(A)$.  Compute the $S$-set of type $s$ for each tuple $\textbf{c}$ of non-negative coefficients obtained as above from the trails which trivialize at $w_j$. Recall the partial order $\leq$ on $N_t$ defined in \ref {2.5}.  With respect to $\leq$, the driving function of an $S$-set of type $s$ defined by $\textbf{c}$ is its unique minimal element.  Then using induction on this partial order discard those tuples $\textbf{c}$ when the corresponding function lies in a previous $S$-set of type $s$.

  \
 %It is easy to see that this is the correct recipe, though it adds a little tedium to the procedure.

 %It is also easy to see (given the conclusion of Corollary \ref {8.3}) that the union of these $S$-sets will contain all those (earlier computed) functions having support in $[j,1]$, that is to say there are ``enough'' functions with non-positive coefficients of the $m_{(s,k)}$.

% Of course we can apply this recipe using $S$-sets in complete ignorance of trails. However without the hypothesis and conclusion of Proposition \ref {7.6}, we cannot be sure that each step is well-defined, that is to say in the set of functions which have support in $[j,1]$, there are enough (or any!) having all their coefficients of $m_{(s,k)}$ non-positive (with $s$ as above). Nor can we be sure that the remarkable independence property of Corollary \ref {8.3} will hold with respect to the corresponding functions $z^K:K \in \sqcup_{e \in \textbf{E}_s}K_\mathbb Z(\textbf{e})$.  In our examples this was a particularly tedious property to verify by hand.

 %Recall that for each  $K \in \mathscr K^{BZ}_t$ there is a linear function $z^K$ on $B_J$.

 \subsection{}\label{8.5}

 Let $\mathbb Q^+$ denote the set of non-negative rational numbers. Following \ref {2.4} we set
 $$(X_t)_{\mathbb Q^+}= z_t^1+ \sum_{(s,k) \in I \times \mathbb N}\mathbb Q^+(r_s^k-r_s^{k+1}),$$
 in which sums are viewed as being finite. In this we may regard the $(r_s^k-r_s^{k+1}):s \in I, k \in \mathbb N^+$ as co-ordinate functions.  Then we may regard $X_t$ as being the set of integer points of $(X_t)_{\mathbb Q^+}$.

 Observe that $(X_t)_{\mathbb Q^+}$ is closed under taking (finite) convex $\mathbb Q$-linear combinations.

 %One may define a filtration $\mathscr F$ on $(N_t)_{\mathbb Q^+}$ by setting
%  $$\mathscr F^k((N_t)_{\mathbb Q^+})= z_t^1+ \sum_{s \in I, \ell \leq k}\mathbb Q^+(r_s^k-r_s^{k+1}).$$
%
%  One may view $\mathscr F^k((N_t)_{\mathbb Q^+})$ as lying in an affine translate of the finite dimensional vector space $\mathbb Q^{k|I|}$.  It is closed under taking convex linear combinations. Again suppose $z \in \mathscr F^k((N_t)_{\mathbb Q^+})$ and is a (finite) convex linear combination of elements $z_i \in (N_t)_{\mathbb Q^+}:i \in F$, (that is $z=\sum_{i \in F}a_iz_i$,  where $F$ is finite and $a_i \in \mathbb Q^+$, $\sum_{i \in F}a_i=1$), then positivity of coefficients ensures that $z_i \in \mathscr F^k((N_t)_{\mathbb Q^+})$, for all $i \in F$.

  Given a subset $S \subset X_t$, let $K(S)$ denote its convex hull in $(X_t)_{\mathbb Q^+}$, that is to say all finite convex linear combinations of its elements, and define $K(S)_\mathbb Z:=K(S) \cap X_t$ to be the $\mathbb Z$ convex hull of $S$.  Let $E(S)$ denote the set of extremal points of $S$, that is to say the elements of $S$ which cannot be written as finite convex linear combinations of elements of $S$.  Obviously $K(E(S))=K(S)$.

  \subsection{}\label{8.6}

  One may interpret Lemma \ref {8.3} in the language of \ref {8.5}.
  \begin {cor}  Assume that $\mathscr K^{BZ}_t$ possesses no false trails.  Then $\mathscr K^{BZ}_t\subset N_t$. Furthermore
  $$E(\mathscr K^{BZ}_t) \subset \sqcup_{\textbf{e} \in \textbf{E}_s}Z(\textbf{e})=:Z_t(s). \eqno {(52)}$$
  \end {cor}

  \begin {proof}  For all $s \in I$ and all $\textbf{e} \in \textbf{E}_s$ one has  $K_\mathbb Z(\textbf{e}) \subset N_t$, , which gives the first part, whilst by \cite [Thm. 1.4]{J5} one has $E(K_\mathbb Z(\textbf{e}))=Z(\textbf{e})$, which gives the second part.
  \end {proof}

  \subsection{}\label{8.7}

  In this section $b$ denotes an element of $B_J(\infty)$ and not an integer as in Sects. \ref {3} - \ref{7}.

  An immediate consequence of Corollaries \ref {8.3} and \ref {8.5} is the following

 \begin {thm}  Assume that $\mathscr K^{BZ}_t$ possesses no false trails. Then for all $b \in B_J(\infty)$ and all $s \in I$ one has
 $$\max_{z \in Z_t(s)}z(b) = \max_{z \in \mathscr K^{BZ}_t}z(b). \eqno {(53)}$$

 In particular the left hand side is independent of $s\in I$.
 \end {thm}

  %Assume that $\mathscr K^{BZ}_t$ possesses no false trails.  Then $\mathscr K^{BZ}_t\subset N_t$. Furthermore
%  $$E(\mathscr K^{BZ}_t) \subset \sqcup_{\textbf{e} \in \textbf{E}_s}Z_\textbf{e}=:Z_t(s). \eqno {(52)}$$

  %We would like to obtain equality in $(52)$ since then $Z_t(s)$ is independent of $s \in I$ and hence a giant $S$-set.  Even so this is refinement is unnecessary to conclude

   \textbf{Remarks.}  The conclusion of this theorem combined with the results in \cite {JZ}, that the dual Kashiwara parameter $\varepsilon^\star_t(b)$ on $B_J(\infty)$ is given by $\varepsilon^\star_t(b)=\max_{z \in \mathscr K^{BZ}_t}z(b)$, for all $b\in B_J(\infty)$.

   Up to our assumption of there being no false trails, this extends the result of Berenstein and Zelevinsky \cite [Thm. 3.9]{BZ} from the finite to the general case.

    \subsection{}\label{8.8}

    We would like to obtain equality in $(52)$ as this would mean that $Z_t(s)$ is independent of $s \in I$ and hence define a giant $S$-set $Z_t$.  Equivalently in the notation of Corollary \ref {8.3} we would have $Z_t=E(\mathscr K_t)$.  It is also an open question if $\mathscr K_t$ is the $\mathbb Z$ convex hull of $Z_t$ in the module of $\mathbb Q$-valued linear functions on the free $\mathbb N$ module $B_J$.

 \section{From Giant $S$-sets to the Absence of False Trails }\label{9}

 \subsection{}\label{9.1}  An immediate consequence of Theorem \ref {7.6} is that if its hypothesis  holds for \textit{all} $s \in I, \textbf{e}\in \textbf{E}_s$ then all the trials in $\mathscr K^{BZ}_t$ can be written down.

 Thus it remains to show that the hypothesis of Theorem \ref {7.6} is always satisfied but this is not so easy.  Nevertheless its conclusion provides a natural induction argument which we describe below.

 Let $\mathscr F^j(\mathscr K^{BZ}_t)$ denote the subset of $\mathscr K^{BZ}_t$ of trails which trivialize at $w_j$.  By Lemma \ref {2.2}, this set is empty if $j<(t,1)$ and reduced to $\{K_t^1\}$ when $j=(t,1)$.  Assume that $\mathscr F^{j-1}(\mathscr K^{BZ}_t)$ has been constructed by applying the conclusion of Theorem \ref {7.6}.  Set $s =i_j$.

% On the other hand it is clear from the conclusion of Theorem \ref {7.6}, that under its hypothesis the $E^j_{i_j}: j \in J$ are constructed inductively.  Thus we can inductively construct the functions $z^K: K \in \mathscr K^{BZ}_t$ with support in $[1,j]$.  In particular those with support in $[1,j-1]$ determine those with support in $[1,j]$.  Now with $s:=i_j$, the condition of Theorem \ref {7.6} is that the set of functions $z^K: K \in \mathscr K^{BZ}_t$ with support in $[1,j-1]$ form the disjoint union $\sqcup_{\textbf{e}\in \textbf{E}_s^j} K^-_\mathbb Z(\textbf{e})$.

  Then $z^K: K \in \mathscr F^{j-1}(\mathscr K^{BZ}_t)$ have support in $[1,j-1]$.

  To show that the hypothesis of Theorem \ref {7.6} is verified at the next induction step we have to verify that there exists a finite subset $\textbf{E}^j_s$ of $\textbf{E}_s$ such that
  $$\{z^K: K \in \mathscr F^{j-1}(\mathscr K^{BZ}_t)\}\subset \sqcup_{\textbf{e}\in \textbf{E}_s^j} K^-_\mathbb Z(\textbf{e}).\eqno{(54)}$$

  Then by Theorem \ref {7.6} we obtain
 $$\{z^K: K \in \mathscr F^{j}(\mathscr K^{BZ}_t)\}=\sqcup_{\textbf{e}\in \textbf{E}_s^j} K_\mathbb Z(\textbf{e}).\eqno{(55)}$$

 We shall say that a giant $S$-set associated to $t$ is constructible if $(54)$ holds for all $j \in J$ (starting in effect from $j =(t,1)$).
 %Notice that in this the algorithm $(A)$ given in \ref {8.4} is exactly what is involved in verifying the constructibility of a giant $S$-set.
 %One may observe that by definition, if a giant $S$-set is constructible then $\sqcup_{\textbf{e}\in \textbf{E}_s} K_\mathbb Z(\textbf{e})$ is independent of $s \in I$. In this holds $E(\sqcup_{\textbf{e}\in \textbf{E}_s} K_\mathbb Z(\textbf{e}))$ is contained in $Z_t(s)$ and is independent of $s \in I$.

 In view of $(55)$ we obtain inductively the following

 \begin {thm}  If a giant $S$-set associated to $t$ is constructible, then there are no false trails in $\mathscr K^{BZ}_t$.
 \end {thm}

 \textbf{Remark}.  Recall the notation of $(52)$.  In principle $Z_t(s)$ should also be independent of $s \in I$ and thus be the desired giant $S$-set $Z_t$.  The latter is a stronger condition and it is how we interpret the statement that a giant $S$-set exists.

  \subsection{}\label{9.2}

  Let us note what is involved in verifying $(54)$.

  The left hand side is supposed known through the induction hypothesis.

  %The advantage of the present formulation is that we do not have to consider trails to verify the constructibility of a giant $S$-set. Indeed the condition we have to verify, namely $(54)$ can be presented independent of the construction of a trail.  For this let us complete the discussion initiated in \ref {8.4}.

  The set $\textbf{E}_s^j$ appearing in right hand side is determined by the trails trivializing at $w_{j-1}$ and these are determined by the left hand side. Then as in \ref {8.4}, we may compute $K_\mathbb Z(\textbf{e})$ for all $\textbf{e} \in \textbf{E}_s$ and obtain $K^-_\mathbb Z(\textbf{e})$ as the subset of $K_\mathbb Z(\textbf{e})$ of functions with support in $[1,j-1]$.
  % For each $\textbf{e} \in \textbf{E}_s^j$ we may compute a tuple $\textbf{a}$ of non-negative integers and then from  corresponding $\ell$-minimal trail, a coefficient set $\textbf{c}$ using $(26)$. This allows one to compute an $S$-set $Z(\textbf{c})$, its subset $Z^-(\textbf{c})$ of elements with support in $[1,j-1]$ and the $\mathbb Z$ convex hull $K^-_\mathbb Z(\textbf{c})$ of the latter which by definition is just $K^-_\mathbb Z(\textbf{e})$.

  From the above data we may then verify $(54)$.
  %Notice that by (i) of Proposition \ref {7.3} we only have to show that the inclusion $\subset$ holds.

  Actually through the discussion in \ref {8.4} we do not have to compute trails to verify $(54)$, that is to say we do not have to compute the elements of the set $\textbf{E}_s^j$.  Indeed we may use the algorithm $(A)$ of \ref {8.4} to compute directly each coefficient set $\textbf{c}$. (In effect only the new coefficient sets $\textbf{c}$ for which not all elements of $Z(\textbf{c})$ have support in $[1,j-1]$ need be considered.)  From these coefficient sets the right hand side of $(54)$ may be computed and then (as noted above) we only have to verify that it contains the left hand side. Moreover we have only to verify that it contains the extremal elements of the left hand side which by the induction hypothesis and \cite [Thm 1.4]{J5} is contained in the union of the $S$-sets of type $s':=i_{j-1}$ obtained from the previous induction step, that is to say from $(55)$ with $j$ replaced by $j-1$.

  Summarizing we may rewrite Theorem \ref {9.1} in the following form

  \begin {thm}  Fix $j_1\in J$. Suppose that
  $$\cup_{\textbf{e}\in \textbf{E}_{i_j}^j}Z(\textbf{e}) \subset \cup_{\textbf{e}\in \textbf{E}_{i_{j+1}}}K^-_\mathbb Z(\textbf{e}), \eqno {(56)}$$
 for all $j\in J$ satisfying $j_1\geq j \geq (t,1)$.  Then there are no false
 trails in $\mathscr K_t^{BZ}$ which trivialize at $j_1$.

  \end {thm}

 \textbf{Remark}.  Thus to show that there are no false trails in $\mathscr K^{BZ}_t$, it suffices to show that $(56)$ holds for all $j \in J$.

 \subsection{}\label{9.3}

 One expects the left hand side of $(56)$ to consist of only extremal elements in which case one may replace condition $(56)$ by
 $$\cup_{\textbf{e}\in \textbf{E}_{i_j}^j}Z(\textbf{e}) \subset \cup_{\textbf{e}\in \textbf{E}_{i_{j+1}}}Z^-(\textbf{e}), \eqno {(57)}$$
 which in any case is stronger.  In our examples it was this that we verified.

 % Notice that by verifying the hypothesis of Theorem \ref {9.2} one obtains from Corollary \ref {8.3} the remarkable fact that  $\mathscr K_t(s)=\sqcup_{e \in \textbf{E}_s}K_\mathbb Z(\textbf{e})$ is independent of $s \in I$.  In our examples this had been a particularly tedious property to verify by hand, though in our examples we had actually verified the stronger statement that $Z_t(s)$ is independent of $s \in I$ thereby giving rise to the common set $Z_t$ which is a giant $S$-set.

   %Actually in these cases it is relatively easy to direct construct $Z_t$ by the method outlined in \ref {8.4}.  Still this does not guarantee that equality holds in $(52)$, though again this can be verified in special cases.

  \section{Index of Notation }\label{10}

  Symbols appearing frequently are given below in the paragraph they are first defined.

  \

  $1.1$. $\mathfrak g, \mathfrak h, \pi, \pi^\vee, W, P, P^+, \mathfrak b, \mathfrak n$.

  $1.2$. $e_\alpha, V(\lambda), v_{w\lambda}$.

  $1.3$. $I,s_i,\varpi_i, J, w_j, B_J, B_J(\infty)$.

  $1.5$. $Z_t, z_t^1, \varepsilon^\star_t$.

  $1.6$. $V(-\varpi_t)$.

  $2.2$.  $w_0, \ell(\cdot), \Omega(V(-\varpi_t)), (s,k), \gamma_j^K, K_t^1, \mathscr K^{BZ}_t, v_j^K$.

  $2.3$.  $z^K,m_j,\delta_j^K,m_s^k,r_s^k$.

  $2.4$.  $F_s^k,r_s^0, X_t$.

  $3.2$.  $n,v_\textbf{k},\textbf{a},\textbf{b},\textbf{k},\textbf{l}$.

  $3.3$.  $N, \hat{N}$.

  $4.2$.  $\overline{v}_\textbf{k}, T_s(\textbf{e}),M_s(\textbf{e})$.

  $4.3$.  $T^-_s(\textbf{e})T^+_s(\textbf{e})$.

  $4.5$.  $K_{\ell \min}, K_{\ell \max}$.

  $4.7.1$.  $[i,j]$.

  $4.7.2$.  $K_{\min}$.

  $4.7.3$.  $K_{\max}$.

  $5.2$.  $<\cdot,\cdot>$.

  $5.3.3$.  $c_j^{\ell \min}(s)$.

  $5.3.4$.  $\overline{\textbf{l}}$.

  $5.4$.  $c_i'^K(s)$.

  $7.1$.  $\prec,K(\textbf{c}),Z(\textbf{c}),K_\mathbb Z(\textbf{c})$.

  $7.2.1$.  $M(B_J),\mathscr G^+,\mathscr G^-,\mathscr G(\textbf{c})$.

  $7.2.2$.  $\mathscr G^j(\textbf{c})$.

  $8.5$.  $E(\cdot)$.

  $9.1$.  $\mathscr F^j(\mathscr K^{BZ}_t)$.


\begin{thebibliography}{m}

 \bibitem {BZ} A. Berenstein and A. Zelevinsky, Tensor product multiplicities, canonical bases and totally positive varieties. Invent. Math. 143 (2001), no. 1, 77–-128.

 \bibitem {GK}   O. Gabber and V. G.  Kac, On defining relations of certain infinite-dimensional Lie algebras. Bull. Amer. Math. Soc. (N.S.) 5 (1981), no. 2, 185–-189.

\bibitem {GP} O. Gleizer and A. Postnikov,  Littlewood-Richardson coefficients via Yang-Baxter equation. Internat. Math. Res. Notices 2000, no. 14, 741–-774.

 \bibitem {J1} A. Joseph, Results and problems in enveloping algebras arising from quantum groups. Representation theory, dynamical systems, and asymptotic combinatorics, 87–-100, Amer. Math. Soc. Transl. Ser. 2, 217, Amer. Math. Soc., Providence, RI, 2006.

 \bibitem {J2} A. Joseph,  Quantum groups and their primitive ideals. Ergebnisse der Mathematik und ihrer Grenzgebiete (3) [Results in Mathematics and Related Areas (3)], 29. Springer-Verlag, Berlin, 1995.

 \bibitem {J3}  A.  Joseph, Consequences of the Littelmann path theory for the structure of the Kashiwara $B(\infty)$ crystal. Highlights in Lie algebraic methods, 25–-64, Progr. Math., 295, Birk\"auser/Springer, New York, 2012.

 \bibitem {J4}   A. Joseph, A Preparation Theorem for the Kashiwara $B(\infty)$ Crystal,  Selecta Mathematica (to appear).

 \bibitem {J5} A, Joseph, Convexity Properties of the Canonical $S$-graphs.

 \bibitem {JL} A. Joseph and P. Lamprou, A new interpretation of the Catalan numbers, arXiv:1512.00406.

 \bibitem {JZ} A, Joseph and S. Zelikson, Dual Kashiwara functions for the $B(\infty)$ crystal.

\bibitem{Ka1} M. Kashiwara, Global crystal bases of quantum groups.
      Duke Math. J.  69  (1993), no. 2, 455--485.

\bibitem {Ka2} M. Kashiwara, The crystal base and Littelmann's refined Demazure
character formula. Duke Math. J.  71  (1993),  no. 3, 839-858.

\bibitem{Lu1} G. Lusztig, Canonical bases arising from quantized
enveloping algebras. II. Common trends in mathematics and quantum
field theories (Kyoto, 1990). Progr. Theoret. Phys. Suppl.  No.
102 (1990), 175-201 (1991).

\end{thebibliography}
\end{document}